\documentclass[12pt]{article}       
\usepackage{amsthm}
\usepackage{amsmath,amssymb}
\usepackage[arrow,matrix]{xy}
\usepackage{color}
\usepackage{comment}
\usepackage{version}
\newenvironment{NB}{
\color{red}{\bf NB}. \footnotesize 
}{}
\excludeversion{NB}

\newtheorem{thm}{Theorem}[section]

\newtheorem{prop}[thm]{Proposition}
\newtheorem{defn}[thm]{Definition}
\newtheorem{lem}[thm]{Lemma}

\def\min{\mathop{\mathrm{min}}\nolimits}

\def\Hom{\mathop{\mathrm{Hom}}\nolimits}
\def\im{\mathop{\mathrm{im}}\nolimits}

\def\rk{\mathop{\mathrm{rk}}\nolimits}
\def\ch{\mathop{\mathrm{ch}}\nolimits}
\def\id{\mathop{\mathrm{id}}\nolimits}
\def\exp{\mathop{\mathrm{exp}}\nolimits}
\def\Spec{\mathop{\mathrm{Spec}}\nolimits}

\def\C{\mathbb{C}}
\def\R{\mathbb{R}}
\def\e{\varepsilon}
\def\mbi#1{\boldsymbol{#1}}
\newcommand{\Res}{\operatornamewithlimits{Res}}
\newcommand{\ba}{\Bar}
\newcommand{\wt}{\widetilde}
\newcommand{\wh}{\widehat}
\newcommand{\mk}{\mathfrak}
\newcommand{\mc}{\mathcal}
\newcommand{\mb}{\mathbb}
\newcommand{\mo}{\mathcal{O}}
\newcommand{\E}{\mathcal{E}}
\newcommand{\F}{\mathcal{F}}
\newcommand{\PP}{\mathbb{P}}
\newcommand{\Z}{\mathbb{Z}}
\newcommand{\GL}{\operatorname{GL}}
\begin{document}

\title{Wall-crossing between stable and co-stable ADHM data%\thanks{This work was supported by JSPS Grant-in-Aid for JSPS Fellows.}
}
\author{Ryo Ohkawa
}

%\institute{
%\at Department of Mathematics, School of Fundamental Science and Engineering, Waseda University, 3--4--1 Okubo, Shinjuku-ku, Tokyo 169--8555, Japan \\
%\email{ohkawa.ryo@aoni.waseda.jp}\\
%Tel.: +81-3-5286-3195
%}

\maketitle

\begin{abstract}
We prove formula between Nekrasov partition functions defined from stable and co-stable ADHM data for the plane following method by Nakajima-Yoshioka \cite{NY3} based on the theory of wall-crossing formula developed by Mochizuki \cite{M}.
This formula is similar to conjectures by Ito-Maruyoshi-Okuda \cite[(4.1), (4.2)]{IMO} for $A_{1}$ singularity.
%\keywords{Nekrasov partition functions \and Framed sheaves}
%\subclass{14D21 \and 57R57}
%14D20  	Algebraic moduli problems, moduli of vector bundles 
%14D21  	Applications of vector bundles and moduli spaces in mathematical physics (twistor theory, instantons, quantum field theory)
%57R57  	Applications of global analysis to structures on manifolds, Donaldson and Seiberg-Witten invariants 
%81T13  	Yang-Mills and other gauge theories
%81T60  	Supersymmetric field theories
\end{abstract}

%%%%%%%%%%%%%%%%%%%%%%%%%%%%%%%%%%%%%%%%%%%%%%%%%%

\section{Intorduction}
\label{sec:intr}
Nekrasov partition functions are introduced in \cite{Nek}.
They are defined by integrations
$$
Z=\sum_{n=0}^{\infty} q^{n} \int_{M(r,n)} \psi
$$ 
on moduli spaces $M(r,n)$ of framed sheaves on the plane $\PP^{2}$ with the rank $r$ and the second Chern class $n$.
Here $\psi$ are various equivariant cohomology classes on $M(r,n)$ corresponding to physical theories.
These integrations are defined by localization for torus actions on moduli spaces (cf. \S \ref{subsec:part}). 
In particular we consider $T^{2}\times T^{r} \times T^{N_{f}}$-actions on $M(r,n)$ with $0 \le N_{f} \le 2r$ and $T=\C^{\ast}$.
Here $T^{2}$-actions are induced by the diagonal action on $\C^2 \subset \PP^2$, $T^{r}$-actions are induced by scale change of framings, and $T^{N_{f}}$ acts on certain vector bundles on $M(r,n)$ called {\it tautological bundles}.

Nekrasov's conjecture states that these partition functions give deformations of the Seiberg-Witten prepotentials for $N=2$ SUSY Yang-Mills theory.
This conjecture is proven in 
\cite{BE}, \cite{NO} and 
\cite{NY1} 
independently.
In \cite{NY1} they study the case 
for $\psi=1$.
%In \cite{NY1} they consider the case where $N_{f}=0$ and $\psi =1$ corresponding to $\mc N=2$ $U(r)$ pure gauge theory.
Furthermore in \cite{GNY} they study the case where $\psi$ is defined by the equivariant Euler class
\begin{eqnarray}
\label{psi}
\psi = e\left( \bigoplus_{f=1}^{N_{f}} \mc V \otimes \left( \frac{e^{m_{f}}}{\sqrt{t_{1}t_{2}}} \right) \right),
\end{eqnarray}
for $(t_{1},t_{2})\in T^{2}, (e^{m_{1}}, \ldots, e^{m_{N_{f}}}) \in T^{N_{f}}$, in particular with $N_{f}=1$ and $r=2$,
%corresponding to $\mc N=2$ $U(2)$ gauge theory with $N_{f}$ matters, 
where $\mc V$ are tautological bundles on $M(r,n)$.
They extend arguments in \cite{NY1} using the theory of perverse coherent sheaves \cite{NY3}.
Combining with Mochizuki's formula \cite{M} they proved the Witten's conjecture \cite{W} relating Donaldson invariants with Seiberg-Witten invariants for complex projective surfaces. 

On the other hand in \cite{IMO}, similar functions are considered on ALE spaces of type $A_{p-1}$ and quotient stacks $[\C^2/\Z_{p}]$ for $p>1$.
They conjecture formulas \cite[(4.1), (4.2)]{IMO} among these functions.
Physical background of Nekrasov partition functions is to compute integrations on moduli spaces of instantons on ALE spaces.
Since moduli of instantons are singular, they use moduli of framed sheaves as resolutions.

In the setting of \cite{IMO} we can consider two resolutions, moduli spaces of framed sheaves on ALE spaces of type $A_{p-1}$ and quotient stacks $[\C^2/\Z_{p}]$.
It is natural to expect that the difference of Nekrasov partition functions defined from two resolutions is small, and to ask how results depend on choices of resolutions.
%This is a motivation of the above conjecture.
These two resolutions are also considered as moduli of stable ADHM data ( or, quiver varieties in more general ) corresponding to different stability parameters.

In this paper we treat the case where $p=1$, and hence both the ALE space and the quotient stack coincide with $\C^2$.
But we still have two stability conditions and corresponding two resolutions from the viewpoint of the ADHM description.
These are moduli of stable ADHM data and co-stable ADHM data, which are isomorphic as manifolds, but having different torus actions.  

We coompare Nekrasov partition functions defined from stable and co-stable ADHM data.
Here we take $\psi$ in \eqref{psi} with $N_{f}=2r$. 
This gives a formula in Theorem \ref{main} similar to the above conjecture by \cite{IMO}. 
This formula could be predicted from AGT correspondence \cite{AGT} (cf. \cite{FL}, \cite{Neg}), but the author does not understand very well.
We also note that for $A_{1}$-singularity the similar phenomenon is studied in \cite{BBFLT} from different point of view.

If we take $\mbi \e=(\e_{1},\e_{2})$, $\mbi a=(a_{1}, \ldots, a_{r})$ and $\mbi m^{N_{f}}=(m_{1}, \ldots, m_{N_{f}})$ corresponding to characters of $T^{2}, T^{r}$ and $T^{N_{f}}$, they form a polynomial ring $\Z[\mbi \e, \mbi a, \mbi m^{N_{f}}]$ isomorphic to the $T^{2}\times T^{r}\times T^{N_{f}}$-equivariant Chow ring of a point.
Nekrasov partition functions $Z=Z(\mbi \e, \mbi a, \mbi m^{N_{f}}, q)$ take values in the quotient field $\mb Q(\mbi \e, \mbi a, \mbi m^{N_{f}})$, and ones defined from integrations on moduli of co-stable ADHM data are equal to $Z(-\mbi \e, \mbi a, \mbi m^{N_{f}},q)$.

As an application of Theorem \ref{main}, we determine the odd degree parts of 
$$
\e_{1} \e_{2}\log Z(\mbi \e, \mbi a, \mbi m^{N_{f}})
$$ 
for $0 \le N_{f} \le 2r$ in \S \ref{subsec:appl}.  
These are equal to zero for $N_{f} \le 2r-2$, and have non-zero coefficients only in $\e_{1}+\e_{2}$ for $N_{f} = 2r-1, 2r$.
It was known that only coefficients for $\e_{1}+\e_{2}$ are determined in the case where $N_{f} \le 2r-1$ ( cf. \cite[Lemma 7.1]{NY1}), and this method can not be applied to the case where $N_{f}=2r$.

Our proof follows the method in \cite{NY3} based on the theory of wall-crossing formula developed in \cite{M}.
We apply the method to the case where skyscraper sheaves destabilize framed sheaves on the wall.
This case is not treated before.
We will show in \cite{O} that similar wall-crossing formulas solve the above conjecture by \cite{IMO}.

The paper is organized as follows.
In \S \ref{sec:nota}, we list our notation and give outline of the paper.
In \S \ref{sec:fram}, we recall ADHM description of framed sheaves, give a statement of main result in Theorem \ref{main}, and an application for Nekrasov partition functions.
In \S \ref{sec:redu}, we recall quiver description of moduli of framed sheaves, and reduce a proof of Theorem \ref{main} to wall-crossing formulas.
In \S \ref{sec:enha}, we introduce enhanced master spaces used in Mochizuki method \cite{M}, \cite{NY3}.
In \S \ref{sec:obst}, we study obstruction theories for moduli stacks.
In \S \ref{sec:hilb}, we describe destabilizing objects in terms of Hilbert schemes.
In \S \ref{sec:wall}, we compute wall-crossing formulas, and complete a proof of Theorem \ref{main}.
%In Appendix A we compute integrations on Hilbert schemes following \cite{Na2}.

The author thanks Hiraku Nakajima for telling him the conjecture \cite[(4.1), (4.2)]{IMO} and advices, which leads him to the similar formula in the setting on the plane, and \cite{Na2} for computations in \S \ref{subsec:inte}, and many other advices. 
The author would like to appreciate referees for comments about AGT correspondence, and work in \cite{BBFLT}.
This work was supported by JSPS Grant-in-Aid for JSPS Fellows.

%%%%%%%%%%%%%%%%%%%%%%%%%%%%%%%%%%%%%%%%%%%%%%%%%%%%%%%%%%%%%%%%%%%%%%%%
%%%%%%%%%%%%%%%%%%%%%%%%%%%%%%%%%%%%%%%%%%%%%%%%%%%%%%%%%%%%%%%%%%%%%%%%
%%%%%%%%%%%%%%%%%%%%%%%%%%%%%%%%%%%%%%%%%%%%%%%%%%%%%%%%%%%%%%%%%%%%%%%%

\section{Notation and outline of the paper}
\label{sec:nota}

%%%%%%%%%%%%%%%%%%%%%%%%%%%%%%%%%%%%%%%%%%%%%%%%%%%%%%%%%%%%%%%%%%%%%%%%

\subsection{Moduli of ADHM data}
\label{subsec:modu0}
Let $[z_{0}, z_{1}, z_{2}]$ be the homogeneous coordinate of $\PP^{2}=\PP(\C \oplus Q)$, where $Q=\C^{2}$, and $\ell_{\infty} = \lbrace z_{0}=0 \rbrace$ the line at infinity.
In \S \ref{subsec:fram}, we consider {\it framed sheaves} on $\PP^{2}$, and describe them in terms of {\it ADHM data} as follows. 

Let $W=\C^{r}$ and $V=\C^{n}$ be vector spaces, where $r$ and $n$ corresponds to ranks and second Chern classes of framed sheaves later (cf. Theorem \ref{barth}).
We put
\begin{eqnarray*}
\mb M(W, V)
&=&
\Hom_{\C}(Q^{\vee} \otimes V, V) \oplus \Hom_{\C}(W,V)\oplus \Hom_{\C}(\wedge Q^{\vee} \otimes V,W),\\
\mb L(V)
&=&
\Hom_{\C}(\wedge Q^{\vee} \otimes V, V).
\end{eqnarray*}
We also write $\mb M=\mb M(r,n)=\mb M(W,V)$ and $\mb L=\mb L(n)=\mb L(V)$.
We consider a map ( cf. Definition \ref{adhm})
\begin{eqnarray}
\label{mu}
\mu \colon \mb M \to \mb L, X=(B, z,w) \mapsto \mu(X)=[B \wedge B]+zw,
\end{eqnarray}
where $B \in \Hom_{\C}(Q^{\vee} \otimes V, V)$, $z \in \Hom_{\C}(W,V)$, and $w \in \Hom_{\C}(\wedge Q^{\vee} \otimes V,W)$.

For elements in $\mb M(W, V)$, we introduce stability and co-stability conditions.
In particular, elements $X=(B, z,w)$ in $\mu^{-1}(0)$ are called ADHM data on $(W, V)$, and we put 
\begin{eqnarray*}
M(r, n)
&=&
\lbrace (B, z,w) \mid \text{ stable ADHM data on }(W, V)\rbrace / \GL(V),\\
M^{c}(r, n)
&=&
\lbrace (B, z,w) \mid \text{ co-stable ADHM data on }(W, V) \rbrace / \GL(V).
\end{eqnarray*}
Our goal, Theorem \ref{main}, is to compare integrations over $M(r,n)$ and $M^{c}(r,n)$ (cf. \S \ref{sec:redu} ).
To this end, in \S \ref{sec:enha} we introduce further moduli spaces as follows.

%%%%%%%%%%%%%%%%%%%%%%%%%%%%%%%%%%%%%%%%%%%%%%%%%%%%%%%%%%%%%%%%%%%%%%%%

\subsection{Enhanced master space}
\label{subsec:enha}

We put $\underbar{n}=\lbrace 1, \ldots, n\rbrace$, and write by $Fl=Fl(V, \underbar{n})$ the full flag variety of $V=\C^{n}$.
We put 
$$
\wt{\mb M}=\wt{\mb M}(r,n)=\wt{\mb M}(W, V)=\mb M(W, V) \times Fl, 
$$
and introduce $\GL(V)$-equivariant line bundles $L_{-}$ and $L_{+}$ on $\wt{\mb M}$, which define stable loci as complements of all zero loci of $\GL(V)$-equivariant sections.
We put 
$$
\wh{\mb M}=\wh{\mb M}(r,n)=\wh{\mb M}(W,V)=\mb P_{\wt{\mb M}}(L_{-} \oplus L_{+}).
$$
We also consider maps $\hat{\mu} \colon \wh{\mb M} \to \mb L$ and $\tilde{\mu} \colon \wt{\mb M} \to \mb L$ defined by  compositions of projections to $\mb M$ and $\mu$.
Objects $(X, F^{\bullet})$ in $\tilde{\mu}^{-1}(0)$ are ADHM data with full flags of $V$, and the quotient stack
$$
\mc M = [\hat{\mu}^{-1}(0)^{ss} / \GL(V) ]
$$ 
is called the {\it enhanced master space}, where $\hat{\mu}^{-1}(0)^{ss}$ is a semi-stable locus with respect to the tautological line bundle $\mo(1)$ of the projective bundle $\wh{\mb M}=\mb P_{\wt{\mb M}}(L_{-} \oplus L_{+})$.

%%%%%%%%%%%%%%%%%%%%%%%%%%%%%%%%%%%%%%%%%%%%%%%%%%%%%%%%%%%%%%%%%%%%%%%%

\subsection{$\C^{\ast}_{\hbar}$-fixed points set}
\label{subsec:cast}

We consider an algebraic torus $\C^{\ast}_{\hbar}$, and $\C^{\ast}_{\hbar}$-action on $\mc M$ induced by a map $[x_{-}, x_{+}] \mapsto [ e^{\hbar} x_{-}, x_{+}]$ for $e^{\hbar} \in \C^{\ast}_{\hbar}$, where $[x_{-}, x_{+}]$ is the homogeneous coordinate of fibers of the projective bundle $\wh{\mb M}=\PP(L_{-} \oplus L_{+})$.
We have $\C^{\ast}_{\hbar}$-fixed point set $\mc M^{\C^{\ast}_{\hbar}}$. 
For each $\ell \in \underbar{n} = \lbrace 1, \ldots, n\rbrace$, if we choose $L_{\pm}$ suitably, then we get a decomposition 
$$
\mc M^{\C^{\ast}_{\hbar}}= \mc M_{+} \sqcup \mc M_{-} \sqcup \bigsqcup_{\mk J \in S^{\ell}} \mc M_{\mk J},
$$
where $\mc M_{\pm}$ is obvious components defined by a zero locus of $x_{\mp}$, and $S^{\ell}$ is a finite set of decomposition data $(I_{\flat}, I_{\sharp})$.
Here $I_{\flat}$ and $I_{\sharp}$ are subsets of $\underbar{n}=\lbrace 1, \ldots, n\rbrace$ such that $\underbar{n}=I_{\flat} \sqcup I_{\sharp}$, $I_{\sharp} \neq \emptyset$, and $\min (I_{\sharp}) \le \ell$.

For each decomposition data $\mk J=(I_{\flat}, I_{\sharp}) \in S^{\ell}$, we fix a decomposition $V=V_{\flat} \oplus V_{\sharp}$ with $\dim V_{\flat} = |I_{\flat}|$ and $\dim V_{\sharp}=|I_{\sharp}|$, and put $p = \dim V_{\sharp}$.
Roughly $\mc M_{\mk J}$ parametrizes $\GL(V_{\flat}) \times \GL(V_{\sharp})$-orbits of objects 
$$
(X_{\flat}, F_{\flat}^{\bullet}, X_{\sharp}, F_{\sharp}^{\bullet}, \rho) .
$$
Here $(X_{\flat}, F_{\flat}^{\bullet})$ are ADHM data on $(W, V_{\flat})$ with full flags of $V_{\flat}$, $(X_{\sharp}, F_{\sharp}^{\bullet})$ are ADHM data on $(0, V_{\sharp})$ with full flags of $V_{\sharp}$, and $\rho= x_{-} / x_{+}$ are called orientations of $V$ since $g \in \GL(V)$ acts by $\det g^{-D} \rho$ for certain integer $D$.

By decomposition data $\mk J$, we allow the indices of $F_{\flat}^{\bullet}$ and $F_{\sharp}^{\bullet}$ repetitions such that
$$
I_{\flat} = \lbrace i \in \underbar{n} \mid F_{\flat}^{i} / F_{\flat}^{i-1} \neq 0 \rbrace, I_{\sharp} = \lbrace i \in \underbar{n} \mid F_{\sharp}^{i} / F_{\sharp}^{i-1} \neq 0 \rbrace.
$$
Then $F_{\flat}^{\bullet} \oplus F_{\sharp}^{\bullet}$ can be regarded as full flags of $V$, denoted by $Fl_{\mk J}(F_{\flat}^{\bullet}, F_{\sharp}^{\bullet})$ in \S \ref{subsec:deco}.
Furthermore $(X_{\flat}, F_{\flat}^{\bullet})$ have only trivial automorphism, but $(X_{\sharp}, F_{\sharp}^{\bullet})$ have automorphism groups $\C^{\ast} \id_{V_{\sharp}}$.
Hence $\GL(V)$-orbits of $(X_{\flat} \oplus X_{\sharp}, F_{\flat}^{\bullet} \oplus F_{\sharp}^{\bullet}, \rho) \in \hat{\mu}^{-1}(0)^{ss}$ contains a fiberwise direction
$$
(X_{\flat} \oplus X_{\sharp}, F_{\flat}^{\bullet} \oplus F_{\sharp}^{\bullet}, t^{-pD} \rho) \in \hat{\mu}^{-1}(0)^{ss}
$$ 
for $t \id_{V_{\sharp}} \in \GL(V_{\sharp})$, and represents a fixed point in $\mc M^{\C^{\ast}_{\hbar}}$.
In \S \ref{subsec:dire}, it will be shown that all $\C^{\ast}_{\hbar}$-fixed points outside $\mc M_{\pm}$ are obtained in this way under suitable stability conditions.

%%%%%%%%%%%%%%%%%%%%%%%%%%%%%%%%%%%%%%%%%%%%%%%%%%%%%%%%%%%%%%%%%%%%%%%%

\subsection{Localizations}
\label{subsec:loca}

In \S \ref{sec:obst}, obstruction theories of these spaces $\mc M, \mc M_{\pm}$, and $\mc M_{\mk J}$ are naturally introduced from their constructions. 
From these obstruction theories, we have virtual fundamental cycles $[\mc M]^{vir}, [\mc M_{\pm}]^{vir}$, and $[\mc M_{\mk J}]^{vir}$, and normal bundles $\mk N (\mc M_{\alpha})$ of $\mc M_{\alpha}$ in $\mc M$ for $\alpha=\pm, \mk J \in S^{\ell}$.
By the main result \cite[(1)]{GP}, we have 
\begin{eqnarray}
\label{local0}
[\mc M]^{vir} = \iota_{\ast} \sum_{\alpha=\pm, \mk J} \frac{[\mc M_{\alpha}]^{vir}}{e(\mk N(\mc M_{\alpha}))} \in A^{\C^{\ast}_{\hbar}}_{\ast}( \mc M ) \otimes \mb Q[\hbar, \hbar^{-1}], 
\end{eqnarray}
where $\iota \colon \mc M^{\C^{\ast}_{\hbar}} \to \mc M$ is the inclusion, and $A^{\C^{\ast}_{\hbar}}_{\ast}( \mc M )$ is the $\C^{\ast}_{\hbar}$-equivariant Chow ring of $\mc M$, and $\hbar$ is the first Chern class of the $\C^{\ast}_{\hbar}$-weight space with eigenvalue $e^{\hbar}$.
We note that we need torus actions other than $\C^{\ast}_{\hbar}$ to define integrations over these moduli spaces, and hence we must modify \eqref{local0} later in \eqref{local}.

The left hand side of \eqref{local0} is a polynomial in $\hbar$, but right hand side contains $\hbar^{-1}$.
Hence taking coefficients of $\hbar^{-1}$, the difference between integrations over $\mc M_{\pm}$ is reduced to integrations over $\mc M_{\mk J}$. 
Furthermore in our choice of $L_{\pm}$, $\mc M_{-}$ is a full flag bundle over $M(r,n)$ while $\mc M_{+}$ parametrizes $\ell$-stable objects, where $\ell$-stability for $(X, F^{\bullet})$ coincides with stability and co-stability for $X$ when $\ell=0$ and $\ell=n$ respectively (cf. Definition \ref{ellstab}).
Moduli of $\ell$-stable objects $(X, F^{\bullet})$ on $(W, V)$ are denoted by $\wt{M}^{\ell}(r,n)$, where $X$ are ADHM data on $(W,V)$ and $F^{\bullet}$ are full flags of $V$.

%%%%%%%%%%%%%%%%%%%%%%%%%%%%%%%%%%%%%%%%%%%%%%%%%%%%%%%%%%%%%%%%%%%%%%%%

\subsection{Further decompositions}
\label{subsec:furt}
We consider decompositions of $\mc M_{\mk J}$ as follows.
Roughly these moduli stacks $\mc M_{\mk J}$ parametrize data
$$
(X_{\flat}, F_{\flat}^{\bullet}, X_{\sharp}, F_{\sharp}^{\bullet}, \rho)
$$
up to $\GL(V_{\flat}) \times \GL(V_{\sharp})$-action as in \S \ref{subsec:cast}.

We break orientations $\rho$ of $V=V_{\flat} \oplus V_{\sharp}$ into pairs $(\rho_{\flat}, \rho_{\sharp})$ of orientations of $V_{\flat}$ and $V_{\sharp}$, and add another torus $\C^{\ast}_{t}$ to symmetry groups, where $t$ acts by $(t \rho_{\flat}, t^{-1} \rho_{\sharp})$.
Furthermore, we introduce $s= t^{1/pD}$, and get \'etale cover $S_{\mk J}$. 
After automorphisms $(g_{\flat}, g_{\sharp}, s) \mapsto (g_{\flat}, s^{-1} g_{\sharp}, s)$ of $\GL(V_{\flat}) \times \GL(V_{\sharp}) \times \C^{\ast}_{s}$, the moduli stacks $S_{\mk J}$ are products of quotient stacks parametrizing $(X_{\flat}, F_{\flat}^{\bullet}, \rho_{\flat})$ and $(X_{\sharp},  F_{\sharp}^{\bullet}, \rho_{\sharp})$ up to $\GL(V_{\flat}) \times \C^{\ast}_{s}$-action and $\GL(V_{\sharp} )$-action.
The first components are \'etale covers of $\wt{M}^{\min(I_{\sharp})-1}(r, n-p)$.

The second components are denoted by $M_{p}^{+}$.
We see that $M_{p}^{+}$ are \'etale covers of full flag bundles over $M(1,p)$ as follows.
From data $(X_{\sharp}, F_{\sharp}^{\bullet})$, we take generators $z_{\sharp}$ of $F_{\sharp}^{\min(I_{\sharp})} \subset V_{\sharp}$, and regard flags $F_{\sharp}^{\bullet}$ as pairs $( z_{\sharp}, \bar{F}_{\sharp}^{\bullet})$, where $\bar{F}_{\sharp}^{\bullet}$ are full flags of $V_{\sharp} / \langle z_{\sharp} \rangle$ induced by $F_{\sharp}^{\bullet}$.
Then $M_{p}^{+}$ are isomorphic to quotient stacks of spaces parametrizing $(X_{\sharp}, z_{\sharp}, \bar{F}_{\sharp}^{\bullet}, \rho_{\sharp})$ by the groups $\GL(V_{\sharp}) \times \C^{\ast}_{u}$, where $u$ acts by $u z_{\sharp}$.
After automorphisms $(g_{\sharp} ,u) \mapsto (u^{-1} g_{\sharp} ,u)$ of $\GL(V_{\sharp}) \times \C^{\ast}_{u}$, pairs $(X_{\sharp}, z_{\sharp})$ can be regarded as ADHM data on $(W_{\sharp}, V_{\sharp})$. 
Here we introduce a vector space $W_{\sharp}=\C$, and it does not have $\GL(W)$-action.
We note that $\C^{\ast}_{u}$ also acts on orientations by the multiplications of $u^{pD}$, and $M_{p}^{+}$ are \'etale covers of full flag bundles over $M(1,p)$.

Hence integrations over $\mc M_{\mk J}$ are essentially reduced to integrations over $\wt{M}^{\min(I_{\sharp})-1}(r, n-p)$, and Hilbert schemes $M(1, p)$.
Hence we can continue computations recursively.
Since $\min(I_{\sharp}) \le \ell$, this process terminates finally, and we get Theorem \ref{main} (cf. \S \ref{sec:wall}).

%%%%%%%%%%%%%%%%%%%%%%%%%%%%%%%%%%%%%%%%%%%%%%%%%%%%%%%%%%%%%%%%%%%%%%%%

\subsection{Ambient stacks}
\label{subsec:ambi}

Technically, it is often convenient to consider ambient spaces of the above moduli spaces, since $\mc M$ is not smooth.
We put $\mc N = \wh{\mb M}^{ss} / \GL(V)$, where $\wh{\mb M}^{ss}$ is the semistable locus with respect to the tautological bundle $\mo(1)$.
We see that $\mc N$ is a smooth Deligne-Mumford stack, and we have $\C^{\ast}_{\hbar}$-action on $\mc N$.
By definition, $\C^{\ast}_{\hbar}$-fixed points set $\mc N^{\C^{\ast}_{\hbar}}$ is the zero locus of vector fields on $\mc N$ generated by the $\C^{\ast}_{\hbar}$-action, and $\mc M^{\C^{\ast}_{\hbar}}$ is defined by the pull-back of $\mc N^{\C^{\ast}_{\hbar}}$ by the inclusion $\mc M \to \mc N$ (cf. \cite[Appendix C]{GP}).
The above decomposition of $\mc M^{\C^{\ast}_{\hbar}}$ is induced by the similar deomposition in \S \ref{subsec:modu}
$$
\mc N^{\C^{\ast}_{\hbar}}= \mc N_{+} \sqcup \mc N_{-} \sqcup \bigsqcup_{\mk J \in S^{\ell}} \mc N_{\mk J}.
$$ 
Furthermore $\hat{\mu} \colon \wh{\mb M} \to \mb L, \GL(V), \mo(1)$ forms {\it Kuranishi chart} of $\mc M$ introduced in \S \ref{subsec:virt}.
This concept helps to compute obstruction theories.

We mainly concern about ADHM data, but from the above technical reason, we consider various stability conditions to data including elements in $\mb M(W,V)$ rather than ADHM data.
However for moduli stacks $M_{p}^{+}$, we consider ambient spaces $\mb M(W,V)/ \Hom_{\C}(\wedge Q^{\vee} \otimes V, W)$ from the construction in the previous section.
This causes integrations over virtual fundamental cycles $[M(1,p)]^{vir}$ different from usual fundamental cycles $[M(1,p)]$.

%%%%%%%%%%%%%%%%%%%%%%%%%%%%%%%%%%%%%%%%%%%%%%%%%%%%%%%%%%%%%%%%%%%%%%%%

%\subsection{Outline of proof}
%\label{subsec:outl}

%%%%%%%%%%%%%%%%%%%%%%%%%%%%%%%%%%%%%%%%%%%%%%%%%%%%%%%%%%%%%%%%%%%%%%%%
%%%%%%%%%%%%%%%%%%%%%%%%%%%%%%%%%%%%%%%%%%%%%%%%%%%%%%%%%%%%%%%%%%%%%%%%
%%%%%%%%%%%%%%%%%%%%%%%%%%%%%%%%%%%%%%%%%%%%%%%%%%%%%%%%%%%%%%%%%%%%%%%%

\section{Framed sheaves on $\PP^2$}
\label{sec:fram}
We recall ADHM descriptions of framed sheaves on $\PP^2$ from \cite{NY1}, and introduce partition functions and our main result Theorem \ref{main}.

%%%%%%%%%%%%%%%%%%%%%%%%%%%%%%%%%%%%%%%%%%%%%%%%%%%%%%%%%%%%%%%%%%%%%%%%

\subsection{Framed sheaves and ADHM data}
\label{subsec:fram}
We consider the projective plane $\PP^2$ over $\C$ and the line $\ell_{\infty}= \lbrace x_{0}=0 \rbrace$, where $[x_0, x_1, x_2]$ is the homogeneous coordinate of $\PP^2$.
\begin{defn}
A framed sheaf on $\PP^2$ is a pair $(E, \Phi)$ of
\begin{enumerate}
\item[-] a torsion free sheaf $E$ on $\PP^2$, and
\item[-] an isomorphism $\Phi \colon E|_{\ell_{\infty}} \cong \mo_{\ell_{\infty}}$.
\end{enumerate}
\end{defn}
We remark that framed sheaves are automatically slope semistable.
Moduli spaces of framed sheaves were constructed in \cite{HL} in more general framework.
In this paper, we construct the moduli spaces via ADHM description.
To introduce ADHM description of framed sheaves, we take finite dimensional vector spaces $Q=\C^2, W=\C^{r}$ and $V=\C^{n}$, and put
$$
\mb M(W, V)
=
\Hom_{\C}(Q^{\vee} \otimes V, V) \oplus \Hom_{\C}(W,V)\oplus \Hom_{\C}(\wedge Q^{\vee} \otimes V,W).
$$\begin{defn}
\label{adhm}
ADHM data on vector spaces $(W, V)$ are collections of linear maps $(B, z,w)$ such that $B \in \Hom_{\C}(Q^{\vee}\otimes V, V), z \in \Hom_{\C}(W, V)$ and $w \in \Hom_{\C}(\wedge Q^{\vee}\otimes V,W)$ satisfying
$$
[B \wedge B] + zw =0 \in \Hom_{\C}(\wedge Q^{\vee} \otimes V, V),
$$ 
where $[B \wedge B]$ is the restriction of $B \circ ( \id_{Q^{\vee}} \otimes B ) \colon Q^{\vee} \otimes Q^{\vee} \otimes V \to V$ to the subspace $\wedge Q^{\vee} \otimes V$ of $Q^{\vee} \otimes Q^{\vee} \otimes V$.
\end{defn}
If we write $Q\otimes V = V \oplus V$, and $B=(B_{1}, B_{2})$ by a linear map $B_{i} \colon V \to V$ for $i=1,2$, then the equation $[B \wedge B] +zw =0$ is equivalent to $[B_1,B_2] + zw =0$.
 
\begin{defn}
\label{stab}
An element $(B, z, w)$ in $\mb M(W,V)$ is said to be stable if there exists no subspace $S \subset V$ other than $S=V$ such that $B_i(S)\subset S$ for $i=1,2$ and $\im z \subset S$.  
\end{defn}

We put 
$$
M(r, n)=\lbrace (B, z,w) \mid \text{ stable ADHM data on }(W, V)\rbrace / \GL(V),
$$
where $\GL(V)$ acts by base change of $V$, that is, 
$$
g(B_1, B_2, z, w) = (gB_1g^{-1}, gB_2g^{-1}, g z, w g^{-1}).
$$
We have natural $\GL(Q) \times \GL(W)$-action on $M(r,n)$.
In particular $M(r,n)$ has $T^{2} \times T^{r}$-equivariant structure, where $T^2$ and $T^r$ are diagonal tori of $\GL(Q)$ and $\GL(W)$ respectively. 

\begin{thm}
[\protect{\cite{B}}]
\label{barth}
We have an isomorphism from $M(r,n)$ to the moduli of isomorphism classes of framed sheaves $(E, \Phi)$ with $\rk(E) =r=\dim W, c_{2}(E)=n=\dim V$.
\end{thm}
A framed sheaf $(E, \Phi)$ corresponding to $(B,z,w) \in M(r,n)$ via the above isomorphism is defined by the complex
$$
C^{\bullet} = \left(\mo_{\PP^2}(-1) \otimes \wedge Q^{\vee} \otimes V \stackrel{\sigma}{\to} \mo_{\PP^2}\otimes Q^{\vee} \otimes V \oplus \mo_{\PP^2} \otimes W
 \stackrel{\tau}{\to} \mo_{\PP^2}(1) \otimes V\right).
 $$
 Here 
 $\sigma= 
\begin{bmatrix}
B_{1}x_0 - x_{1} \id_{V}\\
B_{2}x_0- x_{2} \id_{V}\\
wx_{0}
\end{bmatrix}$, 
$\tau= 
\begin{bmatrix}
-B_{2}x_0 + x_{2} \id_{V} & B_{1}x_0 - x_{1} \id_{V} & zx_{0}
\end{bmatrix}$,
and $\id_{V}$ is the identity map of $V$.
The condition $[B_1, B_2]+ zw =0$ implies $\tau \sigma =0$.
We define $E= \ker \tau /\im \sigma$. 
Then substituting $z_0=0$, we have a natural isomorphism $\Phi \colon E|_{\ell_{\infty}} \cong W \otimes \mo_{\ell_{\infty}}$.
This gives an isomorphism (see \cite[Theorem 2.1]{Na1} for the proof).
Hereafter via this isomorphism we identify $M(r,n)$ and the moduli of isomorphism classes of framed sheaves $(E, \Phi)$ with $\rk(E) =r, c_{2}(E)=n$.

%%%%%%%%%%%%%%%%%%%%%%%%%%%%%%%%%%%%%%%%%%%%%%%%%%%%%%%%%%%%%%%%%%%%%%%%

\subsection{Torus action on $M(r,n)$}
\label{subsec:toru}
To describe torus fixed points of $M(r,n)$ we give a sheaf description of the torus action on $M(r,n)$ introduced in the previous subsection.
For $T=\C^{\ast}$, we put $\tilde{T} = T^2 \times T^r \times T^{2r}$ and 
$$
\mbi t= (t_1,t_2)\in T^{2}, e^{\mbi a}= (e^{a_1}, \ldots, e^{a_r})\in T^{r}, e^{\mbi m} = (e^{m_1}, \ldots, e^{m_{2r}}) \in T^{2r},
$$
where $\mbi a=(a_{1}, \ldots, a_{r}), \mbi m= (m_{1}, \ldots, m_{2r})$.
For $\mbi t=(t_1, t_2)\in T^2$, we consider the morphism $F_{\mbi t} \colon \PP^2 \to \PP^2$ defined by $[x_0,x_1,x_2] \mapsto [x_0, t_1x_1, t_1 x_2]$.
We identify $e^{\mbi a}=(e^{a_1}, \ldots, e^{a_r}) \in T^r$ with the diagonal matrix $\text{diag}(e^{a_1}, \ldots, e^{a_r})$.
We define $(\mbi t, e^{\mbi a}) (E, \Phi) = ( E', \Phi')$ by $E'=(F_{\mbi t}^{-1})^{\ast}E$, and $\Phi'$ is defined by the following commutative diagram:
$$
\xymatrix{
\ar@{=}[d] E' |_{\ell_{\infty}} \ar[rr]^{\Phi'}& & \mo_{\ell_{\infty}}^{\oplus r}\\
(F_t^{-1})^{\ast} E |_{\ell_{\infty}} \ar[r]^{(F_{t}^{-1})^{\ast}\Phi}& (F_t^{-1})^{\ast} \mo_{\ell_{\infty}}^{\oplus r} \ar@{=}[r] & \mo_{\ell_{\infty}}^{\oplus r} \ar[u]_{e^{\mbi a}}
}
$$
Let $T^{2r}$ act on $M(r,n)$ trivially.
These actions on $M(r,n)$ are compatible with ones defined by ADHM data in the previous section via the isomorphism in Theorem \ref{barth}.
 
\begin{prop}
[\protect{\cite[Proposition 2.9]{NY1}}]
\label{ny}
The $\tilde{T}$-fixed points of $M(r,n)$ are given by 
$$
I_{\vec{Y}}=(I_{Y_{1}} \oplus \cdots \oplus I_{Y_{r}}, \Phi), 
$$
where $\vec{Y}=(Y_{1}, \ldots, Y_{r})$ is a $r$-tuple of Young diagrams with $\sum_{i=1}^{r} |Y_{i}| =n$,  $I_{Y_{1}}, \ldots, I_{Y_{r}}$ are corresponding monomial ideals supported at $[1,0,0] \in \PP^2$, and $\Phi$ is a direct sum of natural isomorphisms $I_{Y_{i}}|_{\ell_{\infty}} \cong \mo_{\ell_{\infty}}$ induced by the inclusion $I_{Y_{i}} \subset \mo_{\PP^2}$.
\end{prop}

%We consider some $\tilde{T}$-equivariant sheaves on $M(r,n)$.
By Theorem \ref{barth} we have a universal framed sheaf $(\E, \tilde \Phi)$ on $\PP^2 \times M(r,n)$, that is, for each point $(E, \Phi) \in M(r,n)$  we have a unique isomorphism
$$
\E|_{\PP^2 \times \lbrace (E, \Phi)\rbrace } \cong E
$$
such that $\tilde \Phi$ and $\Phi$ commute on $\ell_{\infty}$.
This unique isomorphism gives a $\tilde T$-equivariant structure of $\E$.
\begin{NB}
For any point $(P, E, \Phi) \in \PP^2 \times M(r,n)$ and $(t,a) \in T^2 \times T^{r}$ we have identities
\begin{eqnarray*}
t^{\ast} \E |_{(P, E, \Phi)} 
&=&
\E|_{(t,a)(P,E, \Phi)}\\
&=&
\E|_{(F_t(P),(F_t^{-1})^{\ast} E, \Phi')}\\
&=&
(F_t^{-1})^{\ast} E |_{F_t(P)}
\end{eqnarray*}
giving $(t,a)^{\ast}\E \to \E$.
For $m\in T^{2r}$ we have $m^{\ast}\E = \E$ completing the definition of $\tilde T$-structure on $\E$, since $T^{2r}$ acts on $M(r,n)$ trivially.
\end{NB}

We consider the {\it tautological bundle} $\mc V = \mc V_{r,n} = (\mb R^1 p_{2} )_{\ast} (\mc E \otimes p_{1}^{\ast} \mo_{\PP^2}(-1)) $, where $p_{1}\colon \PP^2 \times M(r,n) \to \PP^2, p_{2} \colon \PP^2 \times M(r,n) \to M(r,n)$ are projections.
Then $\mc V$ is a vector bundle with the fiber $H^1(\PP^2, E(-1))$ over $(E, \Phi)$.
The $\tilde{T}$-action on $\E$ induces a $\tilde T$-equivariant structure of $\mc V$.
%by the identities
%$$
%\mc V|_{(E, \Phi)} = H^1(\PP^2, E(-1)) = H^1(\PP^2, (t, e, m)^{\ast} E) = \mc V|_{(t,e,m)(E, \Phi)}.
%$$
We also consider the tangent bundle $TM(r,n)$ with the natural $\tilde T$-action.

%%%%%%%%%%%%%%%%%%%%%%%%%%%%%%%%%%%%%%%%%%%%%%%%%%%%%%%%%%%%%%%%%%%%%%%%

\subsection{Main result and Nekrasov Partition functions}
\label{subsec:part}
\begin{NB}
For $i=(i_1, i_2) \in \Z^2, j=(j_1, \ldots, j_r) \in \Z^r, k=(k_1, \ldots, k_{2r}) \in \Z^{2r}$, we write $t^i e^j m^k$ to denote the one dimensional $\tilde T$-representation with eigen-value $t_1^{i_1}t_2^{i_2} e_1^{j_1} \cdots e_r^{j_r} m_1^{k_1} \cdots m_{2r}^{k_{2r}}$. 
For $r$-tuples of Young diagrams $\vec{Y}=(Y_1, \ldots, Y_r)$ we take $I_{\vec{Y}}=(I_{Y_{1}} \oplus \cdots \oplus I_{Y_{r}}, \Phi )$ the $\tilde T$-fixed point as in Proposition \ref{ny}.
For any $\tilde T$-equivariant vector bundle $\mc F=\mc F_{r,n}$ on $M(r,n)$ the fiber $\mc F |_{I_{\vec{Y}}}$ has the induced $\tilde T$-module structure, and we have the eigen space decomposition
$$
\bigoplus_{(i,j,k)\in \Z^{3r+2}} (t^i e^j m^k)^{\oplus a_{i,j,k}}.
$$
Then we put 
$$
e_{\vec{Y}}(\mc F)= \prod_{(i,j,k)\in \Z^{3r+2}} (i_1 \e_1 + i_2 \e_2 + j_1 a_1 + \cdots + j_ra_r + k_1 \mu_1 + \cdots + k_{2r} \mu_{2r}  )^{a_{i,j,k}} \in \C(\e, a, \mu),
$$
where $\e=(\e_1, \e_2), a= (a_1, \ldots, a_r)$ and $\mu =(\mu_1, \ldots, \mu_{2r})$.
\end{NB}

Let $A_{\ast}^{\tilde T}(X)$ be the $\tilde T$-equivariant Chow group of a $\tilde T$-space $X$ with rational coefficients.
They are modules over the $\tilde{T}$-equivariant Chow ring $A^{\ast}_{\tilde T}(\text{pt})$ of a point, isomorphic to $S(\tilde T)=\mb Z[\mbi \e, \mbi a, \mbi m]$, where $\mbi \e=(\e_{1},\e_{2}), \mbi a=(a_{1}, \ldots, a_{r})$ and $\mbi m=(m_{1}, \ldots, m_{2r})$ correspond to characters of $\tilde T$ with eigen-values $\mbi t, e^{\mbi a}, e^{\mbi m}$.
The quotient field of $S(\tilde T)$ is denoted by $\mc S$.

We have a projective morphism $\pi \colon M(r,n) \to M_{0}(r,n)$, where $M_{0}(r,n)$ is the Uhlenbeck (partial) compactification of the moduli space $M_{0}^{\text{reg}}(r,n)$ of framed {\it locally free} sheaves $(E,\Phi)$.
This morphism induces a homomorphism $\pi_{\ast}\colon A_{\ast}^{\tilde T}(M(r,n)) \to A_{\ast}^{\tilde T}(M_{0}(r,n))$.
In the next section we will explain the construction of $M_{0}(r,n)$ and $\pi$ in \cite{Na1} via the quiver description.  
Set-theoretically we have a bijection
$$
M_{0}(r,n) \to \bigsqcup_{n'=0}^{n} M_{r}^{\text{reg}}(r,n') \times S^{n-n'}(\C^2),
$$
and the $\tilde T$-fixed points set $M_{0}(r,n)^{\tilde T}$ consists of only one point $n[0]$.

By the localization theorem \cite[Theorem 1]{EG} we have an isomorphism 
$$
(\iota_{0})_{\ast} \colon \mc S \cong A_{\ast}^{\tilde T}(M_{0}(r,n)) \otimes \mc S,
$$ 
where $\iota_{0} \colon M_{0}(r,n)^{\tilde T}=\lbrace n[0] \rbrace \to M_{0}(r,n)$ is the inclusion.
For $\psi \in A_{\ast}^{\tilde T}(M(r,n))$, we put 
$$
\int_{M(r,n)} \psi = (\iota_{0})_{\ast}^{-1}\pi_{\ast} ( \psi \cap [M(r,n)]) \in \mc S.
$$
By Proposition \ref{ny} we have a commutative diagram
$$
\xymatrix{ 
\ar[d]_{\pi_{\ast}} A_{\ast}^{\tilde T}(M(r,n)) \otimes \mc S \ar[r]^{\hspace{1cm}\cong}_{\hspace{1cm}\oplus (\iota_{\vec{Y}})_{\ast}^{-1}}& \bigoplus_{\vec{Y}} \mc S \ar[d]^{\sum}\\
A_{\ast}^{\tilde T}(M_{0}(r,n)) \otimes \mc S \ar[r]^{\hspace{1cm}\cong}_{\hspace{1cm} (\iota_{0})_{\ast}^{-1}} & \mc S
}
$$
where $\iota_{\vec{Y}} \colon \lbrace I_{\vec{Y}} \rbrace \to M(r,n)$ is the inclusion, and $(\iota_{\vec{Y}})_{\ast}^{-1}(\psi)=\frac{\iota_{\vec{Y}}^{\ast}\psi}{  e(TM(r,n))}$.
Hence we have 
$$
\int_{M(r,n)} \psi = \sum_{\vec{Y}} \frac{\iota_{\vec{Y}}^{\ast} \psi }{ e(TM(r,n))} \in \mc S.
$$

We consider the equivariant Euler class $e(\mc F_{r}(\mc V))$ in $\mc S$ of a $\tilde{T}$-equivariant vector bundle 
\begin{eqnarray}
\label{fr}
\mc F_{r}(\mc V)= \left( \mc V \otimes \frac{e^{m_{1}}}{\sqrt{t_{1}t_{2}}} \right) \oplus \cdots \oplus \left( \mc V \otimes \frac{e^{m_{2r}}}{\sqrt{t_{1}t_{2}}} \right). 
\end{eqnarray}
Here we consider a homomorphism $\tilde{T}' =\tilde{T} \to \tilde{T}$ defined by
$$
(t_{1}', t_{2}', e^{\mbi a'}, e^{\mbi m'}) \mapsto (t_1, t_{2}, e^{\mbi a}, e^{\mbi m})=((t_{1}')^2, (t_{2}')^2, e^{\mbi a'}, e^{\mbi m'}), 
$$
and use identification $t_{1}'=\sqrt{t_{1}}, t_{2}'=\sqrt{t_{2}}$ and $A^{\ast}_{\tilde{T}'}(\text{pt}) \otimes \mc S \cong \mc S$.
For fixed $r>0$, we put
\begin{eqnarray}
\label{comb0}
\alpha_{n} = \alpha_{n}(\mbi \e, \mbi a, \mbi m) = \int_{M(r,n)} e(\mc F_{r}(\mc V)) = \sum_{\vec{Y}} \frac{\iota_{\vec{Y}}^{\ast} e(\mc F_{r}(\mc V)) }{ e(TM(r,n))}\in \mc S,
\end{eqnarray}
where we omit $r$ in notation $\alpha_{n}$, since in this paper we always fix $r$ and no confusion does not occur.

We put $\e_{+}= \e_{1} + \e_{2}$, and consider $\beta_{n}=\beta_{n}(\mbi \e, \mbi a, \mbi m) = \alpha_{n}(\mbi \e, -\mbi a, -\mbi m)$.
We have our main theorem similar to the conjectured relations \cite[(4.1), (4.2)]{IMO}.
\begin{thm}
\label{main}
We have
$$
\beta_{n}-\alpha_{n}=\sum_{k=1}^{n}(-1)^{k(r+1)}\frac{u_{r}(u_{r}-1)\cdots(u_{r}-k+1)}{k!} \alpha_{n-k},
$$
where $u_{r}=\frac{\e_{+}\left( 2 \sum_{\alpha=1}^{r}  a_{\alpha} + \sum_{f=1}^{2r} m_{f}\right)}{ \e_{1} \e_{2}}$.
\end{thm}
We consider the Nekrasov partition function 
$$
Z(\mbi \e, \mbi a, \mbi m, q)=\sum_{n=0}^{\infty} \alpha_{n}(\mbi \e, \mbi a, \mbi m ) q^n
$$ 
as in the introduction.
Then this theorem says that we have
$$
Z(\mbi \e, -\mbi a, -\mbi m, q)=(1-(-1)^rq)^{u_{r}}Z(\mbi \e,\mbi a,\mbi m,q).
$$
Since $Z(\mbi \e, \mbi a, \mbi m, q)=Z(-\mbi \e, -\mbi a, -\mbi m, q)$ as in \cite[Lemma 6.3 (3)]{NY1}, we get
\begin{equation}
\label{even}
Z(-\mbi \e, \mbi a, \mbi m, q)=(1-(-1)^rq)^{u_{r}}Z(\mbi \e,\mbi a,\mbi m,q).
\end{equation}

%%%%%%%%%%%%%%%%%%%%%%%%%%%%%%%%%%%%%%%%%%%%%%%%%%%%%%%%%%%%%%%%%%%%%%%%

\subsection{Application}
\label{subsec:appl}
For an application of Theorem \ref{main}, we recall combinatorial expression of the Nekrasov partition function $Z(\mbi \e,\mbi a,\mbi m,q)$ from \cite{NY1}.
As in Proposition \ref{ny}, each $\tilde{T}$-fixed point $I_{\vec{Y}}$ in $M(r,n)$ corresponds to $r$-tuple $\vec{Y}=(Y_{1}, \ldots, Y_{r})$ of Young diagrams. 
The fibre of the tautological bundle $\mc V$ over $I_{\vec{Y}}$ is isomorphic to 
\begin{eqnarray}
\label{comb1}
H^{1}(\PP^{2}, I_{\vec{Y}}(-1))=\bigoplus_{\alpha=1}^{r} \bigoplus_{(i,j) \in Y_{\alpha}} e^{a_{\alpha}}t_{1}^{-i+1} t_{2}^{-j+1} 
\end{eqnarray}
as $\tilde{T}$-modules. 

We also recall $\tilde{T}$-modules structures of fibres of $TM(r,n)$ from \cite[Theorem 2.11]{NY1}. 
Let $Y_{\alpha}=\lbrace \lambda_{\alpha, 1}, \lambda_{\alpha,2}, \cdots, \rbrace$ be a Young diagram, where $\lambda_{\alpha,i}$ is the height of the $i$-th column.
We set $\lambda_{\alpha, i}=0$ when $i$ is larger than the width of the diagram $Y_{\alpha}$.
Let $Y_{\alpha}^{T}=\lbrace \lambda_{\alpha, 1}', \lambda_{\alpha,2}', \cdots \rbrace $ be its transpose.
For a box $s=(i,j)$ in the $i$-th column and the $j$-th row, we define its arm-length $a_{Y_{\alpha}}(s)$ and leg-length $l_{Y_{\alpha}}(s)$ with respect to the diagram $Y_{\alpha}$ by $a_{Y_{\alpha}}(s)=\lambda_{\alpha, i}-j$ and $l_{Y_{\alpha}}(s)=\lambda_{\alpha, j}'-i$. 
Then the fibre of $TM(r,n)$ over $I_{\vec{Y}}$ is isomorphic to $\bigoplus_{\alpha,\beta=1}^{r} N_{\alpha,\beta}(t_1,t_2)$ as $\widetilde{T}$-modules, where 
\begin{eqnarray}
\label{comb2}
N_{\alpha, \beta}(t_1,t_2) 
&=&
e_{\beta}e_{\alpha}^{-1} \notag \\
&\times& 
\left\lbrace \bigoplus_{s\in Y_{\alpha}}\left(t_1^{-l_{Y_{\beta}}(s)}t_2^{a_{Y_{\alpha}}(s)+1}\right) \oplus 
\bigoplus_{t\in Y_{\beta}}\left(t_1^{l_{Y_{\alpha}}(t)+1}t_2^{-a_{Y_{\beta}}(t)}\right) \right\rbrace.
\end{eqnarray}

By \eqref{comb0} we have
\begin{eqnarray}
\label{comb}
Z(\mbi \e, \mbi a, \mbi m, q) =\sum_{n=0}^{\infty} \alpha_{n}q^{n}= \sum_{n=0}^{\infty} \sum_{|\vec{Y}|=n}\frac{\prod_{f=1}^{2r} \iota_{\vec{Y}}^{\ast} e\left( \mc V \otimes \frac{e^{m_{f}}}{\sqrt{t_{1}t_{2}}} \right)}{\iota_{\vec{Y}}^{\ast} e \left( TM(r,n) \right)} q^{n}.
\end{eqnarray}
Here by \eqref{comb1} and \eqref{comb2}, we substitute
\begin{eqnarray*}
\iota_{\vec{Y}}^{\ast}e \left( \mc V \otimes \frac{e^{m_{f}}}{\sqrt{t_{1}t_{2}}} \right)
&=&
\prod_{\alpha=1}^{r} \prod_{(i,j) \in Y_{\alpha}} \left(a_{\alpha} + \left(-i+\frac{1}{2} \right)\e_{1} + \left( -j+\frac{1}{2} \right) \e_{2} + m_{f} \right)\\
\iota_{\vec{Y}}^{\ast} e \left( TM(r,n) \right)
&=&
\prod_{\alpha, \beta=1}^{r} \left( \prod_{s\in Y_{\alpha}}\left( a_{\beta}-a_{\alpha}-l_{Y_{\beta}}(s)\e_{1}+(a_{Y_{\alpha}}(s)+1)\e_{2} \right) \times \right.\\
&& \left. \prod_{t\in Y_{\beta}}\left(a_{\beta}-a_{\alpha}+(l_{Y_{\alpha}}(t)+1) \e_{1} -a_{Y_{\beta}}(t) \e_{2}\right) \right).
\end{eqnarray*}

For an application we introduce partition functions
$$
Z^{\psi_{N_{f}}}\left( \mbi \e, \mbi a,  \mbi m^{N_{f}}, q \right)= \sum_{n=0}^{\infty} q^{n} \int_{M(r,n)} \psi_{N_{f}}
$$
 for $0 \le N_{f} \le 2r$, where $\mbi m^{N_{f}} = (m_{1}, \ldots, m_{N_{f}})$ and
$$
\psi_{N_{f}}=
e\left( \bigoplus_{f=1}^{N_{f}} \mc V \otimes \frac{e^{m_{f}}}{\sqrt{t_{1}t_{2}}}\right).
$$ 
We note that $\mbi m = \mbi m^{2r}$ and $Z(\mbi \e, \mbi a,  \mbi m, q)=Z^{\psi_{2r}}(\mbi \e, \mbi a,  \mbi m^{2r}, q)$.
For $N_{f}>0$, we have 
\begin{equation}
\label{reduce}
Z^{\psi_{N_{f}-1}}\left( \mbi \e, \mbi a, \mbi m^{N_{f}-1}, q \right)=Z^{\psi_{N_{f}}}\left( \mbi \e, \mbi a, \mbi m^{N_{f}-1}, \frac{q}{q'}, q' \right) \Big|_{q'=0}.
\end{equation}

When $N_{f}=2r$, by \eqref{even} we have
$$
\log \frac{Z^{\psi_{2r}}( - \mbi \e, \mbi a, \mbi m^{2r}, q)}{Z^{\psi_{2r}}(\mbi \e, \mbi a, \mbi m^{2r}, q)} =u_{r} \log (1 - (-1)^{r}q).
$$
Hence by \eqref{reduce} we have
$$
\log \frac{Z^{\psi_{2r-1}}\left( -\mbi \e, \mbi a, \mbi m^{2r-1}, q\right)}{Z^{\psi_{2r-1}}\left(\mbi \e, \mbi a, \mbi m^{2r-1}, q\right)}= (-1)^{r+1}\frac{\e_{+}}{\e_{1}\e_{2}}q.
$$
Similarly, for $0 \le N_{f} \le 2r-2$ we have
\begin{eqnarray}
\label{odd}
\log \frac{ Z^{N_{f}}\left( -\mbi \e, \mbi a, \mbi m^{N_{f}}, q \right)}{ Z^{N_{f}}\left( \mbi \e, \mbi a, \mbi m^{N_{f}}, q \right) }= 0.
\end{eqnarray}

Hence we determined the odd degree part of $\e_{1}\e_{2} \log Z^{N_{f}}(\mbi \e, \mbi a, \mbi m^{N_{f}}, q)$ with respect to $\e_{1},\e_{2}$.
In particular it is equal to zero unless $N_{f} \ge 2r -1$.

%%%%%%%%%%%%%%%%%%%%%%%%%%%%%%%%%%%%%%%%%%%%%%%%%%%%%%%%%%%%%%%%%%%%%%%%
%%%%%%%%%%%%%%%%%%%%%%%%%%%%%%%%%%%%%%%%%%%%%%%%%%%%%%%%%%%%%%%%%%%%%%%%
%%%%%%%%%%%%%%%%%%%%%%%%%%%%%%%%%%%%%%%%%%%%%%%%%%%%%%%%%%%%%%%%%%%%%%%%
%%%%%%%%%%%%%%%%%%%%%%%%%%%%%%%%%%%%%%%%%%%%%%%%%%%%%%%%%%%%%%%%%%%%%%%%

\section{Reduction to wall-crossing}
\label{sec:redu}
In this section we reduce a proof of Theorem \ref{main} to analysis of wall-crossing phenomena between stable and co-stable ADHM data.

%%%%%%%%%%%%%%%%%%%%%%%%%%%%%%%%%%%%%%%%%%%%%%%%%%%%%%%%%%%%%%%%%%%%%%%%

\subsection{Quiver description of moduli spaces of framed sheaves}
\label{subsec:quiv}
For later purpose we modify the definition of ADHM data following \cite{C}.
As in the previous section, we consider ADHM data on $(W,V)$ for vector spaces $W=\C^{r}, V=\C^{n}$.

We introduce a quiver $\Gamma$, consisting of two verteces $0, \infty$, and two arrows $A_{1}, A_{2}$ from $0$ to $0$, $r$ arrows $\lbrace \gamma_{k}\rbrace_{k=1}^{r}$ from $\infty$ to $0$ and another $r$ arrows $\lbrace \delta_{k}\rbrace_{k=1}^{r}$ from $0$ to $\infty$.
$\Gamma$-representations $X$ consist of vector spaces $X_{0}, X_{\infty}$, and linear maps among $X_{0}$ and $X_{\infty}$ corresponding to $A_{1}, A_{2}, \gamma_{k}$ and $\delta_{k}$ for $k=1,\ldots, r$.
Then elements in $\mb M(W, V)$ are identified with $\Gamma$-representations $X$ with 
$$
X_{0}=V, X_{\infty}=
\begin{cases}
\C \text{ if }W \neq 0, \\
0 \text{ if }W=0,
\end{cases}
$$
via $B_{1}=A_{1}, B_{2}=A_{2}, z=\sum \gamma_{k} \mbi w^{\ast}_{k}, w=\sum \delta_{k} \mbi w_{k}$, where $\mbi w_{1}, \ldots, \mbi w_{r}$ is the canonical basis of $W=\C^{r}$.
For ADHM data, we consider the following relation 
\begin{eqnarray}
\label{relation}
[A_{1},A_{2}] + \sum_{k=1}^{r} \gamma_{k} \delta_{k}=0.
\end{eqnarray}
ADHM data are identified with $\Gamma$-representations satisfying \eqref{relation}.
We also write the vector space $X_{0}\oplus X_{\infty}$ by $X$.

For $\Gamma$-representation $X$, $\theta$-stability conditions are introduced in \cite{K}, where $\theta=(\theta_{0}, \theta_{\infty}) \in \R^{2}$ satisfies $\theta_{0} \dim X_{0} + \theta_{\infty} \dim X_{\infty}=0$. 
Then $X$ is called $\theta$-stable if any non-zero proper sub-representation $P=P_{0} \oplus P_{\infty}$ of $X$ satisfies an inequality $\theta_{0} \dim P_{0} + \theta_{\infty} \dim P_{\infty} < 0$.
Elements $(B_{1},B_{2}, z, w)$ in $\mb M(W,V)$ are called {\it co-stable} if $(^tB_{2}, ^tB_{1}, ^tw, ^tz)$ are stable.
Via the above correspondence between elements in $\mb M(W,V)$ and $\Gamma$-representations, the stability (resp. co-stability) is equivalent to $\zeta(1,-n)$-stability with $\zeta < 0$ (resp. $\zeta>0)$.
We put 
$$
M^{c}(r,n)=\lbrace (B, z,w) \colon \text{ co-stable ADHM data on }(W, V)\rbrace / \GL(V).
$$

To construct moduli spaces, we introduce another affine space $\mb L=\mb L(n)=\mb L(V)=\Hom_{\C}(\wedge Q^{\vee} \otimes V, V)$, and consider a map 
\begin{eqnarray}
\label{mu}
\mu \colon \mb M \to \mb L, (B, z,w) \mapsto \mu(B,z,w)=[B \wedge B]+zw.
\end{eqnarray}
For $\zeta \in \R$, we consider an open locus  
$$
\mu^{-1}(0)^{\zeta}=\lbrace X \in \mu^{-1}(0) \mid X \text{ is } \zeta(-1,n) \text{-stable} \rbrace.
$$

For $\zeta>0$, we have $M(r,n)=[\mu^{-1}(0)^{-\zeta}/G]$ and $M^{c}(r,n)=[\mu^{-1}(0)^{\zeta}/G]$, where $G=\GL(V)$.
In this description the tautological vector bundle $\mc V$ on $M(r,n)$ defined in \S \ref{subsec:toru} is isomorphic to $[\mu^{-1}(0)^{-\zeta} \times V /G]$, where $G$ acts on $V$ naturally.
We also write by $\mc V$ the similar vector bundle $[\mu^{-1}(0)^{\zeta} \times V /G]$ on $M^{c}(r,n)$ and call it tautological bundle.
We define the Uhlenbeck compactification by $M_{0}(r,n)=\Spec \C[\mu^{-1}(0)]^{G}$.
Then via the above construction, we have a $\tilde{T}$-equivariant projective morphism $\pi \colon M(r,n) \to M_{0}(r,n)$. 

%%%%%%%%%%%%%%%%%%%%%%%%%%%%%%%%%%%%%%%%%%%%%%%%%%%%%%%%%%%%%%%%%%%%%%%%
\subsection{Relations between stable and co-stable ADHM data}
\label{subsec:rela1}
We define $\tilde{T}$-equivariant morphisms $D \colon M(r,n) \to M^{c}(r,n)$ by a $G\times \tilde{T}$-equivariant morphism
$$
\mu^{-1}(0)^{-\zeta} \to \mu^{-1}(0)^{\zeta}, (B_{1}, B_{2},  z,w) \mapsto (^tB_{2}, ^tB_{1}, ^tw, ^tz),
$$
via a group homomorphism
$$
G\times \tilde{T} \to G\times \tilde{T}, (g, t_{1}, t_{2}, e^{\mbi a}, e^{\mbi m}) \mapsto (t_{1}t_{2} {}^tg^{-1},t_{2}, t_{1}, e^{-\mbi a}, e^{\mbi m}). 
$$
\begin{NB}
Hence we have $D^{\ast} \mc V= \mc V^{\vee} \otimes t_{1}t_{2}$.
\end{NB}

For each $\tilde{T}$-fixed point $I_{\vec{Y}}$ in $M(r,n)$ as in Proposition \ref{ny}, we consider the following embedding of the co-stable $\tilde{T}$-fixed point 
$$
\iota^{c}_{\vec{Y}} \colon \lbrace D(I_{\vec{Y}}) \rbrace \to M^{c}(r,n).
$$
\begin{NB}
For each $\tilde{T}$-fixed point $I_{\vec{Y}} \in M(r,n)^{\tilde{T}}$ we take a ADHM data $X_{\vec{Y}}=(B_{\vec{Y}},z_{\vec{Y}},w_{\vec{Y}})$ representing $I_{\vec{Y}}$.
Then we have a group homomorphism $g \colon T^{2+r} \to G$ such that $(g(t,a), t,e^{a}) X_{\vec{Y}}=X_{\vec{Y}}$.
If we take a co-stable ADHM data $D(X_{\vec{Y}})$ representing the $\tilde{T}$-fixed point $D(I_{\vec{Y}})$ of $M^{c}(r,n)$, we have $(^{t}g(t,e^{-a})^{-1}t_{1}t_{2}, t,e^{a}) D(X_{\vec{Y}})=D(X_{\vec{Y}})$. 
\end{NB}
We take a decomposition $\iota_{\vec{Y}}^{\ast}e(\mc V)=\prod_{k=1}^{n} F_k(\mbi \e, \mbi a)$ in $\mc S$ by linear polynomials $\lbrace F_k \rbrace_{k=1}^{n}$. 
Then at $D(I_{\vec{Y}})$ in $M^{c}(r,n)$ we have 
\begin{eqnarray*}
(\iota^{c}_{\vec{Y}})^{\ast} e\left( \mc V \otimes \frac{e^{m_{f}}}{\sqrt{t_{1}t_{2}}} \right)
&=&\prod_{k=1}^{n} \left( - F_k(\e_{2}, \e_{1}, -\mbi a)+ \frac{\e_{+}}{2} + m_{f}\right)\\
&=&(-1)^n \prod_{k=1}^{n} \left( F_k(\e_{2}, \e_{1}, -\mbi a) - \frac{\e_{+}}{2} - m_{f}\right)\\
%&=&(-1)^n \prod_{k=1}^{n} \left( F_k(\e, a) - m_{\alpha}\right)|_{\mbi a=-\mbi a, m_{\alpha}=\e_{+}-m_{\alpha}}\\
&=& (-1)^{n} \iota_{\vec{Y}^{t}}^{\ast}e\left( \mc V \otimes \frac{e^{m_{f}}}{\sqrt{t_{1}t_{2}}} \right) \Big |_{\substack{\mbi a=-\mbi a \\ m_{f}=-m_{f}}},
\end{eqnarray*}
where $\vec{Y}^{t}=(Y_{1}^{t}, \ldots, Y_{r}^{t})$ and $Y_{\alpha}^{t}$ are transposes of Young diagrams $Y_{\alpha}$.
Similarly we have $( \iota^{c}_{\vec{Y}})^{\ast} e(TM^{c}(r,n)) = \iota_{\vec{Y}^{t}}^{\ast}e(TM(r,n))|_{\mbi a=-\mbi a}$.
Thus we get
$$
\int_{M^{c}(r,n)} e(\mc F_{r}(\mc V)) = \alpha_{n}(\mbi \e, -\mbi a, -\mbi m)=\beta_{n}(\mbi \e, \mbi a, \mbi m).
$$

%%%%%%%%%%%%%%%%%%%%%%%%%%%%%%%%%%%%%%%%%%%%%%%%%%%%%%%%%%%%%%%%%%%%%%%%
%%%%%%%%%%%%%%%%%%%%%%%%%%%%%%%%%%%%%%%%%%%%%%%%%%%%%%%%%%%%%%%%%%%%%%%%
%%%%%%%%%%%%%%%%%%%%%%%%%%%%%%%%%%%%%%%%%%%%%%%%%%%%%%%%%%%%%%%%%%%%%%%%
%%%%%%%%%%%%%%%%%%%%%%%%%%%%%%%%%%%%%%%%%%%%%%%%%%%%%%%%%%%%%%%%%%%%%%%%

\section{Enhanced master spaces}
\label{sec:enha}
We apply Mochizuki method \cite{M} to the quiver description in the previous section following \cite{NY3}.
The argument in this section is totally similar to \cite{NY3} except that we do not give a sheaf description.
Hence we often omit proofs. 
But for some statements, we give proofs for understandings.
See \cite{NY3} for complete proofs.
%%%%%%%%%%%%%%%%%%%%%%%%%%%%%%%%%%%%%%%%%%%%%%%%%%%%%%%%%%%%%%%%%%%%%%%%
\subsection{ADHM data with full flags}
\label{subsec:adhm}
For vector spaces $W=\C^{r}, V=\C^{n}$, we consider pairs $(X, F^{\bullet})$ of elements $X=(B,z,w)$ in $\mb M(W, V)$ and full flags $F^{\bullet}$ of $V$.
We identify $X$ with $\Gamma$-representations $X_{0} \oplus X_{\infty}$ such that 
$$
X_{0}=V, X_{\infty}=
\begin{cases}
\C \text{ if }W \neq 0, \\
0 \text{ if }W=0,
\end{cases}
$$
as in \S \ref{subsec:quiv}.
We note that $X=(B,z,w)$ does not necessarily satisfy the relation $[B \wedge B] +zw=0$, and $F^{\bullet}$ may have repetitions satisfying $\dim F^{i}/F^{i-1} =0$ or $1$ and $F^{i}=V$ for $i \gg 0$. 
\begin{defn}
\label{ellstab}
For a non-negative integer $\ell \le n$, $(X, F^{\bullet})$ is said to be $\ell$-stable if any sub-representation $P=P_{0} \oplus P_{\infty}$ of $X$ satisfies the following two conditions:
\begin{enumerate}
\item[\textup{(1)}] If $P_{\infty}=0$ and $P \neq 0$, we have $P_{0} \cap F^{\ell} =0$.
\item[\textup{(2)}] If $P_{\infty}=\C$ and $P \neq X$, we have $F^{\ell} \not\subset P_{0}$. 
\end{enumerate}
\end{defn}
We write by $\wt{M}^{\ell}(r,n)$ moduli of $\ell$-stable ADHM data with full flags of $V$, which will be constructed in the next subsection.
We remark that when $\ell=0$ (resp. $\ell=n$), an object $(X, F^{\bullet})$ is $\ell$-stable if and only if $X=(B,z,w)$ is stable (resp. co-stable).
Hence we see that $\wt{M}^{0}(r,n)$ and $\wt{M}^{n}(r,n)$ is the full flag bundle of tautological bundles on $M(r,n)$ and $M^{c}(r,n)$ respectively.  

For later analysis, we need stability parameters.
For $\zeta \in \R$ and $\mbi \eta = (\eta_{1}, \ldots, \eta_{n}) \in \mb Q^{n}_{>0}$, we introduce $(\zeta, \mbi \eta)$-stability for $(X, F^{\bullet})$ as follows.
We take a $\Gamma$-representation $X=X_{0} \oplus X_{\infty}$ with a full flag $F^{\bullet}$ of $X_{0}$.
For a non-zero graded subspace $P=P_{0} \oplus P_{\infty}$ of $X=X_{0} \oplus X_{\infty}$, we define
$$
\mu_{(\zeta, \mbi \eta)}(P)= \frac{\zeta(\dim P_{0} -n \dim P_{\infty} )+ \sum_{i} \eta_{i}\dim (P_{0} \cap F^{i})}{\dim P_{0}+\dim P_{\infty}}.
$$ 
We say that $(X,F^{\bullet})$ is $(\zeta, \mbi \eta)$-{\it semistable} if for any non-zero sub-representation $P$ of $X$, we have
$$
\mu_{(\zeta, \mbi \eta)}(P) \le \mu_{(\zeta, \mbi \eta)}(X).
$$ 
If the inequality is always strict unless $P=X$, we say that $(X, F^{\bullet})$ is $(\zeta, \mbi \eta)$-{\it stable}.

For fixed $\ell \ge 1$, we can take $(\zeta, \mbi \eta)$ satisfying the following two conditions
\begin{eqnarray}
\label{conda}
\eta_{k} 
&>&
(n+1) \sum_{i=k+1}^{n}i\eta_{i} 
\end{eqnarray}
for  each  $k=1, \ldots, n$, and
\begin{multline}
\label{condb}
(n+1)\sum_{i=\ell +1}^{n} i\eta_{i} \\
< \min \left( \sum_{i=1}^{\ell} i\eta_{i} -(n+1) \zeta, \right.
\left. -\sum_{i=1}^{\ell} i\eta_{i} + (n+1)\zeta+\eta_{\ell}\right).
\end{multline}
For example, first we take $\eta_{1},\ldots,\eta_{\ell} \in \mb Q_{>0}$ so that it satisfies \eqref{conda} when we put $\eta_{\ell+1}=\cdots=\eta_{n}=0$.
Then we take $\zeta>0$ so that it satisfies
$$
\sum_{i=1}^{\ell} i\eta_{i} - \eta_{\ell} < (n+1) \zeta < \sum_{i=1}^{\ell}i\eta_{i}.
$$
Finally we take $\eta_{\ell+1}, \ldots, \eta_{n}$ satisfying \eqref{conda} and \eqref{condb}.

Here we only use \eqref{condb}.
The condition \eqref{conda} will be used in \S \ref{subsec:dire} Lemma \ref{lem4}.
\begin{prop}
[\protect{
\cite[Proposition 4.2.4]{M}, 
\cite[Lemma 5.6]{NY3}}]
\label{ell}
Assume that $(\zeta, \mbi \eta)$ satifsfy the condition \eqref{condb}, and $F^{\bullet}$ is a full flag. 
Then for $(X,F^{\bullet})$, the $\ell$-stability is equivalent to the $(\zeta, \mbi \eta)$-stability.
Furthermore the $(\zeta, \mbi \eta)$-semistability automatically implies the $(\zeta,\mbi \eta)$-stability.
\end{prop}
\proof
Suppose that we have a non-zero sub-representation $P$ of $X$ with $\dim P_{0}=p>0, P_{\infty}=0$.
Then $\mu_{(\zeta, \mbi \eta)}(P)\le \mu_{(\zeta, \mbi \eta)}(X)$ means
$$
 \zeta + \frac{\sum_{i=1}^{n} \eta_{i} \dim (P_{0} \cap F^{i})}{p} \le \frac{\sum_{i=1}^{n} i\eta_{i}}{n+1}. 
$$
By \eqref{condb}, this holds if and only if $P_{0}\cap F^{\ell}=0$.
Moreover, the equality never holds.

Next suppose that we have a non-zero proper sub-representation $P$ of $X$ with $\dim X_{0}/P_{0} = p>0, P_{\infty}=\C$.
Then $\mu_{(\zeta, \mbi \eta)}(P)\le \mu_{(\zeta, \mbi \eta)}(X)$ means $\mu_{\zeta,\mbi \eta}(X/P) \ge \mu_{\zeta,\mbi \eta}(X)$, that is, 
$$
 \zeta + \frac{\sum_{i=1}^{n} \eta_{i} \dim (F^{i}/P_{0} \cap F^{i})}{p} \ge \frac{\sum_{i=1}^{n} i\eta_{i}}{n+1}. 
$$
By \eqref{condb} this holds if and only if $F^{\ell} \not\subset P_{0}$.
Moreover, the equality never holds.
\endproof

We also consider the following condition on $\mbi \eta$:
\begin{eqnarray}
\label{n}
\sum_{i=1}^{n} k_{i} \eta_{i} \neq 0 \text{ for any  } (k_{1}, \ldots, k_{n}) \in \Z^{n} \setminus \lbrace 0\rbrace \text{ with } |k_{i}| \le n^{2}. 
\end{eqnarray}
We have the following lemma.
%\begin{lem}
%[\protect{\cite[Lemma 5.14]{NY3}}]
%\label{lem1}
%Assume that $\mbi \eta$ satisfies \eqref{n}.
%If $(X, F^{\bullet})$ is strictly $(\zeta, \mbi \eta)$-semistable, then there exists a non-zero proper subrepresentation $P \subset X$ such that
%\begin{enumerate}
%\item[(1)] $\mu_{(\zeta, \mbi \eta)}(P) = \mu_{(\zeta, \mbi \eta)} (X)$,\\
%\item[(2)] $(P, P_{0} \cap F^{\bullet})$ and $(X/P, F^{\bullet}/(P_{0} \cap F^{\bullet}))$ are $(\zeta, \mbi \eta)$-stable.
%\end{enumerate}
%\end{lem}

\begin{lem}
[\protect{
\cite[Lemma 4.3.9]{M}, 
\cite[Lemma 5.16]{NY3}}]
\label{lem2}
Assume that $\mbi \eta$ satisfies \eqref{n}.
If $(X, F^{\bullet})$ is $(\zeta, \mbi \eta)$-semistable, then its stabilizer is either trivial or $\C^{\ast}$.
In the latter case $(X, F^{\bullet})$ has a unique decomposition $(X_{\flat}, F^{\bullet}_{\flat})\oplus (X_{\sharp}, F^{\bullet}_{\sharp})$ such that both $(X_{\flat}, F^{\bullet}_{\flat})$ and $(X_{\sharp}, F^{\bullet}_{\sharp})$ are $(\zeta, \mbi \eta)$-stable, and the equality $\mu_{(\zeta, \mbi \eta)}(X_{\flat})=\mu_{(\zeta, \mbi \eta)}(X_{\sharp})$ holds.
The stabilizer comes from that of the factor $(X_{\sharp}, F^{\bullet}_{\sharp})$ with $(X_{\sharp})_{\infty}=0$.
\end{lem}

%%%%%%%%%%%%%%%%%%%%%%%%%%%%%%%%%%%%%%%%%%%%%%%%%%%%%%%%%%%%%%%%%%%%%%%%%%%
\subsection{Moduli stacks and $\C^{\ast}_{\hbar}$-action}
\label{subsec:modu}
Let $\underbar{n}$ denote the set $\lbrace 1, \ldots, n\rbrace$ of integers from $1$ to $n$, and $Fl=Fl(V, \underbar{n})$ denote the full flag variety of $V$.
We consider natural projections $\rho_{i}\colon Fl \to G_i=Gr(V,i)$ to Grassmanian manifolds $G_i$ of $i$-dimensional subspace of $V$ and pull-backs $\rho_{i}^{\ast}\mo_{G_i}(1)$ of polarizations $\mo_{G_{i}}(1)$ by Plucker embeddings.

In the following, for fixed $\ell$, we take rational numbers $\zeta^{-} < 0, \zeta>0$ and $\mbi \eta \in \mb Q_{>0}^{n}$ such that $(\zeta, \mbi \eta)$ satisfies \eqref{conda}, \eqref{condb}, and $\mbi \eta$ satisfies \eqref{n}, and $|\zeta|, |\mbi \eta|$ are sufficiently smaller than $|\zeta^{-}|$.
We take a positive integer $k$ enough divisible such that $k\zeta, k\zeta^{-}$ and $k\mbi \eta$ are all integer valued, and consider ample $G$-linearizations 
\begin{eqnarray*}
L_{+}
&=&
\left(\mo_{\mb M} \otimes (\det V)^{\otimes k\zeta} \right)\boxtimes \bigotimes_{i=1}^{n}\rho_{i}^{\ast}\mo_{G_{i}}(k\eta_{i}),\\ 
L_{-}
&=&
\left(\mo_{\mb M} \otimes (\det V)^{\otimes k\zeta^{-}} \right)\boxtimes \bigotimes_{i=1}^{n}\rho_{i}^{\ast}\mo_{G_{i}}(k\eta_{i})
\end{eqnarray*}
on $\wt{\mb M}=\wt{\mb M}(r,n)=\wt{\mb M}(W,V)=\mb M(r,n) \times Fl$.
We consider the composition $\tilde{\mu} \colon \wt{\mb M} \to \mb L$ of the projection $\wt{\mb M} \to \mb M$ and $\mu \colon \mb M \to \mb L$ in \eqref{mu}, and semistable loci $\tilde{\mu}^{-1}(0)^{\ell}$ and $\tilde{\mu}^{-1}(0)^{0}$ with respect $L_{+}$ and $L_{-}$ respectively. 
Then by our choice of $\zeta^{-}, \zeta, \mbi \eta$, we see that $[\tilde{\mu}^{-1}(0)^{\ell}/G]=\wt{M}^{\ell}(r,n)$ and $[\tilde{\mu}^{-1}(0)^{0}/G]$ coincides with the full flag bundle $Fl(\mc V, \underbar{n})$ of the tautological bundle $\mc V$ over $M(r,n)$.

We put $\wh{\mb M}=\wh{\mb M}(r,n)= \mb P_{\wt{\mb M}}(L_{-}\oplus L_{+})$ and consider a composition $\hat{\mu}\colon \wh{\mb M} \to \mb L$ of the projection $\wh{\mb M} \to \mb M$ and $\mu \colon \mb M \to \mb L$.
Then we have a natural $G=\GL(V)$-action on $\wh{\mb M}$ compatible with $\hat{\mu}$.
We consider a quotient stack $\mc N=[\wh{\mb M}(r,n)^{ss}/G]$, where $\wh{\mb M}(r,n)^{ss}$ is a semistable locus with respect to the tautological line bundle $\mo(1)$ on the projective bundle $\wh{\mb M}(r,n)= \mb P_{\wt{\mb M}}(L_{-}\oplus L_{+})$.
We also write by $\mo(1)$ the the restriction of $\mo(1)$ to $\hat{\mu}^{-1}(0)$, which defines semistable locus $\hat{\mu}^{-1}(0)^{ss}$.  
We define an {\it enhanced master space} by $\mc M=[\hat{\mu}^{-1}(0)^{ss}/G]$.
The projection $\hat{\mu}^{-1} (0) \to \mu^{-1}(0)$ induces a proper morphism $\Pi \colon \mc M \to M_{0}(r,n)$.

By Lemma \ref{lem2}, stabilizer groups of $G$-action on $\wh{\mb M}(r,n)^{ss}$ are finite. 
Hence $\mc N$ is a smooth Deligne-Mumford stack containing $\mc M$ and having natural $\tilde{T}$-action.
We also have a $\C^{\ast}_{\hbar}$-action on $\mc N$ defined by 
\begin{eqnarray}
\label{act0}
\left(X,F^{\bullet}, [x_{-},x_{+}] \right) \mapsto \left( X,F^{\bullet}, [e^{\hbar} x_{-},x_{+}] \right),
\end{eqnarray}
where $[x_{-}, x_{+}]$ is the homogeneous coordinate of $\PP(L_{-} \oplus L_{+})$.
These actions keep $\mc M$, and induce $\tilde{T}\times\C^{\ast}_{\hbar}$-action on $\mc M$.

We define $\C^{\ast}_{\hbar}$-fixed points sets $\mc N^{\C^{\ast}_{\hbar}}$ by the zero locus of vector fields on $\mc N$ generated by the $\C^{\ast}_{\hbar}$-action, and $\mc M^{\C^{\ast}_{\hbar}}$ by $\mc M^{\C^{\ast}_{\hbar}} =\mc M \times_{\mc N} \mc N^{\C^{\ast}_{\hbar}}$ (cf. \cite[Appendix C]{GP}).

%%%%%%%%%%%%%%%%%%%%%%%%%%%%%%%%%%%%%%%%%%%%%%%%%%%%%%%%%%%%%%%%%%%%%%%%%%%

\subsection{Direct sum decompositions of fixed points sets}
\label{subsec:dire}
As in \S \ref{subsec:modu}, we take $\zeta^{-} < 0, \zeta>0$ and $\mbi \eta \in \mb Q_{>0}^{n}$ such that $(\zeta, \eta)$ satisfy \eqref{conda}, \eqref{condb}, and \eqref{n}, and $|\zeta|, |\mbi \eta|$ are sufficiently smaller than $|\zeta^{-}|$.

We take $x \in \wh{\mb M}(r,n) \setminus \left( \mb P_{\wt{\mb M}}(L_{-}) \sqcup \mb P_{\wt{\mb M}}(L_{+}) \right)$ over $(X,F^{\bullet}) \in \wt{\mb M}(r,n)$.
\begin{lem}
[\protect{cf. \cite[Sections 3, 4]{T}}]
A point $x\in \wh{ \mb M}(r,n) \setminus \left( \mb P_{\wt{\mb M}}(L_{-}) \sqcup \mb P_{\wt{\mb M}}(L_{+}) \right)$ is semistable if and only if $(X,F^{\bullet})$ is $(\zeta',\mbi \eta)$-semistable for some $\zeta'$ on the segment connecting $\zeta^{-}$ and $\zeta$.
\end{lem}
We assume that $x$ represents a $\C^{\ast}_{\hbar}$-fixed point in $\mc M$, then $(X,F^{\bullet})$ has a non-trivial stabilizer group of $G$.
By Lemma \ref{lem2}, we have a direct sum decomposition $(X, F^{\bullet})=(X_{\flat}, F^{\bullet}_{\flat}) \oplus (X_{\sharp}, F^{\bullet}_{\sharp})$ with $(X_{\sharp})_{\infty}=0$, and $\mu_{(\zeta',\mbi \eta)}(X_{\sharp})=\mu_{(\zeta',\mbi \eta)}(X)$.
We put $I_{\alpha}=\lbrace i \in \underbar{n} \mid F^{i}_{\alpha}/ F^{i-1}_{\alpha} \neq 0 \rbrace$ for $\alpha=\flat, \sharp$ so that $\underbar{n}=I_{\flat} \sqcup I_{\sharp}$.
The datum $\mk J=(I_{\flat}, I_{\sharp})$ is called the {\it decomposition type} of the $\C^{\ast}_{\hbar}$-fixed point represented by $x$.

\begin{lem}
[\protect{
\cite[Lemma 4.4.3]{M}, 
\cite[Lemma 5.25]{NY3}}]
\label{lem3}
We have $\min(I_{\sharp}) \le \ell$.
\end{lem}
\begin{NB}
We see that $\mu_{(\zeta, \eta)}(X_{\sharp}, F_{\sharp}^{\bullet}) < \mu_{(\zeta, \eta)}(X_{\flat}, F_{\flat}^{\bullet})$ is equivalent to the following inequality
$$
\frac{(n-2p+1)\zeta}{n-p+1} + \frac{\sum_{i=\min(I_{\sharp})}^{n} \eta_{i} \dim F_{\sharp}^{i}}{p} < \frac{\sum_{i=1}^{n} \eta_{i} \dim F_{\flat}^{i}}{n-p+1}.
$$ 
By \eqref{condb} we have $\zeta + \sum_{i=\ell +1}^{n} i \eta_{i} < \frac{\sum_{i=1}^{\ell} i \eta_{i}}{n+1}$.
Hence if $\min(I_{\sharp}) > \ell$, then we have $\mu_{(\zeta, \eta)}(X_{\sharp}, F_{\sharp}^{\bullet}) < \mu_{(\zeta, \eta)}(X_{\flat}, F_{\flat}^{\bullet})$.
On the other hand, we have $\mu_{(\zeta^{-}, \eta)}(X_{\sharp}, F_{\sharp}^{\bullet}) < \mu_{(\zeta^{-}, \eta)}(X_{\flat}, F_{\flat}^{\bullet})$.
Hence we can not have $\zeta^{-} < \zeta' < \zeta^{-}$ such that $\mu_{(\zeta', \eta)}(X_{\sharp}, F_{\sharp}^{\bullet}) = \mu_{(\zeta', \eta)}(X_{\flat}, F_{\flat}^{\bullet})$.
This is a contradiction.
\end{NB}

Conversely  suppose that an object $(X, F^{\bullet})=(X_{\flat}, F^{\bullet}_{\flat}) \oplus (X_{\sharp}, F^{\bullet}_{\sharp})$ with the decomposition type $\mk J=(I_{\flat}, I_{\sharp})$ with $\min(I_{\sharp}) \le \ell$ is given.
Then the condition \eqref{condb} implies $\mu_{(\zeta, \mbi \eta)}(X_{\sharp}) > \mu_{(\zeta, \mbi \eta)}(X)$ as in the proof of Proposition \ref{ell}.
On the other hand, since $|\mbi \eta|$ is enough smaller than $|\zeta^{-}|$, we have 
$$
\mu_{(\zeta^{-},\mbi \eta)}(X_{\sharp})=\zeta^{-} + \frac{\sum_{i=1}^{n} \eta_{i} \dim (F_{\sharp}^{i})}{\dim (X_{\sharp})_{0}} < \mu_{(\zeta^{-},\mbi \eta)}(X)=\frac{\sum_{i=1}^{n} i\eta_{i}}{n+1}. 
$$
Hence we can find $\zeta'$ on the segment connecting $\zeta$ and $\zeta^{-}$ such that $\mu_{(\zeta',\mbi \eta)}(X_{\sharp}) = \mu_{(\zeta',\mbi \eta)}(X)$.

We introduce another stability condition for pairs $(X_{\sharp}, F_{\sharp}^{\bullet})$ of elements $X_{\sharp}$ in $\mb M(0, V_{\sharp})$ and full flags $F_{\sharp}^{\bullet}$ of $V_{\sharp}$.
\begin{defn}
\label{+}
$(X_{\sharp}, F_{\sharp}^{\bullet})$ is said to be $+$-stable if for any proper sub-representation $S =S_{0} \oplus S_{\infty}$, we have $S_{0} \cap F^{1}=0$.
\end{defn}

\begin{lem}
[\protect{\cite[Proposition 4.4.4]{M},\cite[Lemma 5.26]{NY3}}]
\label{lem4}
We have the following.
\begin{enumerate}
\item[\textup{(1)}] $(X_{\flat}, F_{\flat}^{\bullet})$ is $(\zeta', \mbi \eta)$-stable if and only if it is $(\min(I_{\sharp})-1)$-stable.

\item[\textup{(2)}] $(X_{\sharp}, F_{\sharp}^{\bullet})$ is $(\zeta', \mbi \eta)$-stable if and only if it is $+$-stable.
\end{enumerate}
\end{lem}
\proof
(1) Let $S =S_{0} \oplus S_{\infty}$ be a sub-representation of $X_{\flat}$.
We first suppose $S_{\infty}=0$.
Then the inequality $\mu_{\zeta',\mbi \eta}(S) < \mu_{\zeta', \mbi \eta}(X_{\flat}) = \mu_{\zeta', \mbi \eta}(X_{\sharp})$ is equivalent to
$$
\frac{\sum_{i} \eta_{i} \dim (S_{0} \cap F_{\flat}^{i}) }{\dim S_{0}} < \frac{\sum_{i} \eta_{i} \dim ( F_{\sharp}^{i}) }{\dim (X_{\sharp})_{0}}.
$$
Since $\eta_{i}$ for $i \ge \min(I_{\sharp})$ is much smaller than $\eta_{\min(I_{\sharp}) - 1}$ by \eqref{conda}, if the inequality holds, then we must have $S_{0} \cap F^{\min(I_{\sharp})-1}_{\flat}=0$.
Conversely if we have $S_{0} \cap F^{\min(I_{\sharp})-1}_{\flat}=0$, then since $F^{\min(I_{\sharp}) -1}_{\flat} = F^{\min(I_{\sharp})}_{\flat}$ again by \eqref{conda} the above inequality holds.

Next, suppose $S_{\infty}=\C$.
Then the inequality $\mu_{\zeta',\mbi \eta}(S) < \mu_{\zeta', \mbi \eta}(X_{\flat}) = \mu_{\zeta', \mbi \eta}(X_{\sharp})$ is equivalent to
$$
\frac{\sum_{i} \eta_{i} \dim (F_{\flat}^{i} / S_{0} \cap F_{\flat}^{i}) }{\dim ( X_{0}/S_{0})} > \frac{\sum_{i} \eta_{i} \dim ( F_{\sharp}^{i}) }{\dim (X_{\sharp})_{0}}.
$$
This is equivalent to $F_{\flat}^{\min(I_{\sharp})-1} \not\subset S_{0}$ by the same argument as the above.
Thus $(X_{\flat}, F_{\flat}^{\bullet})$ is $(\zeta', \eta)$-stable if and only if it is $(\min(I_{\sharp})-1)$-stable.

(2) It follows by the same argument as above.
\endproof

We assume that ADHM data $(X_{\flat}, F^{\bullet}_{\flat}), (X_{\sharp},F^{\bullet}_{\sharp})$ satisfy conditions in the above lemma.
Let $V=V_{\flat} \oplus V_{\sharp}$ be a direct sum of underlying vector spaces $V_{\flat}, V_{\sharp}$ of $X_{\flat}, X_{\sharp}$. 
We consider a group action $\C^{\ast}_{\frac{\hbar}{pD}} \times \mc M \to \mc M$ defined by 
\begin{eqnarray}
\label{act}
\left( X,F^{\bullet}, [x_{-},x_{+}] \right) \mapsto \left( \id_{V_{\flat}} \oplus e^{\frac{h}{pD}} \id_{V_{\sharp}} \right) \left(X,F^{\bullet}, [e^{\hbar}x_{-},x_{+}] \right). 
\end{eqnarray}
This action is equal to the original $\C^{\ast}_{\hbar}$-action \eqref{act0}, since the difference is absorbed in $G$-action.
Then $(X_{\flat} \oplus X_{\sharp}, F^{\bullet}_{\flat} \oplus F^{\bullet}_{\sharp}, [x_{-}, x_{+}])$ is fixed by this $\C^{\ast}_{\hbar}$-action, and represents a $\C^{\ast}_{\hbar}$ -fixed point in $\mc M$.

%%%%%%%%%%%%%%%%%%%%%%%%%%%%%%%%%%%%%%%%%%%%%%%%%%%%%%%%%%%%%%%%%%%%%%%%%%%

\subsection{Moduli stacks of fixed points sets}
\label{subsec:fixe}

We put
$$
S^{\ell}=\lbrace (I_{\flat}, I_{\sharp}) \mid I_{\flat} \sqcup I_{\sharp} = [n], I_{\sharp} \neq \emptyset, \min(I_{\sharp}) \le \ell\rbrace.
$$
We mainly consider moduli stacks for ADHM data.
But we note that for all objects in $\wh{\mb M}=\wh{\mb M}(r,n)$, arguments in the previous subsections hold.

In particular, the ambient stack $\mc N= [\wh{\mb M}^{ss}/G]$ of the enhanced master space $\mc M$ is a smooth Deligne-Mumford stack with $\tilde{T} \times \C^{\ast}_{\hbar}$-action compatible with one on $\mc M$.
Furthermore, we have a decomposition
$$
\mc N^{\C^{\ast}_{\hbar}}= \mc N_{+} \sqcup \mc N_{-} \sqcup \bigsqcup_{\mk J \in S^{\ell}} \mc N_{\mk J},
$$
where $\mc N_{\pm}$ is defined by the zero locus of $x_{\mp}$, and $\mc N_{\mk J}$ is the closed substack of $\mc N^{\C^{\ast}_{\hbar}}$ of elements over direct sum $(X_{\flat}, F^{\bullet}_{\flat}) \oplus (X_{\sharp}, F^{\bullet}_{\sharp})$ with the decomposition type $\mk J$ as in the previous subsection.
\begin{NB}
For a scheme $S$, the full sub-category $\mc N_{\mk J}(S)$ consists of objects $\phi \in \mc N(S)$ such that the pull-back of vector fields corresponding to $\C^{\ast}_{\hbar}$-action vanish, and for any closed point $p \in S$, the corresponding data $(X, F^{\bullet}, [\rho, 1])$ satisfies that $(X, F^{\bullet})$ is $(X_{\flat}, F^{\bullet}_{\flat}) \oplus (X_{\sharp}, F^{\bullet}_{\sharp})$ with the decomposition type $\mk J$.
\end{NB}

In the following, we show that for $\alpha=\pm$, or $\mk J \in S^{\ell}$, we can take closed subsets $\wh{\mb M}_{\alpha}$ of $\wh{\mb M}$ and closed subgroups $G_{\alpha}$ of $G$ such that we have \begin{eqnarray}
\label{ambient}
[\wh{\mb M}_{\alpha}^{ss}/G_{\alpha}] \cong \mc N_{\alpha}
\end{eqnarray} 
via the natural morphisms.
Here $\wh{\mb M}_{\alpha}^{ss}=\wh{\mb M}_{\alpha} \times_{\wh{\mb M}} \wh{\mb M}^{ss}$.
For $\pm$, we define $\wh{\mb M}_{\pm}$ as zero loci of $x_{\mp}$ in $\wh{\mb M}$, and $G_{\alpha}=G$.

For $\mk J = (I_{\flat}, I_{\sharp}) \in S^{\ell}$, we fix a direct sum decomposition $V=V_{\flat} \oplus V_{\sharp}$ such that $V_{\flat}=\C^{n-p}, V_{\sharp}=\C^{p}$, where $p=|I_{\sharp}|$.
The decomposition type $\mk J$ gives identifications $\underbar{n-p} \cong I_{\flat}, \underbar{p} \cong I_{\sharp}$, and defines a closed embedding 
$$
Fl_{\mk J} \colon Fl(V_{\flat}, \underbar{n-p}) \times Fl(V_{\sharp}, \underbar{p}) \subset Fl(V, \underbar{n}).
$$
This gives an embedding 
$
\wt{\mb M}(W, V_{\flat}) \times \wt{\mb M}(0, V_{\sharp}) \to \wt{\mb M}(r,n).
$
We put 
$$
\wh{\mb M}_{\mk J} = \left( \wt{\mb M}(W, V_{\flat}) \times \wt{\mb M}(0, V_{\sharp}) \right) \times_{\wt{\mb M}(r,n)} \wh{\mb M}(r,n)
$$ 
and $G_{\mk J}=\GL(V_{\flat}) \times \GL(V_{\sharp})$.
Then we have an isomorphism \eqref{ambient} since elements in $\mc N_{\mk J}$ are included in the image of $\wh{\mb M}_{\mk J}$ up to $G$-actions, and their stabilizers are included in $\C^{\ast} \id_{V_{\sharp}} \subset \GL(V_{\flat}) \times \GL(V_{\sharp})$ by Lemma \ref{lem2}.
\begin{NB}
This follows from the following lemma.
\begin{lem}
Let $U$ be a variety with a group $G$-action, $V$ a closed subvariety of $U$, and $H$ a closed subgroup of $G$ such that $H V \subset V$.
We assume that there exists a substack $\mc Z$ of $[U/G]$ such that the natural inclusion $[V/H] \to [U/H]$ factors through the inclusion $\mc Z \to [U/G]$. 
If all objects in $\mc Z$ are included in the image of $[V/H]$, and 
$$
H \supset \lbrace g\in G \mid g V \cap V \neq \emptyset \rbrace,
$$
then the injection $[V/H] \to \mc Z$ is an isomorphism.
\end{lem}

Indeed, by the stability condition we have an open locus where there exist $m_{1}, \ldots, m_{n} \in (\Z_{\ge 0})^{2}$ such that $\dim \langle B^{m_{1}} F^{1}, \ldots, B^{m_{p}} F^{1} \rangle =p$ and $\dim \langle B^{m_{p+1}} z, \ldots, B^{m_{n}} z\rangle  = n-p$.
Then we have an isomorphism
$
[\bar{U} / \bar{G} ] \cong [U/G],
$
where 
$$
\bar{U}= \lbrace \langle B^{m_{1}} F^{1}, \ldots, B^{m_{p}} F^{1} \rangle = V_{\flat}, \langle B^{m_{p+1}} z, \ldots, B^{m_{n}} z\rangle  = V_{\sharp}\rbrace
$$
and
$\bar{G} = \GL(V_{\flat}) \times \GL(V_{\sharp})$.
Then on $[\bar{U}/\bar{G}]$, we see that the zero-locus of vector fields corresponding to $\C^{\ast}_{\hbar}$-action is equal to $\wh{\mb M}_{\mk J}$.  
\end{NB}

Since $\mc M^{\C^{\ast}_{\hbar}} = \mc M \times_{\mc N} \mc N^{\C^{\ast}_{\hbar}}$, we have a decomposition, 
$$
\mc M^{\C^{\ast}_{\hbar}}= \mc M_{+} \sqcup \mc M_{-} \sqcup \bigsqcup_{\mk J \in S^{\ell}} \mc M_{\mk J},
$$
where $\mc M_{\alpha}=\mc M \times_{\mc N} \mc N_{\alpha}$ for $\alpha=\pm$, or $\mk J \in S^{\ell}$. 
If we write $\hat{\mu}_{\alpha} \colon \wh{\mb M}_{\alpha} \to \mb L$ the restriction of $\hat{\mu}$ to $\wh{\mb M}_{\alpha}$, then we have $\mc M_{\alpha} = [\hat{\mu}_{\alpha}^{-1}(0) \cap \wh{\mb M}_{\alpha}^{ss} / G_{\alpha}]$.
In \S \ref{subsec:loca1}, we will introduce a terminology {\it Kuranishi charts} to describe $\mc M_{\alpha}$ again.

%%%%%%%%%%%%%%%%%%%%%%%%%%%%%%%%%%%%%%%%%%%%%%%%%%%%%%%%%%%%%%%%%%%%%%%%%%%

\subsection{Decomposition of $\mc M_{\mk J}$}
\label{subsec:deco}
For $\mk J = (I_{\flat}, I_{\sharp}) \in S^{\ell}$, we further decompose $\mc M_{\mk J}$ as follows.
We fix a direct sum decomposition $V=V_{\flat} \oplus V_{\sharp}$ such that $V_{\flat}=\C^{n-p}, V_{\sharp}=\C^{p} $, where $p=|I_{\sharp}|$.

We consider a moduli stack 
$$
M^{+}_{p} = [\left( \tilde{\mu}^{-1}(0)^{+} \times \C_{\rho_{\sharp}}\right) / \GL(V_{\sharp})]
$$ 
parametrizing tuples $(X_{\sharp}, F^{\bullet}_{\sharp}, \rho_{\sharp})$ of $+$-stable pairs $(X_{\sharp}, F_{\sharp}^{\bullet})$ and orientations $\rho_{\sharp} \colon \det V_{\sharp}^{\otimes D} \cong \C$ for $D=k(\zeta - \zeta^{-})\in \Z$.
Here $\tilde{\mu} \colon \wt{\mb M}(0, V_{\sharp}) \to \mb L(V_{\sharp})$ is defined in \S \ref{subsec:modu}, $\tilde{\mu}^{-1}(0)^{+}$ is the $+$-stable locus, and $g_{\sharp} \in \GL(V_{\sharp})$ acts by $g_{\sharp} \cdot \rho_{\sharp}=(\det g_{\sharp})^{-D}\rho_{\sharp}$.

We also consider a moduli stack 
$$
\bar{M}=[\wt{\mu}^{-1}(0)^{\min(I_{\sharp})-1} \times \C^{\ast}_{\rho_{\flat}}/\GL(V_{\flat})], 
$$
where $\wt{\mu} \colon \wt{\mb M}(W,V_{\flat}) \to \mb L(V_{\flat})$, and $\wt{\mu}^{-1}(0)^{\min(I_{\sharp})-1}$ is the $\left( \min(I_{\sharp}) -1 \right)$-stable locus.
Here $g_{\flat} \in \GL(V_{\flat})$ acts by $g_{\flat}\cdot\rho_{\flat}=(\det g_{\flat})^{-D}\rho_{\flat}$.
We consider a quotient $\left( \bar{M} \times M^{+}_{p} \right)/ \C^{\ast}_{t}$,  where $\C^{\ast}_{t}$-action is defined by $(t \rho_{\flat}, t^{-1}\rho_{\sharp})$, and the other data are fixed.
By Lemma \ref{lem4}, we have an isomorphism 
\begin{eqnarray}
\label{mkj}
\left( \bar{M} \times M^{+}_{p} \right)/ \C^{\ast}_{t} \cong \mc M_{\mk J},
\end{eqnarray}
sending $(X_{\flat}, F^{\bullet}_{\flat}, \rho_{\flat}, X_{\sharp}, F^{\bullet}_{\sharp}, \rho_{\sharp})$ to $(X_{\flat} \oplus X_{\sharp}, Fl_{\mk J}( F^{\bullet}_{\flat}, F^{\bullet}_{\sharp}), [\rho_{\flat} \rho_{\sharp}, 1])$.

We write by $\mc V^{\mk J}_{\flat}, \mc V^{\mk J}_{\sharp}$ tautological bundles on $\mc M_{\mk J}$ corresponding to $V_{\flat}=\C^{n-p}, V_{\sharp}=\C^{p}$ via the isomorphism \eqref{mkj}, such that we have $\mc V|_{\mc M_{\mk J}} = \mc V^{\mk J}_{\flat} \oplus \mc V^{\mk J}_{\sharp}$ for the tautological bundle $\mc V$ on $\mc M$.
We also write by $\mc V_{\flat}, \mc V_{p}^{+}$ tautological bundles on $M^{\min(I_{\sharp})-1}(r,n-p), M^{+}_{p}$ corresponding to $V_{\flat}=\C^{n-p}, V_{\sharp}=\C^{p}$ and their pull-backs to $M^{\min(I_{\sharp})-1}(r,n-p) \times M^{+}_{p}$ by projections.

\begin{thm}
[\protect{\cite[Theorem 5.18]{NY3}}]
\label{decomp}
For $\C^{\ast}_{\hbar}$-action on $\mc M$ defined by \eqref{act0}, we have a decomposition
$$
\mc M^{\C^{\ast}_{\hbar}}=\mc M_{+} \sqcup \mc M_{-} \sqcup \bigsqcup_{\mk J \in S^{\ell}}\mc M_{\mk J}
$$
such that the followings hold.\\
(i) We have $\mc M_{+} \cong \wt{M}^{\ell}(r,n)$ and $\mc M_{-} \cong Fl(\mc V, \underbar{n})$ the full flag bundle of the tautological bundle $\mc V$ over $M(r,n)$.\\
(ii) For each $\mk J=(I_{\flat}, I_{\sharp}) \in S^{\ell}$, we have finite \'etale morphisms 
$$
F \colon \mc S_{\mk J} \to \mc M_{\mk J}, G \colon \mc S_{\mk J} \to \widetilde{M}^{\min(I_{\sharp})-1}(r,n-p) \times M^{+}_{p}
$$ 
of degree $\frac{1}{pD}$, where $p=|I_{\sharp}|, D=k(\zeta - \zeta^{-})$. \\
(iii) There exists a line bundle $L_{\mc S_{\mk J}}$ on $\mc S_{\mk J}$ such that we have isomorphisms $L_{\mc S_{\mk J}}^{\otimes pD} \cong G^{\ast} \det \mc V_{\flat}^{D}$, $F^{\ast} \mc V^{\mk J}_{\flat} \cong G^{\ast} \mc V_{\flat}$, and $F^{\ast} \mc V^{\mk J}_{\sharp} \cong G^{\ast} \mc V_{p}^{+} \otimes e^{\frac{\hbar}{pD}} \otimes L_{\mc S_{\mk J}}^{\vee}$ as $\C^{\ast}_{\hbar}$-equivariant vector bundles on $\mc S_{\mk J}$.
\end{thm}

\proof
We have to only show assertions (ii) and (iii).
For each $\mk J=(I_{\flat}, I_{\sharp}) \in S^{\ell}$, we fix a decomposition $V=V_{\flat} \oplus V_{\sharp}$ as above.

In the following, we consider quotient stacks $[U/H]$ for various groups $H$, and $\tilde{T} \times \C^{\ast}_{\hbar}$-actions on $[U/H]$ induced from $H \times \tilde{T} \times \C^{\ast}_{\hbar}$-actions on $U$.
We identify $H \times \tilde{T} \times \C^{\ast}_{\hbar}$-representations $A$ with associated $\tilde{T} \times \C^{\ast}_{\hbar}$-equivariant vector bundles $[(U\times A)/H]$ on $[U/H]$.

By $\C^{\ast}_{\frac{\hbar}{pD}}$-action \eqref{act} on $\mc M$, we have $V|_{\mc M_{\mk J}} = V_{\flat} \oplus \left( V_{\sharp} \otimes e^{\frac{\hbar}{pD}} \right)$.
We put $\mc S_{\mk J}=\left( \bar{M} \times M^{+}_{p} \right)/ \C^{\ast}_{s}$, where $s=t^{\frac{1}{pD}}$.
Then $\C^{\ast}_{s} \to \C^{\ast}_{t}, s \mapsto t=s^{pD}$ induces a finite morphism $F \colon \mc S_{\mk J} \to \mc M_{\mk J}$ of degree $\frac{1}{pD}$. 
Here we have $(X_{\sharp}, s^{-pD} \rho_{\sharp})= s \id_{V_{\sharp}}(X_{\sharp}, \rho_{\sharp})$ and $\C^{\ast}_{s}$ acts on $M^{+}_{p}$ trivially.
Hence we have an isomorphsim $\mc S_{\mk J}\cong \bar{M}/\C^{\ast}_{s} \times M^{+}_{p}$ induced by a group automorphism of $\GL(V_{\flat}) \times \GL(V_{\sharp}) \times \C^{\ast}_{s}$ defined by
\begin{eqnarray}
\label{gphom}
 (g_{\flat}, g_{\sharp}, s) \mapsto (g_{\flat}, s^{-1} g_{\sharp}, s).
\end{eqnarray} 
Hence we see that $F^{\ast} \left(V|_{\mc M_{\mk J}}\right)$ corresponds to $V_{\flat} \oplus  \left(V_{\sharp} \otimes e^{\frac{\hbar}{pD}} \otimes s^{-1} \right)$ via the above isomorphism. 

On the other hand, we have $\bar{M}/\C^{\ast}_{u}  \cong \wt{M}^{\min(I_{\sharp})-1}(r,n-p)$ and a finite morphism $G \colon \mc S_{\mk J} \to \wt{M}^{\min(I_{\sharp})-1}(r,n-p)\times M^{+}_{p}$ induced by $\C^{\ast}_{s} \to \C^{\ast}_{u}, s \mapsto u=s^{pD}$.
The line bundle $(\det V_{\flat})^{D}$  on $\bar{M}/\C^{\ast}_{u}  \cong \wt{M}^{\min(I_{\sharp})-1}(r,n-p)$ corresponds to $\C \otimes u$, hence $G^{\ast}(\det V_{\flat})^{D}$ corresponds to $\C \otimes s^{pD}$ on $\bar{M}/\C^{\ast}_{s} \times M^{+}_{p}$.
If we put $L_{\mc S} = \C \otimes s$ on $\mc S_{\mk J}$, we have $L_{\mc S}^{\otimes pD} = G^{\ast} (\det V_{\flat})^{D}$, and get two isomorphisms $F^{\ast}V_{\flat} \cong G^{\ast}V_{\flat}$ and $F^{\ast} V_{\sharp} \cong G^{\ast} V_{\sharp}\otimes e^{\frac{\hbar}{pD}} \otimes L_{\mc S}^{\vee}$ on $\mc S_{\mk J}$.
\endproof

%%%%%%%%%%%%%%%%%%%%%%%%%%%%%%%%%%%%%%%%%%%%%%%%%%%%%%%%%%%%%%%%%%%%%%%%%%%
%%%%%%%%%%%%%%%%%%%%%%%%%%%%%%%%%%%%%%%%%%%%%%%%%%%%%%%%%%%%%%%%%%%%%%%%%%%
%%%%%%%%%%%%%%%%%%%%%%%%%%%%%%%%%%%%%%%%%%%%%%%%%%%%%%%%%%%%%%%%%%%%%%%%%%%

\section{Obstruction theories}
\label{sec:obst}
In this section we compute and compare obstruction theories among moduli stacks.
We recall that an obstruction theory for a Deligne-Mumford stack $\mc Z$ is a homomorphism $ob_{\mc Z} \colon Ob_{\mc Z} \to L_{\mc Z}$ to the cotangent complex $L_{\mc Z}$ of $\mc Z$ in the derived category $D(\mc Z)$ of quasi-coherent sheaves on $\mc Z$ satisfying the following conditions. 
The cohomology $\mc H^{i}(Ob_{\mc Z})$ of the complex $Ob_{\mc Z}$ is coherent for $i=-1,0,1$, $\mc H^{i}(ob_{\mc Z})$ are isomorphisms for $i \ge 0$, and $\mc H^{-1}(ob_{\mc Z})$ is surjective.
It is called perfect, if it is quasi-isomorphic to a complex of locally free sheaves $F^{-1} \to F^{0} \to F^{1}$ in the derived category $D(\mc Z)$.
For more details, see \cite[\S 2.4]{M}.

%%%%%%%%%%%%%%%%%%%%%%%%%%%%%%%%%%%%%%%%%%%%%%%%%%%%%%%%%%%%%%%%%%%%%%%%%%%

\subsection{Virtual fundamental cycle}
\label{subsec:virt}
We recall constructions of virtual fundamental cycles \cite{BF} in the following setteing.
In this paper, constructions of moduli stacks are essentially fitted in \cite[\S 2]{GP}, but slightly different point is that we include group actions as in \cite[\S 2.4]{M}.

We consider a tuple $(Y, E, e, H, M)$. 
Here $Y$ is a smooth scheme, $H$ is a group scheme acting on $Y$, $E$ is a $H$-equivariant vector bundle over $Y$, $e$ is a $H$-invariant section of $E$, and $M$ is an ample $H$-linearization.
We take a semistable locus $Y^{ss}$ of $Y$ with respect to $M$, and put $Z=e^{-1}(0) \cap Y^{ss}$.
We assume that any point in $Y^{ss}$ has a finite stabilizer group.
We define a moduli stack $\mc Z=[Z/H]$ as the quotient stack, and call $(Y, E, e, H, M)$ a {\it Kuranishi chart} of $\mc Z$.
In our case, $E$ is always a trivial bundle $Y \times \C^{l}$, and we have a $H$-invariant function $\varphi \colon Y \to \C^{l}$ to write $e=\id_{Y} \times \varphi \in \Gamma(Y, E)$. 

To construct a perfect obstruction theory $ob_{\mc Z} \colon Ob_{\mc Z} \to L_{\mc Z}$, we identify coherent sheaves on $\mc Z$ with $H$-equivariant coherent sheaves on $Z$.
We define a complex $Ob_{\mc Z}$ by the following $H$-equivariant complex:
\begin{eqnarray}
\label{cot}
0 \to E^{\vee}|_{Z}  \stackrel{e^{\vee}}{\to} \Omega_{Y}|_{Z} \stackrel{\nu}{\to} \mathfrak{h}^{\vee} \otimes \mo_{Z} \to 0,
\end{eqnarray}
where $\Omega_{Y}|_{Z}$ stands on degree $0$, $\mathfrak{h}$ is the Lie algebra of $H$, and $\nu$ is induced by $H$-action.

We have a natural homomorphism $ob_{\mc Z} \colon Ob_{\mc Z} \to L_{\mc Z}$.
This gives a perfect obstruction theory.
When $e^{-1}(0)$ is a complete intersection, then by \cite[Corollaire 3.2.7 in chapter III]{I} and \cite[Proposition 2.3.3]{M}, the obstruction theory $ob_{\mc Z}$ is a quasi-isomorphism.
In this case we write $ob_{\mc Z}=\id_{L_{\mc Z}}$.
This complex is the dual of the tangent complex in \cite[\S 4.2]{NY3}, when we consider moduli of framed sheaves on the blow-up of the plane.

We take another Kuranishi chart $(Y_{1}, E_{1}, e_{1}, H_{1}, M_{1})$ of another stack $\mc Z_{1}=[\left( e_{1}^{-1}(0) \cap Y_{1}^{ss} \right) / H_{1}]$, where $Y_{1}^{ss}$ is a semi-stable locus with respect to a $H_{1}$-linearization $M_{1}$. 
We assume that there exists a bundle morphism 
$$
\xymatrix{
\ar[d] E_{1} \ar[r] & E \ar[d]\\
Y_{1} \ar[r] & Y
}
$$ 
equivariant via a group homomorphism $H_{1} \to H$ such that sections $e_{1}, e$ are equal via the isomorphism $V_{1} \cong V|_{Y_{1}}$.
This collection of morphisms are called a morphism of Kuranishi charts.
Furthermore we assume that these morphisms induces an \'{e}tale morphism $[Y_{1}^{ss}/H_{1}] \to [Y^{ss}/H]$. 
In this case, we call such a collection of morphisms a {\it coordinate change} of Kuranishi charts. 

Then any coordinate change induces \'{e}tale morphism $\Psi \colon \mc Z_{1} \to \mc Z$ such that the homomorphism $\Psi^{\ast} Ob_{\mc Z} \to Ob_{\mc Z_{1}}$ induced by the bundle morphism is quasi-isomorphism.
This is because that complexes \eqref{cot} are  restrictions to $\mc Z$ of complexes quasi-isomorphic to duals of differentials of sections $e_{1}$ and $e$ of vector bundles on smooth ambient spaces $[Y_{1}^{ss}/H_{1}]$ and $[Y^{ss}/ H]$.
Hereafter, we freely use this fact and identify complexes $Ob_{\mc Z_{1}}$ and $\Psi^{\ast} Ob_{\mc Z}$.
On the other hand, if we remove the assumption that $\Psi$ is \'{e}tale, then obstruction theories and their defining virtual fundamental cycles can be different. 
This will happen in \S \ref{subsec:obst4}.

We suppose that $Y$ and $V$ admit $\tilde{T}\times\C^{\ast}_{\hbar}$-actions compatible with $H$-actions such that the section $e$ is also $\tilde{T}\times\C^{\ast}_{\hbar}$-equivariant.
Then $\tilde{T}\times\C^{\ast}_{\hbar}$ acts on $\mc Z$, and $Ob_{\mc Z}$ and $L_{\mc Z}$ have $\tilde{T} \times \C^{\ast}_{\hbar}$-equivariant structures such that $ob_{\mc Z}$ is $\tilde{T}\times\C^{\ast}_{\hbar}$-equivariant.
We take a zero section $0 \in \Gamma(Y^{ss}, V)$, and consider the Gysin homomorphism
$$
0^{!} \colon A^{H \times \tilde{T} \times \C^{\ast}_{\hbar}}_{\ast}(Y^{ss} ) \to A^{H \times \tilde{T} \times \C^{\ast}_{\hbar}}_{\ast}(Z) \cong A_{\ast}^{\tilde{T} \times \C^{\ast}_{\hbar}}(\mc Z),
$$
where $Z=e^{-1}(0) \cap Y^{ss}$, and the last isomorphism follows from \cite[Remark 2.1.17]{Kr}.
The virtual fundamental cycle $[\mc Z]^{vir}$ is defined by 
$$
[\mc Z]^{vir} = 0^{!} [Y^{ss}] \in A_{\ast}^{\tilde{T} \times \C^{\ast}_{\hbar}}(\mc Z),
$$ 
where $[Y^{ss}]$ is the fundamental cycle of $Y^{ss}$.
When we have a morphism $(Y_{1}, E_{1}, e_{1}, H_{1}, M_{1}) \to (Y, E, e, H, M)$ of Kuranishi charts, such that $Y_{1} \to Y$ is flat, then we have $[\mc Z_{1}] = \Psi^{\ast} [\mc Z]$.
\begin{NB}
By \cite[Theorem 2.6 (b), (c)]{F} we have,
\begin{eqnarray*}
[\mc Z_{1}]^{vir}
&=&
0_{1}^{!}[Y_{1}]\\
&=&
0_{1}^{!}[Y|_{Y_{1}}]\\
&=&
0^{!}[Y|_{Y_{1}}]\\
&=&
\Psi^{\ast}0^{!}[Y]\\
&=&
\Psi^{\ast}[\mc Z]^{vir},
\end{eqnarray*}
where $0_{1} \colon Y_{1} \to E_{1}$ is the zero-section.
\end{NB}

For example, if we put $Y=\wh{\mb M} , V=\wh{\mb M} \times \mb L, v=\id_{\wh{\mb M}} \times \hat{\mu}, H=G$, and $M=\mo(1)$ as in the previous section, then we get $\mc Z=\mc M$ the enhanced master space.
Since $\Hat{\mu}$ is $\tilde{T} \times \C^{\ast}_{\hbar}$-equivariant, and $\Hat{\mu}^{-1}(0)$ is a complete intersection by \cite[Theorem 1.2]{C}, we have a $\tilde{T} \times \C^{\ast}_{\hbar}$-equivariant perfect obstruction thoery $ob_{\mc M} = \id_{L_{\mc M}} \colon Ob_{\mc M} = L_{\mc M}$ described by \eqref{cot}.

%%%%%%%%%%%%%%%%%%%%%%%%%%%%%%%%%%%%%%%%%%%%%%%%%%%%%%%%%%%%%%%%%%%%%%%%%%%

\subsection{Localization formula}
\label{subsec:loca1}
In the following, we consider a Kuranishi chart 
$$
\left( Y=\wh{\mb M}(r,n), E=\wh{\mb M} \times \mb L, e=\id_{\wh{\mb M}} \times \hat{\mu}, H=G, M= \mo(1) \right)
$$ 
of the enhanced master space $\mc Z=\mc M$.
We also have $\tilde{T} \times \C^{\ast}_{\hbar}$-action on $(Y, E)$ compatible with $G$-action such that $e$ is also $\tilde{T} \times \C^{\ast}_{\hbar}$-invariant. 
We recall the localization formula in \cite{GP} in this setting.

For $\alpha=\pm$, or $\mk J \in S^{\ell}$, as in \S \ref{subsec:fixe} we take closed subsets $\wh{\mb M}_{\alpha}$ in $\wh{\mb M}$ and subgroups $G_{\alpha}$ of $G$ such that $[\wh{\mb M}_{\alpha}^{ss}/G_{\alpha}] \cong \mc N_{\alpha}$ in \eqref{ambient}, where $\wh{\mb M}_{\alpha}^{ss}=\wh{\mb M}_{\alpha} \times_{\wh{\mb M}} \wh{\mb M}^{ss}$.
Furthermore, we can modify $\C^{\ast}_{\hbar}$-actions up to $G_{\alpha}$-actions and taking finite covers of $\C^{\ast}_{\hbar}$ as in \eqref{act} such that $\C^{\ast}_{\hbar}$ trivially acts on $\wh{\mb M}_{\alpha}$. 
Hence we have decompositions $E|_{\wh{\mb M}_{\alpha}}=E_{\alpha} \oplus (E|_{\wh{\mb M}_{\alpha}})^{mov}$ by invariant parts $E_{\alpha}=(E|_{\wh{\mb M}_{\alpha}})^{inv}$ and moving parts $(E|_{\wh{\mb M}_{\alpha}})^{mov}$ with respect to $\C^{\ast}_{\hbar}$-action.
Since the section $e=\id_{\wh{\mb M}} \times \hat{\mu}$ is $\C^{\ast}_{\hbar}$-invariant, we see that $e_{\alpha}=e|_{\wh{\mb M}_{\alpha}}$ belongs to the invariant part $ \Gamma(\wh{\mb M}_{\alpha}, E_{\alpha})$.

The enhanced master space $\mc M$ is the fiber product $\mc N \times_{\mc E} \mc N$ of the zero section and the section $\mc N \to \mc E= [E / G]$ induced by $e$.
Pulling back this by the embedding
$$
\E_{\alpha}=[E_{\alpha} / G_{\alpha}] \to \E, 
$$
we see that $\mc M_{\alpha} =\mc M \times_{\mc N} \mc N_{\alpha}$ is also the zero-locus $\mc N_{\alpha} \times_{\mc E_{\alpha}} \mc N_{\alpha}$ of the section $\mc N_{\alpha} \to \mc E_{\alpha}$ induced by $e_{\alpha}$.  
Hence we see that $(\wh{\mb M}_{\alpha}, E_{\alpha}, e_{\alpha}, G_{\alpha}, M_{\alpha})$ is a Kuranishi chart of $\mc M_{\alpha}$, where $M_{\alpha} = \mo(1)|_{\wh{\mb M}_{\alpha}}$.

For the embedding $\wh{\mb M}_{\alpha} \to \wh{\mb M}$ and the modified $\C^{\ast}_{\hbar}$-action as above,  we have 
%decompositions 
%$$
%E|_{\wh{\mb M}_{\alpha}}=E_{\alpha} \oplus (E|_{\wh{\mb M}_{\alpha}})^{mov}, T \wh{\mb M}|_{\wh{\mb M}_{\alpha}}= T \wh{\mb M}_{\alpha} \oplus (T \wh{\mb M}|_{\wh{\mb M}_{\alpha}})^{mov}, \mk g= \mk g_{\alpha} \oplus (\mk g_{\alpha})^{mov}, 
%$$ and
a decomposition $\iota_{\alpha}^{\ast} Ob_{\mc M}= (\iota_{\alpha}^{\ast} Ob_{\mc M})^{inv} \oplus (\iota_{\alpha}^{\ast} Ob_{\mc M})^{mov}$, where $\iota_{\alpha} \colon \mc M_{\alpha} \to \mc M$ is the embedding.
The invariant part $(\iota_{\alpha}^{\ast} Ob_{\mc M})^{inv}$ coincides with the complex \eqref{cot} constructed from the Kuranishi chart $(\wh{\mb M}_{\alpha}, E_{\alpha}, e_{\alpha}, G_{\alpha}, M_{\alpha})$, and we have obstruction theories
$\iota_{\alpha}^{\ast} ob_{\mc M} \colon (\iota_{\alpha}^{\ast} Ob_{\mc M})^{inv} \to L_{\mc M_{\alpha}}$. 
Hence virtual funcamental cycles $[\mc M_{\alpha}]^{vir}$ in $A^{\tilde{T} \times \C^{\ast}_{\hbar}}_{\ast}(\mc M^{\C^{\ast}_{\hbar}} )$ are defined similarly.
We define {\it virtual normal bundles} $\mk N(\mc M_{\alpha})$ of $\mc M_{\alpha}$ in $\mc M$ by $\mk N(\mc M_{\alpha})=(( \iota_{\alpha}^{\ast} Ob_{\mc M})^{mov})^{\vee}$ in the $K$-group of $\mc M^{\C^{\ast}_{\hbar}}$.

By the main result \cite[(1)]{GP}, we have 
\begin{eqnarray}
\label{local}
[\mc M]^{vir} = \iota_{\ast} \sum_{\alpha=\pm, \mk J} \frac{[\mc M_{\alpha}]^{vir}}{e(\mk N(\mc M_{\alpha}))} \in A^{\C^{\ast}_{\hbar}}_{\ast}( \mc M \times_{\tilde{T}} E_{m}) \otimes \mb Q[\hbar, \hbar^{-1}], 
\end{eqnarray}
where $E_{m} \to E_{m}/\tilde{T}$ is any approximation of the universal bundle $E \tilde{T} \to B \tilde{T}$ over the classifying space, and $\iota \colon \mc M^{\C^{\ast}_{\hbar}} \to \mc M$ is the inclusion.
In our setting, this is shown similarly as in \cite[\S 2]{GP}, since we have intersection theories on algebraic stacks by \cite{Kr} with the similar property as in \cite{F}.

\begin{NB}
We have $[\mc N] = j_{\ast} \sum \frac{[\mc N_{\alpha}]}{e(\mk N(\mc N_{\alpha}))}$ for a smooth DM stack with $\C^{\ast}$-action, where $j \colon \mc N^{\C^{\ast}_{\hbar}} \to \mc N$ is the embedding.
For the diagram
$$
\xymatrix{
\mc N_{\alpha} \ar[r]^{j_{\alpha}} & \mc N \ar[r]^{e} & \E \\
\ar[u] \mc M_{\alpha}^{\iota_{\alpha}} \ar[r] &\ar[u] \mc M \ar[r] & \mc N \ar[u]^{0}},
$$
we have $0^{!} j_{\alpha \ast} = \iota_{\alpha \ast} 0^{!}$, where $0 \colon \mc N \to \E$ is the zero-section.
Furthermore for the diagram
$$
\xymatrix{
\mc N_{\alpha} \ar[r]^{j_{\alpha}} & \E_{\alpha} \ar[r] & E_{\alpha} \oplus (\E|_{\mc N_{\alpha}})^{mov} \\
\ar[u] \mc M_{\alpha}^{\iota_{\alpha}} \ar[r] &\ar[u]^{0_{\alpha}} \mc N_{\alpha} \ar@{=}[r] & \mc N_{\alpha} \ar[u]^{0|_{\mc N_{\alpha}}}},
$$
we have $0^{!}= e((\E|_{\mc N_{\alpha}})^{mov}) \cap 0_{\alpha}^{!}$ by the excess intersection formula, where $0_{\alpha} \colon \mc N_{\alpha} \to \E_{\alpha}$ is the zero-section.
Hence applying the Gysin homomorphism $0^{!}$ to the equation $[\mc N] = j_{\ast} \sum \frac{[\mc N_{\alpha}]}{e(\mk N(\mc N_{\alpha}))}$, we get
\begin{eqnarray*}
[\mc M]^{vir} = \iota_{\ast} \sum \frac{[\mc M_{\alpha}]^{vir}}{e(\mk N(\mc M_{\alpha}))}. 
\end{eqnarray*}.
\end{NB}

In the following computations, we cite corresponding results in \cite[\S 5]{M}, but they are much simpler in our case.

%%%%%%%%%%%%%%%%%%%%%%%%%%%%%%%%%%%%%%%%%%%%%%%%%%%%%%%%%%%%%%%%%%%%%%%%%%%

\subsection{Obstruction theories for $\mc M_{\pm}$}
\label{subsec:obst1}
For $\alpha=\pm$, as Kuranishi charts we can take $\wh{\mb M}_{\pm}=\PP(L_{\pm})=\wt{\mb M}$ embedded into $\wh{\mb M}=\PP(L_{-} \oplus L_{+})$ by $x_{\mp} = 0$, $E_{\pm}=\wh{\mb M}_{\pm} \times \mb L, e_{\pm}=\id_{\wh{\mb M}_{\pm}} \times \tilde{\mu}, G_{\pm}=G$ and $M_{\pm}=L_{\pm}$ as in \S \ref{subsec:modu}.
In our case, it is easy to see that the following lemma holds.
\begin{lem}[\protect{\cite[Proposition 5.9.2]{M}}]
\label{ob1}
We have 
$$
\iota_{\pm}^{\ast} ob_{\mc M} = \id_{L_{\mc M_{\pm}}}, \mk N(\mc M_{\pm}) = \mc L_{\mp} \otimes \mc L_{\pm}^{\vee} \otimes e^{\pm \hbar}, 
$$
where $\mc L_{\pm}$ are line bundles on $\mc M_{\pm}$ defined from $L_{\pm}$ in \S \ref{subsec:modu}.
\\
\end{lem}
\proof
Since $\tilde{\mu}^{-1}(0)$ is a complete intersection, the complex \eqref{cot} is quasi-isomorphic to the cotangent complex $L_{\mc M_{\pm}}$.
Hence we can identify $\iota_{\pm}^{\ast} ob_{\mc M}$ with $\id_{L_{\mc M_{\pm}}}$.

For $\mk N(\mc M_{\pm})$, we consider the original $\C^{\ast}_{\hbar}$-action \eqref{act0}.
Then moving parts $(\iota_{\pm}^{\ast} Ob_{\mc M})^{mov}$ are $L_{\mp} \otimes L_{\pm}^{\vee} \otimes e^{\pm \hbar}|_{e_{\pm}^{-1}(0)}$ as $G$-equivariant bundles on $e_{\pm}^{-1}(0)$.
\endproof

%%%%%%%%%%%%%%%%%%%%%%%%%%%%%%%%%%%%%%%%%%%%%%%%%%%%%%%%%%%%%%%%%%%%%%%%%%%

\subsection{Obstruction theories for $\mc M_{\mk J}$}
\label{subsec:obst2}
We consider a decomposition type $\mk J=(I_{\flat}, I_{\sharp}) \in S^{\ell}$, and fix a decomposition $V=V_{\flat} \oplus V_{\sharp}$ with $\dim V_{\sharp}=|I_{\sharp}|=p$. 
Then we use the identification \eqref{mkj} to describe $\iota_{\mk J}^{\ast} ob_{\mc M}$ as follows.
 
We consider a Kuranishi chart $(\wh{\mb M}_{\mk J}, E_{\mk J}, e_{\mk J}, G_{\mk J}, M_{\mk J})$ of $\mc M_{\mk J}$ in \S \ref{subsec:loca1}, where 
$$
\wh{\mb M}_{\mk J}= \left( \wt{\mb M}(W, V_{\flat}) \times \wt{\mb M}(0, V_{\sharp}) \right) \times_{\wt{\mb M}(r,n)} \wh{\mb M}(r,n), G_{\mk J}=\GL(V_{\flat}) \times \GL(V_{\sharp})
$$ 
as in \S \ref{subsec:fixe}.
Furthermore we have $E_{\mk J}= \wh{\mb M}_{\mk J} \times \mb L(V_{\flat}) \times \mb L(V_{\sharp})$, $e_{\mk J}= \id_{\wh{\mb M}_{\mk J}} \times \varphi$,  where
$$
\varphi  \colon \wh{\mb M}_{\mk J} \to \mb L(V_{\flat}) \times \mb L(V_{\sharp})
$$
is the composition of the projection to $\wt{\mb M}(W,V_{\flat}) \times \wt{\mb M}(0, V_{\sharp})$ and $\tilde{\mu} \times \tilde{\mu}$.
By lemma \ref{lem4}, $M_{\mk J}=\mo(1)|_{\wh{\mb M}_{\mk J}}$ gives the product of $(\min(I_{\sharp}) -1)$-stability and $+$-stability.
As explained in \S \ref{subsec:loca1}, $\iota_{\mk J}^{\ast} ob_{\mc M}$ coincides with one constructed from the Kuranishi chart $(\wh{\mb M}_{\mk J}, E_{\mk J}, e_{\mk J}, G_{\mk J}, M_{\mk J})$.

More explicitly we have $\wh{\mb M}_{\mk J}^{ss} \subset \wt{\mb M}(W,V_{\flat}) \times \wt{\mb M}(0, V_{\sharp}) \times \C^{\ast}_{\rho}$, where $\rho = \frac{x_{-}}{x_{+}}$.
To get a Kuranishi chart $(\wh{\mb M}_{\mk J}', E_{\mk J}', e_{\mk J}', G_{\mk J}', M_{\mk J}')$ of $\left( \bar{M} \times M^{+}_{p} \right)/ \C^{\ast}_{t}$, we put 
$$
\wh{\mb M}_{\mk J}'=\wt{\mb M}(W,V_{\flat}) \times \wt{\mb M}(0, V_{\sharp}) \times \C^{\ast}_{\rho_{\flat}} \times \C^{\ast}_{\rho_{\sharp}}, G_{\mk J}'=G_{\mk J} \times \C^{\ast}_{t}.
$$
We have a natural map $\wh{\mb M}_{\mk J}' \to \wh{\mb M}_{\mk J}$ induced by $\id_{\wt{\mb M}(W,V_{\flat})} \times \id_{\wt{\mb M}(0, V_{\sharp})}$, multiplication $(\rho_{\flat}, \rho_{\sharp}) \to \rho=\rho_{\flat} \rho_{\sharp}$, and inclusion $\wt{\mb M}(W,V_{\flat}) \times \wt{\mb M}(0, V_{\sharp}) \times \C^{\ast}_{\rho} \to \wh{\mb M}(r,n)$ by $\rho = \frac{x_{-}}{x_{+}}$.
Other data $E_{\mk J}', e_{\mk J}', M_{\mk J}'$ are obtained by pull-back of $E_{\mk J}, e_{\mk J}, M_{\mk J}$ by this map.

We write by $ob_{\mc M_{\mk J}}$ the obstruction theory constructed from this Kuranishi chart $(\wh{\mb M}_{\mk J}', E_{\mk J}', e_{\mk J}', G_{\mk J}', M_{\mk J}').$ 
The isomorphism \eqref{mkj} is obtained from the coordinate change of Kuranishi charts $(\wh{\mb M}_{\mk J}, E_{\mk J}, e_{\mk J}, G_{\mk J}, M_{\mk J})$ and $(\wh{\mb M}_{\mk J}', E_{\mk J}', e_{\mk J}', G_{\mk J}', M_{\mk J}')$ induced by the above map $\wh{\mb M}_{\mk J}' \to \wh{\mb M}_{\mk J}$ and the projection $G_{\mk J}' \to G_{\mk J}$. 
Hence as we remarked in \S \ref{subsec:virt}, we have the same obstruction theory $\iota_{\mk J}^{\ast} ob_{\mc M} = ob_{\mc M_{\mk J}}$ via the isomorphism \eqref{mkj}.

%%%%%%%%%%%%%%%%%%%%%%%%%%%%%%%%%%%%%%%%%%%%%%%%%%%

\subsection{Obstruction theories for decompositions}
\label{subsec:obst3}
For a fixed $\mk J = (I_{\flat}, I_{\sharp}) \in S^{\ell}$ and a decomposition $V=V_{\flat} \oplus V_{\sharp}$, we consider \'{e}tale morphisms 
$$
F \colon \mc S_{\mk J} \to \mc M_{\mk J}, G \colon \mc S_{\mk J} \to \widetilde{M}^{\min(I_{\sharp})-1}(r,n-p) \times M^{+}_{p}
$$
in Theorem \ref{decomp}.
We have a Kuranishi chart of $\mc S_{\mk J}$ by replacing $\C^{\ast}_{t}$ in $G_{\mk J}' = G_{\mk J} \times \C^{\ast}_{t}$ with $\C^{\ast}_{s}$, where $s=t^{\frac{1}{pD}}$, and $F \colon \mc S_{\mk J} \to \mc M_{\mk J}$ is defined from the coordinate change of Kuranishi charts induced by a map $s \mapsto t=s^{pD}$.
Hence we have $ob_{\mc S_{\mk J}}=F^{\ast} ob_{\mc M_{\mk J}}$.

On the other hand, we have the isomorphism $\mc S_{\mk J} \cong \ba{M} / \C^{\ast}_{s} \times M_{p}^{+}$ in the proof of Theorem \ref{decomp}.
This is also defined from the coordinate change of Kuranishi charts induced from the identity of bundles and the automorphism \eqref{gphom} of the group $\GL(V_{\flat}) \times \GL(V_{\sharp}) \times \C^{\ast}_{s}$.
Since the morphism $G$ is obtained by the isomorphism $\mc S_{\mk J} \cong \ba{M} / \C^{\ast}_{s} \times M_{p}^{+}$ and a map $s \mapsto u=s^{pD}$, this is also induced from a coordinate change of Kuranishi charts.
Hence we have 
$$
ob_{\mc S_{\mk J}}=G^{\ast} \left( p_{1}^{\ast} ob_{\wt{M}^{\min(I_{\sharp})-1}(r,n-p)} \oplus p_{2}^{\ast} ob_{M_{p}^{+}} \right),
$$ 
where $p_{1}, p_{2}$ are projections from $\wt{M}^{\min(I_{\sharp})-1}(r,n-p) \times M_{p}^{+}$ to 
$$
\wt{M}^{\min(I_{\sharp})-1}(r,n-p), M_{p}^{+}.
$$

%%%%%%%%%%%%%%%%%%%%%%%%%%%%%%%%%%%%%%%%%%%%%%%%%%%%%%%%%%%%%%%%%%%%%%%%%%%

\subsection{Normal bundle $\mk N(\mc M_{\mk J})$}
\label{subsec:norm}
In the following, for vector bundles $\mc E, \mc F$ on a stack $\mc Z$, we write by $\mc Hom(\mc E, \mc F)$ the vector bundle $\mc E^{\vee} \otimes \mc F$ on $\mc Z$.
We write by $\mc W$ the tautological bundle corresponding to $W$ on various moduli stacks $\mc M$, $\wt{M}^{\min(I_{\sharp})-1}(r,n-p)$, and so on.

We consider $\C^{\ast}_{\frac{\hbar}{pD}}$-action on $\mc M$ as in \eqref{act}.
Then we have 
\begin{eqnarray*}
\label{virtual}
\mk N(\mc M_{\mk J})
&=&
\sum_{(\heartsuit, \spadesuit)} \mc Hom( Q^{\vee} \otimes \mc V^{\mk J}_{\heartsuit}, \mc V^{\mk J}_{\spadesuit}) - \mc Hom(\wedge Q^{\vee} \otimes \mc V^{\mk J}_{\heartsuit}, \mc V^{\mk J}_{\spadesuit}) - \mc Hom( \mc V^{\mk J}_{\heartsuit}, \mc V^{\mk J}_{\spadesuit}) \notag \\
&+&
\mc Hom( \mc W, \mc V^{\mk J}_{\sharp}) +  \mc Hom(\wedge Q^{\vee} \otimes \mc V^{\mk J}_{\sharp}, \mc W)+N_{0},
\end{eqnarray*}
where the sum is taken for $(\heartsuit, \spadesuit)=(\flat, \sharp), (\sharp, \flat)$, and $N_{0}$ corresponds to the normal bundle for the embedding $Fl_{\mk J} \colon F(V_{\flat}, \underbar{n-p}) \times F(V_{\sharp}, \underbar{p}) \to F(V, \underbar{n})$.
If we consider the universal full flag $\mc F^{\bullet}$ of $\mc V$ on $\mc M$, then in the $K$-group, we have
$$
N_{0} =\bigoplus_{(\heartsuit, \spadesuit)} \bigoplus_{\substack{i > j \\ i \in I_{\heartsuit}, j\in I_{\spadesuit}}}  \mc Hom\left(\mc F^{j} / \mc F^{j-1}, \mc F^{i}/\mc F^{i-1} \right) \Big |_{\mc M_{\mk J}},
$$ 
where the sum is also taken for $(\heartsuit, \spadesuit)=(\flat, \sharp), (\sharp, \flat)$.
 
On the other hand, on $\wt{M}^{\min(I_{\sharp})-1}(r,n-p) \times M^{+}_{p}$ we put
\begin{eqnarray}
\mk N_{1} 
&=& 
\label{n1}
\mc Hom( Q^{\vee} \otimes \mc V_{\flat}, \mc V_{\sharp}) - \mc Hom(\wedge Q^{\vee} \otimes \mc V_{\flat}, \mc V_{\sharp}) -
\mc Hom( \mc V_{\flat}, \mc V_{\sharp}) \\
&+& 
\mc Hom( \mc W, \mc V_{\sharp}), \notag \\
\mk N_{2}
&=&
\mc Hom(Q^{\vee} \otimes \mc V_{\sharp}, \mc V_{\flat}) - \mc Hom(\wedge Q^{\vee} \otimes \mc V_{\sharp}, \mc V_{\flat}) - \mc Hom(\mc V_{\sharp}, \mc V_{\flat}) \\
&+& 
\label{n2} 
\mc Hom(\wedge Q^{\vee} \otimes \mc V_{\sharp}, \mc W). \notag
\end{eqnarray}
We consider universal full flags $\mc F_{\flat}^{\bullet}, \mc F_{\sharp}^{\bullet}$ of tautological bundles $\mc V_{\flat}, \mc V_{\sharp}$ on moduli stacks $\wt M^{\min(I_{\sharp})-1}(r,n-p), M_{p}^{+}$.
By the decomposition data $\mk J=(I_{\flat}, I_{\sharp})$, these are considered as full flags with repetitions indexed by $\underbar{n}$.
We also write their pull-backs to the product $\wt{M}^{\min(I_{\sharp})-1}(r,n-p) \times M^{+}_{p}$
by the same letter $\mc F_{\flat}^{\bullet}, \mc F_{\sharp}^{\bullet}$.

By computations in the previous subsections and Theorem \ref{decomp}, we have the following proposition, which corresponds to \cite[Proposition 5.8.1 and remark after Proposition 5.9.3]{M}.
\begin{prop}
\label{ob2}
We have the following.
\begin{enumerate}
\item[\textup{(1)}]
We have $F^{\ast} \iota_{\mk J}^{\ast} ob_{\mc M} = G^{\ast} \left( p_{1}^{\ast} ob_{\wt{M}^{\min(I_{\sharp})-1}(r,n-p)} \oplus p_{2}^{\ast} ob_{M_{p}^{+}} \right).$
\item[\textup{(2)}]
We have 
$
F^{\ast} \mk N(\mc M_{\mk J}) 
=
G^{\ast} \mk N_{1} \otimes e^{\frac{\hbar}{pD}} \otimes L_{\mc S}^{\vee} +G^{\ast} \mk N_{2} \otimes e^{\frac{-\hbar}{pD}} \otimes L_{\mc S} 
+
G^{\ast} \Theta^{rel}_{0}
$
in the equivariant $K$-group of $S_{\mk J}$, where
$$
\Theta^{rel}_{0} = \sum_{(\heartsuit, \spadesuit)} \sum_{i > j}  \mc Hom\left( \mc F^{j}_{\heartsuit}/\mc F^{j-1}_{\heartsuit}, \mc F^{i}_{\spadesuit}/\mc F^{i-1}_{\spadesuit} \right),
$$ 
and the sum is also taken for $(\heartsuit, \spadesuit)=(\flat, \sharp), (\sharp, \flat)$.
\end{enumerate}
\end{prop}

%%%%%%%%%%%%%%%%%%%%%%%%%%%%%%%%%%%%%%%%%%%%%%%%%%%%%%%%%%%%%%%%%%%%%%%%%%%
%%%%%%%%%%%%%%%%%%%%%%%%%%%%%%%%%%%%%%%%%%%%%%%%%%%%%%%%%%%%%%%%%%%%%%%%%%%
%%%%%%%%%%%%%%%%%%%%%%%%%%%%%%%%%%%%%%%%%%%%%%%%%%%%%%%%%%%%%%%%%%%%%%%%%%%

\section{Hilbert schemes parametrizing destabilizing objects}
\label{sec:hilb}
In this section, we reduce integrals on $M_{p}^{+}$ introduced in \S \ref{subsec:deco} to computations on Hilbert schemes.

\subsection{Morphisms to Hilbert schemes}
\label{subsec:morp}
Here, we construct a morphism $h \colon M_{p}^{+} \to M(1,p)$.
We consider a Kuranishi chart $(Y, E, e, H, M)$ of $M_{p}^{+}$, where
$
Y= \wt{\mb M}(0, V_{\sharp}) \times \C^{\ast}_{\rho_{\sharp}}, E=Y \times \mb L(V_{\sharp}), e=\id_{Y} \times \varphi, H=\GL(V_{\sharp}),
$
and $M$ corresponds to $+$-stability condition.
We can take such $M$ as in Proposition \ref{ell} taking $\eta_{1}$ enough large.
Here $\varphi \colon Y \to \mb L(V_{\sharp})$ is the composition of the projection to $\wt{\mb M}(0, V_{\sharp})$ and $\tilde{\mu}$.
The obstruction theory $ob_{M_{p}^{+}}$ in Proposition \ref{ob2} is defined from this Kuranishi chart as described in \S \ref{subsec:obst3}.

On the other hand, we introduce another vector space $W_{\sharp}=\C$, and consider the natural projection
$$
Fl(V_{\sharp}, \underbar{p}) \to \mb P(V_{\sharp}) \cong \left[ \frac{\Hom_{\C}(W_{\sharp}, V_{\sharp}) \setminus \lbrace 0 \rbrace }{\C^{\ast}_{u}} \right], F_{\sharp}^{\bullet} \mapsto F_{\sharp}^{1}. 
$$
If we put $U=\mb M(0, V_{\sharp}) \times \left( \Hom_{\C}(W_{\sharp}, V_{\sharp}) \setminus \lbrace 0 \rbrace \right) \times \C^{\ast}_{\rho_{\sharp}}$, then $z_{\sharp} \in \Hom_{\C}(W_{\sharp}, V_{\sharp}) \setminus \lbrace 0 \rbrace$ gives an inclusion $\mo_{U}=\mo_{U}\otimes W_{\sharp} \to \mo_{U} \otimes V_{\sharp}$ of vector bundles on $U$, and the quotient $\mo_{U} \otimes V_{\sharp} / \mo_{U}$.
The ambient space $Y=\wt{\mb M}(0, V_{\sharp}) \times \C^{\ast}_{\rho_{\sharp}}$ is the quotient of the full flag bundle $F(\mo_{U} \otimes V_{\sharp}/\mo_{U}, \underbar{p-1})$ over $U$ by $\C^{\ast}_{u}$. 
Hence, we get a coordinate change $(Y_{1}, E_{1}, e_{1}, H_{1}, M_{1}) \to (Y, E, e, H, M)$ induced by projections $Y_{1} \to Y$ and $H_{1} \to H$, where 
$$
Y_{1}=Fl(\mo_{U} \otimes V_{\sharp}/\mo_{U}, \underbar{p-1}), E_{1}=E|_{Y_{1}}, e_{1}=e|_{Y_{1}}, H_{1}=\GL(V_{\sharp}) \times \C^{\ast}_{u}
$$
$M_{1}$ is the pull-back of $M$. 

We also have a morphism $(Y_{1}, E_{1}, e_{1}, H_{1}, M_{1}) \to (\bar{Y}_{1}, \bar{E}_{1}, \bar{e}_{1}, \bar{H}_{1}, \bar{M}_{1})$, where
$$
\bar{Y}_{1} = \mb M(0, V_{\sharp}) \times \Hom(W_{\sharp}, V_{\sharp}), \bar{E}_{1}= \bar{Y}_{1} \times \mb L(V_{\sharp}), \bar{e}_{1}=\id_{\bar{Y}_{1}} \times \bar{\varphi}, \bar{H}_{1}=\GL(V_{\sharp}),
$$
induced by the projection $Y_{1} \to \bar{Y}_{1}$ and a map $H_{1} \to \bar{H}_{1}, (g_{\sharp}, u) \mapsto g_{\sharp} u$.
Here, $\bar{\varphi}$ is the restriction of $\mu$ to the subset $\bar{Y}_{1} = \mb M(0, V_{\sharp}) \times \Hom(W_{\sharp}, V_{\sharp}) \subset \mb M(W_{\sharp}, V_{\sharp}).$
We also consider the restriction of the stability condition on $\mb M(W_{\sharp}, V_{\sharp})$ defined in Definition \ref{stab}. 

Then $(\bar{Y}_{1}, \bar{E}_{1}, \bar{e}_{1}, \bar{H}_{1}, \bar{M}_{1})$ is a Kuranishi chart of the moduli $M(1,p)$ of ADHM data on $(W_{\sharp}, V_{\sharp})$, which is different from the usual one using $\mb M(W_{\sharp}, V_{\sharp})$ instead of $\bar{Y}_{1}$.
We have a morphism $h \colon M_{p}^{+} \to M(1,p)$ induced from the above morphism of Kuranishi charts.
%We consider the composition $\mu_{\sharp} \colon \mb{M}_{p}^{+} \to \mb{L}(V_{\sharp})$ of the projection $\mb{M}_{p}^{+} \to \mb{M}(0, V_{\sharp})$ and $\mu \colon \mb{M}(0, V_{\sharp}) \to \mb{L}(V_{\sharp})$.
%We have $\GL(V_{\sharp}) \times \C^{\ast}_{u}$-action on $\mb{M}_{p}^{+}$ defined by 
%$$
%g \cdot (X_{\sharp}, z_{\sharp}, \ba{F}^{\bullet}_{\sharp},\rho_{\sharp})=(g X_{\sharp}, ug z_{\sharp}, g \ba{R}^{\bullet}, (\det g)^{-D}\rho_{\sharp})
%$$ 
%for $g \in \GL(V_{\sharp}), u \in \C^{\ast}_{u}$, where $\ba{F}^{\bullet}_{\sharp}$ is a full flag of $V_{\sharp} / \langle z_{\sharp} \rangle$. 
%
%For the $+$-stable locus $\mu_{\sharp}^{-1}(0)^{s}$, we have an isomorphism
%$
%[\mu_{\sharp}^{-1}(0)^{s}/ \GL(V_{\sharp}) \times \C^{\ast}_{u}] \cong M_{p}^{+}
%$
%induced from a isomorphism 
%\begin{eqnarray}
%\label{morph}
%\mb M_{p}^{+} \to \wt{\mb M}(0, V_{\sharp}) \times \C^{\ast}_{\rho_{\sharp}}, (X_{\sharp}, z_{\sharp}, \ba{F}^{\bullet}_{\sharp},\rho_{\sharp}) \to (X_{\sharp}, F^{\bullet}_{\sharp}, \rho_{\sharp})
%\end{eqnarray}
%and projections $\GL(V_{\sharp}) \times \C^{\ast}_{u} \to \GL(V_{\sharp})$.
%Here $F_{\sharp}^{\bullet}$ is defined by $F^{1}_{\sharp} = \langle z_{\sharp} \rangle$, and $F^{i}_{\sharp} = \pi_{z_{\sharp}}^{-1} ( \ba{F}^{i+1}_{\sharp} )$ for  $i>1$, where $\pi_{z_{\sharp}} \colon V_{\sharp} \to V_{\sharp} / \langle z_{\sharp} \rangle$ is the natural surjection, and $z_{\sharp} \in \Hom_{\C}(W_{\sharp}, V_{\sharp})$ is identified with the vector $z_{\sharp}(1) \in V_{\sharp}$ for $1 \in W_{\sharp}=\C$.
%
We consider a line bundle $\left( \det \mc V_{\sharp} \right)^{D}$ on $M(1,p)$, where $\mc V_{\sharp}$ is the tautological bundle on $M(1,p)$ corresponding to $V_{\sharp}$.

\begin{prop}
[\protect{\cite[Proposition 5.9]{NY3}}]
\label{hilb}
We have the following.
\begin{enumerate}
\item[\textup{(1)}] $M^{+}_{p}$ is a full flag bundle $Flag(\mc V_{\sharp}/\mo_{M(1,p)}, \underline{p-1})$ over the quotient stack 
$$
\left[ \left( \left( \det \mc V_{\sharp} \right) ^{-D} \right) ^{\times} / \C^{\ast}_{u} \right],
$$ 
where $\C^{\ast}_{u}$ acts on $\left( \left( \det \mc V_{\sharp} \right) ^{-D} \right) ^{\times}$ by fiber-wise multiplication of $u^{pD}$.

\item[\textup{(2)}] The homomorphism $\C^{\ast}_{u} \to \C^{\ast}_{s}$ given by $s=u^{pD}$ induces an \'{e}tale and finite morphism $\left[ \left( \left( \det \mc V_{\sharp} \right) ^{-D} \right) ^{\times} / \C^{\ast}_{u} \right] \to M(1,p)$.
\end{enumerate}
\end{prop}
\proof 
(1) For any element $(X_{\sharp}, F^{\bullet}_{\sharp}, \rho_{\sharp}) \in M^{+}_{p}$,  the $+$-stability is equivalent to the condition that for any non-zero element $z_{\sharp} \in F^{1}_{\sharp}$, a tuple $(X_{\sharp}, z_{\sharp}, 0)$ is a stable ADHM data on $(W_{\sharp}, V_{\sharp})$.
Hence by forgetting $F^{i}_{\sharp}$ for $i > 1$, we get elements in the quotient stack of 
$$ 
\lbrace (X_{\sharp}, z_{\sharp}, 0, \rho_{\sharp} ) \in \mb M(W_{\sharp}, V_{\sharp} ) \times \C^{\ast}_{\rho_{\sharp}} \mid (X_{\sharp}, z_{\sharp}, 0) \text{ is a stable ADHM data }\rbrace 
$$
by the group $\GL(V_{\sharp}) \times \C^{\ast}_{u}$-action, where $(g_{\sharp},u) \in \GL(V_{\sharp}) \times \C^{\ast}_{u}$ acts by 
$$
\left( g_{\sharp} X_{\sharp} \id_{Q^{\vee}} \otimes g_{\sharp}^{-1}, ug_{\sharp} z_{\sharp}, 0, (\det g_{\sharp} )^{-D} \rho_{\sharp} \right).
$$
Then this is isomorphic to the quotient stack 
$$
\left[ \left( \left( \det \mc V_{\sharp} \right) ^{-D} \right) ^{\times} / \C^{\ast}_{u} \right]
$$ 
by a group homomorphism $\GL(V_{\sharp}) \times \C^{\ast}_{u} \to \GL(V_{\sharp}) \times \C^{\ast}_{u}$ defined by $(g_{\sharp}, u) \mapsto (ug_{\sharp}, u)$.

(2) It follows from (1), since $\left[ \left( \left( \det \mc V_{\sharp} \right) ^{-D} \right) ^{\times} / \C^{\ast}_{s} \right]=M(1,p)$.
\endproof
From the proof, we easily see that $h \colon M^{+}_{p} \to M(1,p)$ is the composition of the morphisms $M^{+}_{p} \to \left[ \left( \left( \det \mc V_{\sharp} \right) ^{D} \right) ^{\times} / \C^{\ast}_{u} \right]$ and $\left[ \left( \left( \det \mc V_{\sharp} \right) ^{D} \right) ^{\times} / \C^{\ast}_{u} \right] \to M(1,p)$ in the above proposition.
We have an isomorphism
\begin{eqnarray}
\label{vp+}
h^{\ast}\mc V_{\sharp} \cong \mc V_{p}^{+} \otimes L_{M_{p}^{+}}, 
\end{eqnarray}
where $L_{M_{p}^{+}}$ is a line bundle on $M_{p}^{+}$ such that $L_{M_{p}^{+}}^{\otimes pD} \cong h^{\ast} \det \mc V_{\sharp}^{-D}$.

We write by $ob_{M(1,p)}^{\sharp}$ the obstruction theory on $M(1,p)$ defined from the Kuranishi chart $(\bar{Y}_{1}, \bar{E}_{1}, \bar{e}_{1}, \bar{H}_{1}, \bar{M}_{1})$.

\subsection{Obstruction theories $ob_{M(1,p)}^{\sharp}$}
\label{subsec:obst4}

For later computations, we compare obstruction theories $ob_{M(1,p)}^{\sharp}$ and $ob_{M(1,p)}$ defined from Kuranishi charts $(\bar{Y}_{1}, \bar{E}_{1}, \bar{e}_{1}, \bar{H}_{1}, \bar{M}_{1})$ in the previous section and the usual one $(\bar{Y}_{0}, \bar{E}_{0}, \bar{e}_{0}, \bar{H}_{0}, \bar{M}_{0})$, where 
$$
\bar{Y}_{0} = \mb M(W_{\sharp}, V_{\sharp}), \bar{E}_{0} = \bar{Y}_{0} \times L(V_{\sharp}), \bar{e}_{0} = \id_{\bar{Y}_{0}} \times \mu, \bar{H}_{0} = \GL(V_{\sharp}),
$$
and $\bar{M}_{0}$ is the stability condition on $\mb M(W_{\sharp}, V_{\sharp})$ in Definition \ref{stab}.

The complex \eqref{cot} for $(\bar{Y}_{0}, \bar{E}_{0}, \bar{e}_{0}, \bar{H}_{0}, \bar{M}_{0})$ is quasi-isomorphic to the dual of the tangent bundle $TM(1,p)$.
Since we have an isomorphism
$$
\mb M(W_{\sharp}, V_{\sharp}) /\left( \mb M(0, V_{\sharp}) \times \Hom(W_{\sharp}, V_{\sharp}) \right) \cong \Hom(\wedge Q^{\vee} \otimes V_{\sharp}, W_{\sharp})
$$ 
as vector spaces, we have the following distinguished triangle on $M(1,p)$:
\begin{eqnarray}
\label{sharp}
0 \to \left( Ob_{M(1,p)}^{\sharp} \right)^{\vee} \to TM(1,p) \to \mc Hom(\wedge Q^{\vee} \otimes \mc V_{\sharp}, \mc W_{\sharp}) \to 0,
\end{eqnarray}
where $\mc W_{\sharp} = \mo_{M(1,p)}$ is the universal bundle on $M(1,p)=M(W_{\sharp}, V_{\sharp})$ corresponding to $W_{\sharp}$.

%%%%%%%%%%%%%%%%%%%%%%%%%%%%%%%%%%%%%%%%%%%%%%%%%%%%%%%%%%%%%%%%%%%%%%%%%%%

\subsection{Relative tangent bundles for flags}
 \label{subsec:rela2}
We consider the pull-back $\Theta^{rel}$ to $\mc M$ of the relative tangent bundle of $[\wt{\mu}^{-1}(0)/G]$ over $[\mu^{-1}(0)/G]$. 
After restricting to $\mc M_{\mk J}$,  we have a decomposition 
\begin{eqnarray}
\label{rel}
F^{\ast} \Theta^{rel} = G^{\ast} ( p_{1}^{\ast} \Theta^{rel}_{\flat} \oplus p_{2}^{\ast} \Theta^{rel}_{\sharp} \oplus \Theta^{rel}_{0} )
\end{eqnarray} 
on $\mc S_{\mk J}$, where $p_{1}$ and $p_{2}$ are the projections from $\wt{M}^{\min(I_{\sharp})-1}(r,n-p) \times M_{p}^{+}$ to the first and second components respectively.
Here, in the $K$-group, $\Theta^{rel}_{\alpha}$ for $\alpha=\flat, \sharp$ is equal to $\sum_{i > j}  \mc Hom\left( \mc F^{j}_{\alpha}/\mc F^{j-1}_{\alpha}, \mc F^{i}_{\alpha}/\mc F^{i-1}_{\alpha} \right)$, and
$$
\Theta^{rel}_{0} = \sum_{i > j}  \mc Hom\left( \mc F^{j}_{\sharp}/\mc F^{j-1}_{\sharp}, \mc F^{i}_{\flat}/\mc F^{i-1}_{\flat} \right) \oplus \mc Hom\left( \mc F^{j}_{\flat}/\mc F^{j -1}_{\flat},  \mc F^{i}_{\sharp}/\mc F^{i -1}_{\sharp} \right) 
$$ 
as in Proposition \ref{ob2}.

%$=\bigoplus_{i > j} \left( \Hom_{\C}\left(F^{j}_{\sharp}/F^{j-1}_{\sharp}, F^{i}_{\flat}/F^{i-1}_{\flat}\right) \oplus \Hom_{\C}\left(F^{j}_{\flat}/F^{j -1}_{\flat}, F^{i}_{\sharp}/F^{i -1}_{\sharp}\right) \right)$.

We also consider the relative tangent bundle $\Theta'$ of $h \colon M^{+}_{p} \to M(1, p)$.
Then we have an exact sequence
$$
0 \to \Theta' \to \Theta_{\sharp} \to h^{\ast} \left( \mc V_{\sharp} / \mo_{M(1,p)} \right) \to 0,
$$
where $\mc V_{\sharp} / \mo_{M(1,p)}$ is the quotient by the tautological homomorphism $\mo_{M(1,p)} \to \mc V_{\sharp}$.
The fibre of $\mc V_{\sharp} / \mo_{M(1,p)}$ over $h(\wh{X}_{\sharp})$ is identified with $\Hom_{\C}(F_{\sharp}^{1}, V_{\sharp}/ F_{\sharp}^{1})$, and the above homomorphism $\Theta_{\sharp} \to h^{\ast} \left( \mc V_{\sharp} / \mo_{M(1,p)} \right)$ is naturally defined.
By the constructions of $Ob_{M_{p}^{+}}$ and $Ob^{\sharp}_{M(1,p)}$、 we have
\begin{eqnarray}
\label{ob4}
\Theta_{\sharp} = \Theta' + h^{\ast} \left( \mc V_{\sharp} / \mo_{M(1,p)} \right)
\end{eqnarray}
in the $K$-group of $M_{p}^{+}$.

\begin{NB}
For any full flag bundle $p \colon Y =Fl(\E, \underbar{n}) \to X$ for a rank $n$ vector bundle $\E$ on $X$, we can check $p_{\ast} \left( e(T_{Y/X}) \cap p^{\ast} \alpha \right) = n! \alpha$ for $\alpha \in A_{\ast}(X)$ as follows.
Since full flag bundles are obtained by iterations of projective bundles, it is enough to show that for $\pi \colon Z =\PP(\E, \underbar{n}) \to X$, we have $\pi_{\ast} \left( e(T_{Z/X}) \cap \pi^{\ast} \alpha \right)  = n \alpha$ for $\alpha \in A_{\ast}(X)$.
This follows from the description of $\pi_{\ast} \colon A_{\ast}(Z) \to A_{\ast}(X)$ in \cite{F}.
\end{NB}

%%%%%%%%%%%%%%%%%%%%%%%%%%%%%%%%%%%%%%%%%%%%%%%%%%%%%%%%%%%%%%%%%%%%%%%%%%%

\subsection{Integrations on Hilbert schemes}
\label{subsec:inte}
Here, we compute integrations on Hilbert schemes $M(1,n)$ for the next section.
We do not consider the torus action on $M(1,n)$ by multiplication of $e^{a_{1}}$ with framings defined in \S \ref{subsec:toru}.
Hence in the following proposition, integrations take values in $\mb Q(\e_{1}, \e_{2}, \phi_{1}, \phi_{2})$. 

\begin{prop}
\label{hilb2}
We have
$$
\sum_{n=0}^{\infty} q^{n} \int_{M(1,n)} e(\mc V \otimes e^{\phi_{1}} \oplus \mc V^{\vee} \otimes e^{\phi_{2}})=(1-q)^{-\frac{\phi_{1}\phi_{2}}{\e_{1}\e_{2}}},
$$
that is, $\int_{M(1,n)} e(\mc V \otimes e^{\phi_{1}} \oplus \mc V^{\vee} \otimes e^{\phi_{2}})=\frac{\prod_{i=1}^{n} (\frac{\phi_{1}\phi_{2}}{\e_{1}\e_{2}} + i -1 )}{n!}$. 
\end{prop}
\proof
We have
\begin{eqnarray*}
&&
\int_{M(1,n)} e(\mc V \otimes e^{\phi_{1}} \oplus \mc V^{\vee} \otimes e^{\phi_{2}})\\
&=&
\sum_{|Y|=n} \frac{\iota_{Y}^{\ast}e(\mc V \otimes e^{\phi_{1}}) \iota_{Y}^{\ast}e(\mc V^{\vee} \otimes e^{\phi_{2}})}{\iota_{Y}^{\ast}e(TM(1,n))}\\
&=&
\lim_{t \to 0} \sum_{|Y|=n} \frac{\ch \iota_{Y}^{\ast}\wedge_{-e^{\phi_{1}t}} \left( \mc V^{\vee} \right) \ch \iota_{Y}^{\ast} \wedge_{-e^{\phi_{2}t}} \left( \mc V  \right)}{\ch \iota_{Y}^{\ast} \wedge_{-1}(T^{\ast}M(1,n))} \Big |_{ \substack{t_{1}=e^{t\e_{1}} \\ t_{2}=e^{t\e_{2}}}},\\
\end{eqnarray*}
where $\ch$ denotes the Hilbert series (cf. \cite[\S 4]{NY1}), and $\wedge_{u} (E) = \sum u^{i} (\wedge^{i} E)$ in $K_{\tilde{T}}(M(1,p))[u]$.
To compute the sum we follow the notation and the method in \cite{Na2}, where we substitute $t=t_{1}^{-1}, q=t_{2}^{-1}$.
We use plethystic substitution (cf. \cite[1 (i)]{Na2}) by symmetric functions
$$
\Omega=\prod_{i=1}^{\infty} \frac{1}{1- x_{i}}=\exp\left(\sum_{n=1}^{\infty} \frac{p_{n}}{n}\right)=\sum_{n=0}^{\infty} h_{n}, 
$$
where $p_{n}$ and $h_{n}$ are the $n^{\text{th}}$ power and complete symmetric functions respectively.

From \cite[VI (6.11')]{Ma} or \cite[(3.5.20)]{H} we have $\tilde{H}_{\mu}[1-u;t_{1},t_{2}] = \ch \iota_{Y}^{\ast} ( \wedge_{-u} \mc V_{\sharp}^{\vee})$, where $\tilde{H}_{\mu}(x;t_{1},t_{2})$ is the modified Macdonald polynomials (cf. \cite[1 (ii) ]{Na2}), and $\mu$ is the partition corresponding to the Young diagram $Y$.
Hence by the Cauchy formula \cite[(1.7)]{Na2} we have
\begin{eqnarray*}
&&
\sum_{n=0}^{\infty}\sum_{|Y|=n} \frac{\ch \iota_{Y}^{\ast}\wedge_{-u_{1}} \left( \mc V^{\vee} \right) \ch \iota_{Y}^{\ast} \wedge_{-u_{2}} \left( \mc V  \right)}{\ch \iota_{Y}^{\ast} \wedge_{-1}(T^{\ast}M(1,n))} q^{n} \\
&=&
\Omega \left[ \frac{(1-u_{1})(1-u_{2})}{(1-t_{1}^{-1})(1-t_{2}^{-1})} q \right]\\
&=&
\exp \left( \sum_{n=1}^{\infty} \frac{(1-u_{1}^{n})(1-u_{2}^{n})}{(1-t_{1}^{-n})(1-t_{2}^{-n})} \frac{q^{n}}{n} \right).
\end{eqnarray*}
Hence we have
\begin{eqnarray*}
&&
\sum_{n=0}^{\infty} q^{n} \int_{M(1,n)} e(\mc V \otimes e^{\phi_{1}} \oplus \mc V^{\vee} \otimes e^{\phi_{2}}) \\
&=& 
\lim_{t \to 0} \exp \left( \sum_{n=1}^{\infty} \frac{(1-e^{n\phi_{1}t})(1-e^{n\phi_{2}t})}{(1-e^{-n\e_{1}t})(1-e^{-n\e_{2}t})} \frac{q^{n}}{n} \right)\\
&=& 
\exp \left(  \frac{\phi_{1}\phi_{2}}{\e_{1}\e_{2}} \sum_{n=1}^{\infty} \frac{q^{n}}{n} \right)\\
&=&
(1-q)^{-\frac{\phi_{1}\phi_{2}}{\e_{1}\e_{2}}}.
\end{eqnarray*}
\endproof

We consider the Nekrasov partition function 
$$
Z(\mbi \e, \mbi a, \mbi m)=\sum_{n=0}^{\infty} q^{n} \int_{M(1,n)} e(\mc F_{1}(\mc V))
$$ 
for rank $r=1$.
By the explicit form in \eqref{comb}, we can substitute $\phi_{1}=a_{1}-\frac{\e_{+}}{2}+m_{1}, \phi_{2}= -(a_{1} - \frac{\e_{+}}{2} + m_{2})$, and replacing $q$ with $-q$ in Proposition \ref{hilb2} to get
$$
Z(\mbi \e, \mbi a, \mbi m)=(1+q)^{\alpha_{1}},
$$  
where $\alpha_{1}=\frac{(a-\frac{\e_{+}}{2}+m_{1})(a-\frac{\e_{+}}{2}+m_{2})}{\e_{1}\e_{2}}$.

%%%%%%%%%%%%%%%%%%%%%%%%%%%%%%%%%%%%%%%%%%%%%%%%%%%%%%%%%%%%%%%%%%%%%%%%%%%
%%%%%%%%%%%%%%%%%%%%%%%%%%%%%%%%%%%%%%%%%%%%%%%%%%%%%%%%%%%%%%%%%%%%%%%%%%%
%%%%%%%%%%%%%%%%%%%%%%%%%%%%%%%%%%%%%%%%%%%%%%%%%%%%%%%%%%%%%%%%%%%%%%%%%%%

\section{Wall-crossing formulas}
\label{sec:wall}
%%%%%%%%%%%%%%%%%%%%%%%%%%%%%%%%%%%%%%%%%%%%%%%%%%%%%%%%%%%%%%%%%%%%%%%%%%%
In this section we derive wall-crossing formulas following \cite[Chapter 7]{M} and \cite[\S 6]{NY3}, and give a proof of Theorem \ref{main}.

\subsection{Localization}
\label{subsec:loca2}
Let us consider the Euler class $\psi \in A^{\ast}_{\tilde{T} \times \C^{\ast}_{\hbar}}(M(r,n))$ of $\tilde{T} \times \C^{\ast}_{\hbar}$-equivariant $K$-theory class defined by a linear combination of tensor products of $\mc V, \mc V^{\vee}$,  and $\tilde{T}$-modules.
For exmaple, we take $\psi=1$, $\psi=e(\mc F_{r}(\mc V))$, or $\psi=e(TM(r,n) \otimes e^{m_{1}})$. 
We also write by $\psi$ the class defined by the same formula on $M^{c}(r,n)$, $\wt{M}^{\ell}(r,n)$, $\mc M$, and so on.

Let us consider the pull-back $\Theta^{rel}$ to $\wt{M}^{\ell}(r,n)$ of the relative tangent bundle of $[\tilde{\mu}^{-1}(0)/G]$ over $[\mu^{-1}(0)/G]$.
We also write by $\Theta^{rel}$ the pull-back to the enhanced master space $\mc M$ as in the previous section.
We introduce $\tilde{\psi}= \frac{1}{n!}\psi \cup e(\Theta^{rel})$ on $\wt{M}^{\ell}(r,n)$ and $\mc M$ so that 
$$
\int_{\wt{M}^{n}(r,n)} \tilde{\psi} = \int_{M^{c}(r,n)} \psi, \int_{\wt{M}^{0}(r,n)} \tilde{\psi} = \int_{M(r,n)} \psi.
$$
We consider integrations $\int_{\mc M} \tilde{\psi}$ over enhanced master spaces.
This integration is defined by 
$$
\int_{\mc M} \tilde{\psi} =(\iota_{0})_{\ast}^{-1} \Pi_{\ast} \left( \tilde{\psi} \cap [\mc M]^{vir} \right) \in \mc S[\hbar],
$$ 
where $\Pi \colon \mc M \to M_{0}(r,n)$ as in \S \ref{subsec:modu}, $\iota_{0} \colon \lbrace n[0] \rbrace \to M_{0}(r,n)$ as in \S \ref{subsec:part}, and $[\mc M]^{vir}=[\mc M]$ is the (virtual) fundamental cycle defined by the obstruction theory $ob_{\mc M}=\id_{L_{\mc M}}$.
We also use similar push-forward homomorphisms from Chow groups of various moduli stacks to define integrals, for example $\int_{\mc M_{\pm}}, \int_{\mc M_{\mk J}}$.

By \eqref{local}, we have the following commutative diagram:
$$
\xymatrix{
\varprojlim_{m} A^{\ast}_{\C^{\ast}_{\hbar}} (\mc M \times_{\tilde{T}} E_{m}) \otimes_{\C[\hbar]} \C[\hbar^{\pm 1}] \ar[d]_{\int_{\mc M}} \ar[r]^{\cong}& \varprojlim_{m} A^{\ast}_{\C^{\ast}_{\hbar}} (\mc M^{\C^{\ast}_{\hbar}} \times_{\tilde{T}} E_{m}) \otimes_{\C[\hbar]} \C[\hbar^{\pm 1}] \ar[d]_{\int_{\mc M_{+}} + \int_{\mc M_{-}} + \sum_{\mk J} \int_{\mc M_{\mk J}}}\\
\varprojlim_{m} A_{\ast} (M_{0}(r,n) \times_{\tilde{T}} E_{m})[\hbar, \hbar^{-1}] \ar@{=}[r] & \varprojlim_{m} A_{\ast} (M_{0}(r,n) \times_{\tilde{T}} E_{m})[\hbar,\hbar^{-1}]
}
$$
where the upper horizontal arrow is given by 
$$
\frac{\iota_{+}^{\ast}}{e(\mk N(\mc M_{+}))} + \frac{\iota_{-}^{\ast}}{e(\mk N(\mc M_{-}))} + \sum_{\mk J \in S^{\ell} }\frac{\iota_{\mk J}^{\ast}}{e(\mk N(\mc M_{\mk J}))}.
$$
Hence we have
\begin{eqnarray*}
\int_{\mc M} \tilde{\psi} 
& = &
\int_{\mc M_{+}} \frac{ \tilde{\psi} }{e(\mk N(\mc M_{+}))}  + \int_{\mc M_{-}} \frac{ \tilde{\psi} }{ e(\mk N(\mc M_{-})) } + 
\sum_{\mk J\in S^{\ell}}\int_{\mc M_{\mk J}} \frac{ \tilde{\psi} }{e(\mk N(\mc M_{\mk J}))}.
\end{eqnarray*}

The left hand side is a polynomial in $\hbar$.
Hence taking the coefficient in $\hbar^{-1}$, we have 
\begin{eqnarray}
\label{wc1}
\int_{\wt M^{\ell}(r,n)} \tilde{\psi}  - \int_{M(r,n)} \psi = \sum_{\mk J\in S^{\ell}} \Res_{\hbar=\infty} \int_{\mc M_{\mk J}} \frac{ \iota_{\mk J}^{\ast}\tilde{\psi} }{e(\mk N(\mc M_{\mk J}))},
\end{eqnarray}
by Lemma \ref{ob1} and $\mc M_{+}=\wt M^{\ell}(r,n) $, $\mc M_{-} = \wt{M}^{0}(r,n)$.
Here $\Res_{\hbar=\infty}$ means taking the minus of the coefficient in $\hbar^{-1}$, since $\hbar^{-1}d\hbar = - w^{-1} dw$ for $w=\hbar^{-1}$.

By Theorem \ref{decomp}, Proposition \ref{ob2}, and \eqref{rel}, we have
\begin{eqnarray}
\label{wc2}
&&
\int_{\mc M_{\mk J}} \frac{ \iota_{\mk J}^{\ast} \tilde{\psi} }{e \left( \mk N(\mc M_{\mk J}) \right) } \notag \\
&=&
\frac{1}{n!}\int_{\wt M^{\min(I_{\sharp}) -1}(r,n-p)} e(\Theta^{rel}_{\flat}) \notag \\
&\cup& 
\int_{[M^{+}_{p}]^{vir}} \frac{\psi \left( p_{1}^{\ast} \mc V_{\flat} \oplus p_{2}^{\ast} \mc V_{p}^{+} \otimes e^{\frac{\hbar}{pD}} \otimes (\det \mc V_{\flat})^{\frac{1}{p}} \right) \cup e(p_{1}^{\ast} \Theta^{rel}_{\sharp})}{e\left( \mk N_{1} \otimes e^{\frac{\hbar}{pD}} \otimes (\det \mc V_{\flat})^{\frac{1}{p}} \oplus \mk N_{2} \otimes e^{\frac{-\hbar}{pD}} \otimes (\det \mc V_{\flat})^{\frac{-1}{p}} \right) }.
\end{eqnarray}
Here $\int_{[M^{+}_{p}]^{vir}} (\cdot)$ is the push-forward  
$$
(p_{1})_{\ast} \left( (\cdot) \cap [\wt{M}^{\min(I_{\sharp})-1}(r,n-p) \times M_{p}^{+}]^{vir} \right)
$$ 
by the projection $p_{1} \colon \wt{M}^{\min(I_{\sharp})-1}(r,n-p) \times M_{p}^{+} \to \wt{M}^{\min(I_{\sharp})-1}(r,n-p)$, and 
$$
[\wt{M}^{\min(I_{\sharp})-1}(r,n-p) \times M_{p}^{+}]^{vir}
$$ 
is the virtual fundamental cycle defined by the obstruction theory 
$$
p_{1}^{\ast}ob_{\wt{M}^{\min(I_{\sharp})-1}(r,n-p)} \oplus p_{2}^{\ast}ob_{M^{p}_{+}}
$$ 
as in Proposition \ref{ob2}.
\begin{NB}
We consider $\tilde{T}$-equivariant morphisms $\pi_{i} \colon X \to Y_{i}$ for $i=1,2$, such that $Y_{i}^{T}=\lbrace p_{i} \rbrace$ where $p_{i}$ are points in $Y_{i}$.
Then two integrations $(\iota_{i \ast})^{-1} \circ \pi_{i \ast}$ are same, since by localization theorem they are both equal to the integration over $\tilde{T}$-fixed locus on $X$ weighted by the Euler class of the normal bundle.
\end{NB}

In the following, we write by $\mc V_{\flat}, \mc V_{\sharp}$ the pull-backs $p_{1} \mc V_{\flat}, p_{2}^{\ast} \mc V_{\sharp}$.
%%%%%%%%%%%%%%%%%%%%%%%%%%%%%%%%%%%%%%%%%%%%%%%%%%%%%%%%%%%%%%%%%%%%%%%%%%%
\subsection{Computations of residues}
\label{subsec:comp}
We simplify integrations as in \cite[Proof of Theorem 7.2.4]{M} and \cite[\S 6.3]{NY3}.
The integrant of $\int _{[M^{+}_{p}]^{vir}}$ in the above equation \eqref{wc2} is of the form $f(\hbar) = \sum_{j=-\infty}^{m} A_{j} \hbar^{j}$ for some $m \in \Z$, where 
$$
A_{j} \in A^{\ast}_{\tilde{T}} \left( \wt{M}^{\min(I_{\sharp})-1}(r,n-p) \times M_{p}^{+} \right)
$$ 
does not depend on $\hbar$.
Since we have 
$$
\Res_{\hbar =\infty} f(\hbar) = pD\Res_{\hbar_{s}=\infty} f \left( pD \hbar_{s} -  D c_{1}(\mc V_{\flat}) \right),
$$
the residue of the integral is equal to 
\begin{eqnarray*}
&&
pD\Res_{\hbar=\infty} \int_{[M^{+}_{p}]^{vir}} \frac{\psi (\mc V_{\flat} \oplus \mc V_{p}^{+} \otimes e^{\hbar}) \cup e(\Theta^{rel}_{\sharp})}{e\left( \mk N_{1} \otimes e^{\hbar} \oplus \mk N_{2} \otimes e^{-\hbar} \right) } \\
&=&
(p-1)!\Res_{\hbar=\infty} \int_{[M(1,p)]^{vir}} \frac{\psi (\mc V_{\flat} \oplus \mc V_{\sharp} \otimes e^{\hbar} ) \cup e(\mc V_{\sharp} / \mo_{M(1,p)})}{e\left( \mk N_{1} \otimes e^{\hbar} \oplus \mk N_{2} \otimes e^{-\hbar}\right)}.
\end{eqnarray*}
%Here $\mc V_{\flat}, \mc V_{\sharp}$ are tautological bundles on $M^{\min (I_{\sharp})-1}(r,n-p), M(1,p)$, and $\mk N_{1}, \mk N_{2}$ are defined on $M^{\min (I_{\sharp})-1}(r,n-p) \times M(1,p)$ by the similar formulas as in \eqref{n1}, \eqref{n2}.
The last equality follows from Proposition \ref{hilb}, \eqref{vp+} and \eqref{ob4}.
We note that we can ignore $L_{M_{p}^{+}}$ in \eqref{vp+} by the similar computations of residue as above.
By localization theorem this is
\begin{eqnarray*}
(p-1)!\Res_{\hbar=\infty} \sum_{|Y|=p} \frac{\iota_{Y^{\flat}}^{\ast} \psi (\mc V_{\flat} \oplus  \mc V_{\sharp} \otimes e^{\hbar})}{\iota_{Y^{\flat}}^{\ast} e\left( \mk N_{1} \otimes e^{\hbar} \oplus \mk N_{2} \otimes e^{-\hbar} \right)}  \cup \frac{\iota_{Y}^{\ast}e(\mc V_{\sharp} / \mo_{M(1,p)})}{ e\left( \mk{N}^{\sharp}_{Y} \right) },
\end{eqnarray*}
where $\iota_{Y^{\flat}} \colon\wt{M}^{\min(I_{\sharp})-1}(r,n-p) \times \lbrace I_{Y} \rbrace \to \wt{M}^{\min(I_{\sharp})-1}(r,n-p) \times M(1,p)$ is the inclusion.
Here the sum is taken over the set of Young diagrams $Y$ with the weight $|Y|=p$, and $\mk{N}^{\sharp}_{Y}$ is the virtual normal bundle of $\lbrace I_{Y} \rbrace$ in $M(1,p)$ defined by the obstruction theory $ob^{\sharp}_{M(1,p)}$ in \S \ref{subsec:obst4}. 

\begin{NB}
This computation is verified in the following way.
We consider torus $T$-equivariant morphisms $M \to M_{0}$ and $N \to N_{0}$ of Deligne-Mumford stacks.
Then we have the following diagrams
$$
\xymatrix{
A^{T}_{\ast}(M \times N) \ar[d] &\ar[l] A^{T}_{\ast}(M \times N^{T}) \ar[d] & \ar[l] A^{T}_{\ast}(M^{T} \times N^{T}) \ar[d] \\
A^{T}_{\ast}(M_{0} \times N_{0}) &\ar[l] A^{T}_{\ast}(M_{0} \times N_{0}^{T})& \ar[l] A^{T}_{\ast}(M_{0}^{T} \times N_{0}^{T})}
$$
$$
\xymatrix{
A^{T}_{\ast}(M) \otimes A^{T}_{\ast}(N) \ar[d] &\ar[l] A^{T}_{\ast}(M) \otimes A^{T}_{\ast}(N^{T}) \ar[d] & \ar[l] A^{T}_{\ast}(M^{T}) \otimes A^{T}_{\ast}(N^{T}) \ar[d] \\
A^{T}_{\ast}(M_{0}) \otimes A^{T}_{\ast}(N_{0}) &\ar[l] A^{T}_{\ast}(M_{0}) \otimes A^{T}_{\ast}(N_{0}^{T}) & \ar[l] A^{T}_{\ast}(M_{0}^{T}) \otimes A^{T}_{\ast}(N_{0}^{T}),}
$$
where all horizontal maps are isomorphisms after localizations.
We have morphisms from the latter to the former diagram such that every square commutes.

Furthermore we have the following diagram
$$
\xymatrix{
A^{T}_{\ast}(M \times N) \ar[d]_{p_{\ast}} &\ar[l] A^{T}_{\ast}(M) \otimes A^{T}_{\ast}(N^{T}) \ar[d] \\
A^{T}_{\ast}(M) & \ar[l]_{\cong} A^{T}_{\ast}(M) \otimes A^{T}_{\ast}(N_{0}^{T}) ,}
$$
where $p \colon M \times N \to M$ is the projection.
By applying these diagrams for $M=\wt{M}^{\min(I_{\sharp}) -1}(r,n-p), N = M(1,p)$, we can verify the above computations.
\end{NB}

In the following we consider the case where $\psi(\mc V)=e(\mc F_{r}(\mc V))$.
Then we have 
$$
\psi(\mc V_{\flat} \oplus \mc V_{\sharp} \otimes e^{\hbar}) = e(\mc \F_{r}(\mc V_{\flat})) \cup e(\mc F_{r}(\mc V_{\sharp} \otimes e^{\hbar})).
$$
Since the degree of $\frac{\iota_{Y^{\flat}}^{\ast}e \left( \mc F_{r}(\mc V_{\sharp} \otimes e^{\hbar}) \right)}{ \iota_{Y^{\flat}}^{\ast} e\left( \mk N_{1} \otimes e^{\hbar} \oplus \mk N_{2} \otimes e^{-\hbar} \right) }$ with respect to $\hbar$ is equal to $0$, by the definitions \eqref{fr}, \eqref{n1} and \eqref{n2} of $\mc F_{r}(\mc V_{\sharp} \otimes e^{\hbar})$, $\mk N_{1}$ and $\mk N_{2}$, we have
\begin{eqnarray*}
&&
\Res_{\hbar=\infty} \frac{ \iota_{Y^{\flat}}^{\ast}e \left( \mc F_{r}(\mc V_{\sharp} \otimes e^{\hbar}) \right)}{\iota_{Y^{\flat}}^{\ast} e\left( \mk N_{1} \otimes e^{\hbar} \oplus \mk N_{2} \otimes e^{-\hbar} \right)}\\
&=&
\Res_{\hbar=\infty} \frac{\iota_{Y^{\flat}}^{\ast}e \left( \mc F_{r}(\mc V_{\sharp} \otimes e^{\hbar}) \right)}{\iota_{Y^{\flat}}^{\ast} e \left( \mc Hom(\mc W, \mc V_{\sharp}) \otimes e^{\hbar} \oplus \mc Hom(\wedge Q^{\vee} \otimes \mc V_{\sharp}, \mc W) \otimes e^{-\hbar} \right)}\\
&=&(-1)^{rp+1} p \left( 2 \sum_{\alpha=1} ^{r} a_{\alpha} + \sum_{f=1} ^{2r} m_{f} \right).
\end{eqnarray*}
The last equality follows from the similar combinatorial desctiptions as in \S \ref{subsec:appl}.
But for computations of $\iota_{Y^{\flat}}^{\ast}e \left( TM(1,p) \right)$, we must delete $a_{\alpha}$, since we see that $M(1,p)$ does not have $\GL(W)$-action here.
This is because $M(1,p)$ parametrizes ADHM data on $(W_{\sharp}, V_{\sharp})$ in the construction of $h \colon M_{p}^{+} \to M(1,p)$ in \S \ref{subsec:morp}.

\begin{NB}
We have
\begin{eqnarray*}
\iota_{Y}^{\ast}e\left( \mc F_{r} (\mc V_{\sharp} \otimes e^{\hbar} ) \right)
=
\prod_{(i,j) \in Y} 
\prod_{f=1}^{2r}
\left( \hbar - \left( i-\frac{1}{2} \right) \e_{1} - \left(j-\frac{1}{2} \right) \e_{2} + m_{f} \right).
\end{eqnarray*}
The numerator of $\iota_{Y^{\flat}}^{\ast} e\left( \mk N_{1} \otimes e^{\hbar} \oplus \mk N_{2} \otimes e^{-\hbar} \right)^{-1}$ is equal to 
\begin{multline*}
\prod_{(i,j)\in Y} 
\left( \hbar - (i-2) \e_{1} - (j-2) \e_{2} - \mc V_{\flat} \right) \times
\left( \hbar -i \e_{1} - j \e_{2} - \mc V_{\flat} \right) \times \\
\left( \hbar - (i-1) \e_{1} - (j-1) \e_{2} - \mc V_{\flat} \right) \times 
\left (\hbar - (i - 1) \e_{1} - (j - 1) \e_{2} - \mc V_{\flat} \right).
\end{multline*}
The denominator of $\iota_{Y^{\flat}}^{\ast} e\left( \mk N_{1} \otimes e^{\hbar} \oplus \mk N_{2} \otimes e^{-\hbar} \right)^{-1}$ is equal to 
\begin{multline*}
(-1)^{rp} \prod_{(i,j) \in Y} 
\left( \prod_{\alpha=1}^{r} 
( \hbar -(i-1) \e_{1} - (j-1) \e_{2} - a_{\alpha} ) \times
(\hbar -i \e_{1} -j \e_{2} - a_{\alpha} ) \right) \times \\
(\hbar - (i-2) \e_{1} - (j-1) \e_{2} - \mc V_{\flat} ) \times
( \hbar - (i-1) \e_{1} - (j-2) \e_{2} - \mc V_{\flat} ) \times \\
( \hbar -i \e_{1} -(j-1) \e_{2} - \mc V_{\flat} ) \times 
( \hbar - (i - 1) \e_{1} - j \e_{2} - \mc V_{\flat}). 
\end{multline*}
\end{NB}

Thus we have 
\begin{eqnarray}
\label{exc}
\Res_{\hbar=\infty} \int_{\mc M_{\mk J}} \frac{\iota_{\mk J}^{\ast} \tilde{\psi}}{\mk N(\mc M_{\mk J})} 
&=&
\frac{(-1)^{rp+1} p!}{n !}  \left( 2 \sum_{\alpha=1} ^{r} a_{k} + \sum_{f=1} ^{2r} m_{f} \right) \notag \\
&& 
\sum_{|Y|=p} \frac{\iota_{Y}^{\ast}e(\mc V_{\sharp} / \mo_{M(1,p)})}{\iota_{Y}^{\ast} e\left( \mk{N}^{\sharp}_{Y} \right)}   \notag \\
&&
\int_{\wt M^{\min(I_{\sharp}) -1}(r,n-p)} e(\mc F_{r}(\mc V_{\flat}) \oplus\Theta^{rel}_{\flat}).
\end{eqnarray}
Combining this equation with the following proposition and \eqref{wc1}, we have
\begin{multline}
\label{wc3}
\frac{1}{n!}\int_{\wt M^{\ell}(r,n)} e(\mc F_{r}(\mc V) \oplus \Theta^{rel})  - \int_{M(r,n)} e(\mc F_{r}(\mc V)) =\\
(-1)^{rp+1}  \frac{(p-1)!}{n!} u_{r}\sum_{\mk J \in S^{\ell}} \int_{\wt M^{\min(I_{\sharp}) -1}(r,n-p)} e(\mc F_{r}(\mc V_{\flat}) \oplus \Theta^{rel}_{\flat}), 
\end{multline}
where $u_{r}=\frac{\e_{+}\left( 2 \sum_{\alpha=1} ^{r} a_{k} + \sum_{f=1} ^{2r} m_{f} \right)}{\e_{1}\e_{2}}$ as in Theorem \ref{main}.

\begin{prop}
\label{hilb1}
We have
$$
\sum_{|Y|=p} \frac{\iota_{Y}^{\ast}e(\mc V_{\sharp} / \mo_{M(1,p)})}{ e\left( \mk{N}^{\sharp}_{Y} \right)}= \frac{\e_{+}}{p\e_{1}\e_{2}}.
$$
\end{prop}
\proof
By \eqref{sharp}, we have 
$$
e\left( \mk{N}^{\sharp}_{Y} \right) = \frac{ \iota^{\ast}_{Y} e(TM(1,p))}{ \iota^{\ast}_{Y}e(\mc Hom(\wedge Q^{\vee} \otimes \mc V_{\sharp}, \mc W_{\sharp}))},
$$ 
where $\mc W_{\sharp}=\mo_{M(1,p)}$ is the tautological bundle corresponding to $W_{\sharp}$.
Since $\wedge Q = t_{1}t_{2}$ as a $\tilde{T}$-module, we have
\begin{eqnarray*}
\sum_{ |Y|=p} \frac{\iota_{Y}^{\ast}e(\mc V_{\sharp} / \mo_{M(1,p)})}{ e\left( \mk{N}^{\sharp}_{Y} \right) }
&=&
\sum_{|Y|=p} \frac{\iota_{Y}^{\ast}e(\mc V_{\sharp}/\mo_{M(1,p)}) \iota_{Y}^{\ast}e(\mc V_{\sharp}^{\vee} \otimes t_{1}t_{2})}{\iota_{Y}^{\ast}e(TM(1,p))}\\
&=&
\int_{M(1,p)} e\left( (\mc V_{\sharp}/\mo_{M(1,p)}) \oplus \mc V_{\sharp}^{\vee} \otimes t_{1}t_{2} \right).
\end{eqnarray*}
Using new variables $\phi_{1}, \phi_{2}$ and corresponding $T$-modules $e^{\phi_{1}}, e^{\phi_{2}}$, this is equal to
\begin{eqnarray*}
\left( \frac{1}{\phi_{1}}\int_{M(1,p)} e(\mc V_{\sharp}\otimes e^{\phi_{1}} \oplus \mc V_{\sharp}^{\vee} \otimes e^{\phi_{2}} ) \right) \Big|_{\substack{\phi_{1}=0 \\ \phi_{2}=\e_{+}}}
&=&
\frac{\e_{+}}{p\e_{1}\e_{2}}.
\end{eqnarray*}
The last equality follows from Proposition \ref{hilb2}.
\endproof

%%%%%%%%%%%%%%%%%%%%%%%%%%%%%%%%%%%%%%%%%%%%%%%%%%%%%%%%%%%%%%%%%%%%%%%%%%%

\subsection{Proof of Theorem \ref{main}}
\label{subsec:proo}
Here we complete wall-crossing formula by the equation \eqref{wc3}, and prove Theorem \ref{main}.
For $0 \le q \le n$, we put 
$$
\gamma_{q,n}=\frac{1}{n!}\int_{\wt M^{q}(r,n)} e(\mc F_{r}(\mc V) \oplus \Theta^{rel}).
$$
We have $\gamma_{0,n}=\alpha_{n}, \gamma_{n,n}=\beta_{n}$.
%\begin{eqnarray}
%\label{gamma}
%\gamma_{\ell, n}-\gamma_{0,n}
%&=&
%\sum_{q=0}^{\ell - 1}\sum_{p=1}^{n-q} (-1)^{rp+1} \frac{(n-q-1)!(n-p)!}{(n-p-q)!n!} u_{r} \gamma_{q,n-p}\notag\\
%&=&
%\frac{1}{n}\sum_{q=0}^{\ell - 1}\sum_{p=1}^{n-q} (-1)^{rp+1} \frac{(n-q-1)(n-q-2) \cdots (n-q-p+1)}{(n-1)(n-2)\cdots (n-p+1)} u_{r} \gamma_{q,n-p}.
%\end{eqnarray}

We consider the set $\text{Dec}^{i}(n)$ of $i$ tuple $\mk J^{i}=(I_{\sharp}^{1}, \ldots, I_{\sharp}^{i})$ of non-empty disjoint subsets of $\underbar{n}$ such that  $\min(I_{\sharp}^{1}) > \cdots > \min(I_{\sharp}^{i})$.
We put $|\mk J^{i}|=|I_{\sharp}^{1}| \cdots + |I_{\sharp}^{i}|$ for each $\mk J^{i} \in \text{Dec}^{i}(n)$.
\begin{lem}[\protect{\cite[Lemma 7.6.5]{M}, \cite[Lemma 6.6]{NY3}}]
We have
\begin{multline*}
\beta_{n} - \alpha_{n}=\sum_{1\le i < j}\sum_{\mk J^{i}\in \text{Dec}^{i}(n)} (-1)^{r|\mk J^{i}| + i}
\frac{(|I_{\sharp}^{1}|-1)! \cdots (|I_{\sharp}^{i}|-1)!(n-|\mk J^{i}|)!}{n!} u_{r}^{i} \alpha_{n-|\mk J^{i}|}\\
+\sum_{\mk J^{j}\in \text{Dec}^{j}(n)} (-1)^{r|\mk J^{j}| + j}
\frac{(|I_{\sharp}^{1}|-1)! \cdots (|I_{\sharp}^{j}|-1)!(n-|\mk J^{j}|)!}{n!} u_{r}^{j} \gamma_{\min(I_{\sharp}^{j})-1,n-|\mk J^{j}|}.\\
\end{multline*}
\end{lem}
\proof
We use induction on $j>0$. 
When $j=1$, this equation is nothing but \eqref{wc3} for $\ell=n$, since we have $\text{Dec}^{1}(n) \cong S^{n}$.
Then applying \eqref{wc3} repeatedly we get the assertion for any $j>0$.
\endproof
Since $\text{Dec}^{j}(n)=\emptyset$ for $j>n$, by the above lemma we have
$$
\beta_{n} - \alpha_{n} = \sum_{1 \le i \le n}\sum_{\mk J^{i}\in \text{Dec}^{i}(n)} (-1)^{r|\mk J^{i}| + i}
\frac{(|I_{\sharp}^{1}|-1)! \cdots (|I_{\sharp}^{i}|-1)! (n-|\mk J^{i}|)!}{n!} u_{r}^{i} \alpha_{n-|\mk J^{i}|}.
$$
We consider a map $\rho_{i} \colon \text{Dec}^{i}(n) \to \mb N^{i}$ sending $\mk J^{i}$ to $\left( |I_{\sharp}^{1}|, \ldots, |I_{\sharp}^{i}| \right).$
By \cite[Lemma 7.6.6]{M} or \cite[Lemma 6.8]{NY3}, for any $\vec{p}_{(i)}=(p_{1}, \ldots, p_{i}) \in \mb N^{i}$ such that $|\vec{p}_{(i)}|=p_{1} + \cdots + p_{i} \le n$ we have
$$
|\rho_{i}^{-1}(\vec{p}_{(i)})| \frac{(p_{1}-1)! \cdots (p_{i}-1)! (n-|\vec{p}_{(i)}|)!}{n!} = \frac{1}{\prod_{j=1}^{i} \sum_{1 \le h \le j} p_{h}}, 
$$
and hence 
$$
\beta_{n} - \alpha_{n} = \sum_{1 \le i \le n}\sum_{|\vec{p}_{(i)}| \le n} 
(-1)^{r|\vec{p}_{(i)}| + i} \frac{u_{r}^{i}}{\prod_{j=1}^{i} \sum_{1 \le h \le j} p_{h}}  \alpha_{n-|\vec{p}_{(i)}|}.
$$

To prove Theorem \ref{main}, for each $k=1,\ldots, n$, we must show
\begin{multline*}
\label{goal}
\sum_{i=1}^{k} \sum_{p_{1}+\cdots +p_{i}=k}(-1)^{i} \frac{u_{r}^{i}}{\prod_{j=1}^{i} \sum_{1 \le h \le j} p_{h}} 
%\\=
%(-1)^{k} \frac{u_{r}(u_{r}-1) \cdots (u_{r}-k+1)}{k!}.
%\end{multline*}
%The right hand side of \eqref{goal} is equal to 
%\begin{multline*}
=\sum_{i=1}^{k} \sum_{r_{1}<\cdots<r_{i}=k}(-1)^{i} \frac{u_{r}^{i}}{\prod_{j=1}^{i} r_{j}}.
\end{multline*}
This follows from the bijection
$$
(p_{1}, \ldots, p_{i}) \mapsto (r_{1}, \ldots, r_{i})=(p_{1}, p_{1}+p_{2}, \ldots, p_{1}+\cdots+p_{i}).
$$
%%%%%%%%%%%%%%%%%%%%%%%%%%%%%%%%%%%%%%%%%%%%%%%%%%%%%%%%%%%%%%%%%%%%%%%%%%%

\subsection{Partition functions defined from other classes}
\label{subsec:part}
Finally we consider Nekrasov partition functions 
$$
Z^{\psi}(\mbi \e, \mbi a, \mbi m^{N_{f}}, q)=\sum_{n=0}^{\infty} q^{n} \int_{M(r,n)} \psi
$$ 
defined from $\psi \in A^{\ast}_{\tilde{T}}(M(r,n))$ other than $e(\mc F_{r}(\mc V))$.

Here we consider the case where $N_{f}=0$ and $\psi=1$, or $N_{f}=1$ and 
$$
\psi= e\left( TM(r,n) \otimes e^{\sqrt{t_{1}t_{2}}m_{1}} \right).
$$
Then in both cases, residues in \eqref{wc1} with respect to $\hbar$ are equal to zero. 
\begin{NB}
\begin{multline*}
\int_{\mc M_{\mk J}} \frac{e(\mc F_{r}(\mc V) \oplus \Theta^{rel})}{e \left( \mk N(\mc M_{\mk J}) \right) } \\
=
\int_{\wt M^{\min(I_{\sharp}) -1}(r,n-p)} e(\mc F_{r}(\mc V_{\flat}) \oplus \Theta^{rel}_{\flat})\\
 \times \int_{[M^{+}_{p}]^{vir}} e\left( \mk N(\mc V'_{\sharp}, \mc V_{\flat}, \mc W) \otimes e^{\sqrt{t_{1}t_{2}}m_{1}}\right)  \cup e\left( \mk N(\mc W, \mc V_{\flat}, \mc V_{\sharp}') \right) \cup e(\Theta^{rel}_{\sharp}),
\end{multline*}
\end{NB}
Hence the right hand side of \eqref{wc1} is equal to zero.
\begin{NB}
we have
\begin{eqnarray*}
\frac{1}{n!}\int_{\wt M^{\ell}(r,n)} e(TM(r,n) \otimes e^{m_{1}} \oplus \Theta^{rel})  - \int_{M(r,n)} e(TM(r,n) \otimes e^{m_{1}}) =0
\end{eqnarray*}
for any $\ell =0, \ldots, n$.
\end{NB}
As a result we have 
$$
Z^{\psi}(-\mbi \e, \mbi a, \mbi m^{N_{f}}, q)=Z^{\psi}(\mbi \e, \mbi a, \mbi m^{N_{f}}, q).
$$
When $N_{f}=0$ and $\psi =1$, this formula also implies \eqref{odd}.
When $N_{f}=r=1$ and $\psi=e\left( TM(1,n)\otimes e^{\sqrt{t_{1}t_{2}}m_{1}} \right)$, by \cite[Corollary 1]{CO} we have 
$$
Z^{\psi}(\mbi \e, \mbi a, m_{1} ,q)= \prod_{n}(1-q^{n})^{\frac{\e_{+}^{2}}{4}-m_{1}^{2}-1}
$$ 
and we can also check this formula.

%\begin{acknowledgements}
%If you'd like to thank anyone, place your comments here
%and remove the percent signs.
%\end{acknowledgements}

% BibTeX users please use one of
%\bibliographystyle{spbasic}      % basic style, author-year citations
%\bibliographystyle{spmpsci}      % mathematics and physical sciences
%\bibliographystyle{spphys}       % APS-like style for physics
%\bibliography{}   % name your BibTeX data base

% Non-BibTeX users please use

\noindent
              R. Ohkawa \\
              Department of Mathematics, School of Fundamental Science and Engineering, Waseda University, 3--4--1 Okubo, Shinjuku-ku, Tokyo 169--8555, Japan\\
              Tel.: +81-3-5286-3195\\
              ohkawa.ryo@aoni.waseda.jp           
              
\end{document}